\documentclass[preprint,12pt]{elsarticle}

\usepackage{amssymb}
\usepackage{amsmath}
\usepackage{graphicx}      
\usepackage{multirow}      
\usepackage{booktabs}      
\usepackage{algorithm}
\usepackage{algpseudocode}
\usepackage{subcaption}
\usepackage{xcolor}
\usepackage{bm}
\journal{}

\begin{document}

\begin{frontmatter}



\title{A coupled finite element-virtual element method for thermomechanical analysis of electronic packaging structures}


\author[address1]{Yanpeng Gong}
\ead{yanpeng.gong@bjut.edu.cn}
\author[address1]{Sishuai Li}
\author[address3]{Yue Mei\corref{mycorrespondingauthor}}
\ead{meiyue@dlut.edu.cn}
\author[address2]{Bingbing Xu}
\author[address1]{Fei Qin}
\author[address4,address5]{Xiaoying Zhuang}
\author[address6]{Timon Rabczuk}

\address[address1]{Department of Mechanics, Beijing University of Technology, Beijing, 100124, China}
\address[address2]{Institute of Continuum Mechanics, Leibniz University Hannover, Hannover, Germany}
\address[address3]{State Key Laboratory of Structural Analysis for Industrial Equipment, Department of Engineering Mechanics, International Research Center for Computational Mechanics, Dalian University of Technology,Dalian 116023, China}
\address[address4]{Chair of Computational Science and Simulation Technology, Institute of Photonics, Department of Mathematics and Physics, Leibniz University Hannover, 30167 Hannover, Germany}
\address[address5]{Department of Geotechnical Engineering, College of Civil Engineering, Tongji University, Shanghai, 200092, China}
\address[address6]{Institute of Structural Mechanics, Bauhaus University Weimar, 99423 Weimar, Germany}

\cortext[mycorrespondingauthor]{Corresponding author}


\begin{abstract}
    This study presents a finite element and virtual element (FE-VE) coupled method for thermomechanical analysis in electronic packaging structures. The approach partitions computational domains strategically, employing FEM for regular geometries to maximize computational efficiency and VEM for complex shapes to enhance geometric flexibility. Interface compatibility is maintained through coincident nodal correspondence, ensuring solution continuity across domain boundaries while reducing meshing complexity and computational overhead. Validation through electronic packaging applications demonstrates reasonable agreement with reference solutions and acceptable convergence characteristics across varying mesh densities. The method effectively captures thermal distributions and stress concentrations in multi-material systems, establishing a practical computational framework for electronic packaging analysis involving complex geometries. Source codes are available at https://github.com/yanpeng-gong/FeVeCoupled-ElectronicPackaging.
\end{abstract}




\end{frontmatter}


\section{Introduction}
\label{sec:introduction}

Modern electronic devices face increasing thermal management challenges due to rising chip integration density and power density, which generate elevated operating temperatures that degrade performance and reduce service life~\cite{xiao2017effective}. The multi-material nature of electronic packaging exacerbates these challenges, as thermal expansion coefficient mismatches between heterogeneous materials create substantial thermal stresses during temperature cycling~\cite{li2024efficient}. These thermomechanical stresses contribute to critical reliability failures, including solder joint cracking~\cite{gong2024application}, wire bond fatigue, and substrate delamination~\cite{su2014investigation}. Recent studies indicate that thermomechanical effects account for over 55\% of electronic device failures in advanced 3D integrated circuits~\cite{he2021thermal}, underscoring the critical importance of accurate thermal and thermomechanical analysis for electronic packaging reliability assessment.

Numerical simulation methods have become essential for electronic packaging reliability analysis, offering cost-effective and reliable assessment capabilities. However, thermomechanical coupling analysis in electronic packaging faces significant computational challenges due to multi-scale geometric features and complex material interfaces in modern designs. Researchers have developed various numerical approaches to address these challenges, including the Finite Element Method (FEM)~\cite{ye2022practical}, Boundary Element Method (BEM)~\cite{VALLEPUGAESPINOSA202028}, Isogeometric BEM~\cite{yu2021thermal,gong2022isogeometric}, FEM-BEM coupling~\cite{gong2025coupled,qin2022application}, Phase Field scheme~\cite{gong2025phase,LONG2025110753}, Scientific Machine Learning~\cite{gong2025pinn}, Finite Difference Method~\cite{raszkowski2020analysis} and Finite Volume Method~\cite{khor2010fvm}. Recent advances have also introduced specialized approaches: Inamdar et al.~\cite{inamdar2024modelling} developed physics-based Digital Twin modeling for thermomechanical degradation prediction; Feng et al.~\cite{feng2023engineered} proposed smoothed finite element algorithms for electro-thermal-mechanical coupling; and Wang et al.~\cite{wang2023structural} integrated machine learning approaches for thermal stress reduction strategies.

With the advancement of commercial finite element software packages (ANSYS, Abaqus, COMSOL), FEM has become the predominant numerical method for thermal and thermomechanical analysis in electronic devices. Recent applications demonstrate its effectiveness across diverse packaging technologies: Wang et al.~\cite{wang2025thermal} integrated finite element analysis with Taguchi optimization to investigate thermal stress in chip-scale package LEDs; Liu et al.~\cite{liu2022efficient} developed coupled thermomechanical simulation techniques using Spectral Element Time Domain methods; Gong et al.~\cite{gong2020thermal} employed COMSOL Multiphysics for three-dimensional multi-chip module BGA modeling to address thermal hotspot issues; and advanced constitutive models such as the Anand model have been applied~\cite{xia2023effect} to characterize SAC305 solder behavior in 3D packaging under thermal cycling from -100 $^\circ$C to 120 $^\circ$C.

Despite FEM's mature theoretical framework and high computational accuracy for regular structures with quality meshes, it faces significant challenges when analyzing multi-scale geometric configurations. Complex geometries require extensive mesh refinement and careful element partitioning to prevent distortion, while geometric scale variations necessitate transitional meshes that can generate millions of elements. This dramatically increases computational cost and may cause convergence issues. The Virtual Element Method (VEM) addresses these limitations through enhanced geometric flexibility~\cite{gong2025VEMthermomechanical}. VEM accommodates arbitrary polygonal elements that simplify complex geometry discretization and enables local refinement through non-matching meshes without compromising overall mesh quality. This capability reduces element count and computational cost while preserving solution accuracy, making VEM particularly suitable for multi-scale problems where traditional meshing approaches become prohibitively expensive. Non-matching meshes refer to the situation where, during the independent discretization of two or more adjacent regions within a computational domain, the positions of interface nodes between these regions do not align, leading to misaligned nodes and edges along their shared interface. During the assembly process, nodes from one region’s interface are inserted into the interface elements of the adjacent region based on their coordinates, and any duplicate nodes are merged. Consequently, the interface nodes become shared between both regions' interface elements, thereby enforcing the continuity of physical quantities such as displacements at the interface.

The VEM was first introduced by Beirao da Veiga et al.~\cite{beirao2013basic} for Poisson's equation and has since been expanded across multiple engineering applications. Early validation studies by Artioli et al.~\cite{artioli2017arbitrary} demonstrated VEM's accuracy through computational benchmarks including Cook's membrane tests. Subsequent research~\cite{sorgente2022role,sorgente2022polyhedral,van2020virtual,van2021virtual} has confirmed VEM's geometric robustness and convergence properties for arbitrary element shapes, including non-convex configurations.
VEM applications have extended beyond linear elasticity to encompass hyperelastic materials~\cite{xu2024high,chi2017some}, contact mechanics~\cite{shen20222,cihan2022virtual}, elastoplasticity~\cite{hudobivnik2019low,wriggers2017low}, and phase field fracture modeling~\cite{liu2023virtual}. Dhanush et al.~\cite{dhanush2019implementation} pioneered VEM implementation for thermomechanical coupling using Abaqus UEL framework. Recently, Gong et al.~\cite{gong2025VEMthermomechanical} analyzed thermoelastic problems using VEM, presenting both the standard virtual element method and a stabilization-free virtual element method~\cite{XU2024116826} for thermomechanical behavior analysis in electronic packaging structures with geometric multi-scale features. Their work demonstrated that VEM's polygonal mesh flexibility enables localized mesh modifications without affecting global mesh structure, making it particularly effective for electronic packaging reliability analysis involving complex geometries.
However, VEM requires projection operators and stabilization terms that increase algorithmic complexity and may result in higher local computational costs compared to conventional FEM.


This work proposes an FE-VE coupled methodology that strategically combines the computational advantages of both methods for problems involving multi-scale geometric features. The approach employs VEM with polygonal discretization in geometrically complex regions, leveraging its tolerance for arbitrary element shapes, while utilizing standard quadrilateral finite elements in regular domains to exploit FEM's computational efficiency. This hybrid strategy balances accuracy and computational cost by applying each method where it performs optimally. The coupled approach is applied to thermomechanical analysis of electronic packaging structures and validated through analytical benchmarks and practical case studies, demonstrating its effectiveness for capturing thermal distributions and stress concentrations in multi-material systems. The framework addresses computational challenges in multi-scale thermomechanical problems by strategically combining FEM's efficiency in structured regions with VEM's geometric flexibility in irregular domains.

The remainder of this paper is organized as follows: Sec.~\ref{sec:mathematical_formulation} presents the mathematical formulation for thermomechanical coupling, including governing equations and weak forms. Sec.~\ref{sec:numerical_discretization} establishes the discrete framework for both FEM and VEM approaches, followed by a comparative analysis. Sec.~\ref{sec:fem_vem_coupling_strategy} describes the FE-VE coupling methodology, including domain decomposition and interface treatment. Sec.~\ref{sec:numerical_examples} validates the proposed approach through numerical examples, and Sec.~\ref{sec:conclusions} provides concluding remarks.

\section{Mathematical formulation of thermomechanical coupling}
\label{sec:mathematical_formulation}

\begin{figure}[htbp]
    \centering
    \includegraphics[width=0.35\textwidth]{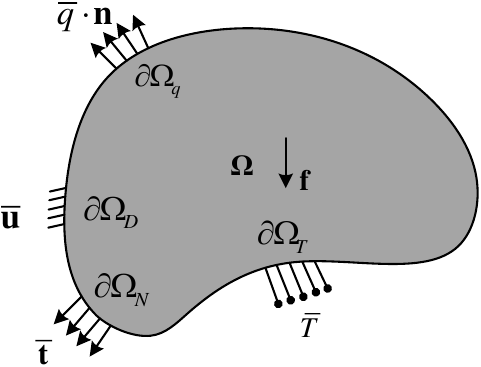}
    \caption{Schematic of the thermoelastic coupling problem domain.}
    \label{fig:thermoelastic_domain}
\end{figure}

We consider a continuous domain $\Omega$ subjected to body force $\mathbf{f}$, as illustrated in Fig.~\ref{fig:thermoelastic_domain}. The domain boundary $\partial\Omega$ is partitioned into two pairs of disjoint regions for the thermomechanical problem
\begin{align}
\partial\Omega &= \partial\Omega_D \cup \partial\Omega_N, \quad \partial\Omega_D \cap \partial\Omega_N = \emptyset, \\
\partial\Omega &= \partial\Omega_T \cup \partial\Omega_q, \quad \partial\Omega_T \cap \partial\Omega_q = \emptyset,
\end{align}
where $\partial\Omega_D$ and $\partial\Omega_N$ represent the displacement and traction boundaries for the mechanical field, while $\partial\Omega_T$ and $\partial\Omega_q$ denote the temperature and heat flux boundaries for the thermal field, respectively. The governing equations for both physical fields are established below.

\subsection{Governing equations}
\label{sec:governing_equations}

For the steady-state heat conduction problem without heat sources, the governing equation takes the form
\begin{equation}
\label{eq:heat_conduction}
\nabla \cdot (\lambda\nabla T) +Q= 0,
\end{equation}
where $\lambda$ is the thermal conductivity and $T$ is the temperature field, $Q$ is the internal heat generation per unit volume.
The boundary conditions on $\partial\Omega = \partial\Omega_T \cup \partial\Omega_q$ are expressed as
\begin{align}
T(\mathbf{x}) &= \overline {T}, \quad \forall \mathbf{x} \in \partial\Omega_T \label{eq:temp_bc}, \\
-\lambda \nabla T(\mathbf{x}) \cdot \mathbf{n} &= \overline {q}, \quad \forall \mathbf{x} \in \partial\Omega_q, \label{eq:flux_bc}
\end{align}
where $\overline {q}$ represents the prescribed heat flux (positive for outward flux) and $\mathbf{n}$ is the outward unit normal vector.

For the thermoelastic problem, the governing equation is
\begin{equation}
\nabla \cdot \boldsymbol{\sigma} + \mathbf{f} = 0, \label{eq:equilibrium}
\end{equation}
where $\boldsymbol{\sigma}$ is the Cauchy stress tensor. The constitutive relation considering thermal strain is
\begin{equation}
\boldsymbol{\sigma} = \mathbb{D} : (\boldsymbol{\varepsilon} - \boldsymbol{\varepsilon}_{\mathrm{th}}), \label{eq:constitutive}
\end{equation}
where $\mathbb{D}$ is the fourth-order elasticity tensor, $\boldsymbol{\varepsilon}$ is the strain tensor, and $\boldsymbol{\varepsilon}_{\mathrm{th}}$ is the thermal strain. The strain tensor is defined as
\begin{equation}
\boldsymbol{\varepsilon}(\mathbf{u}) = \frac{1}{2}\left(\nabla \mathbf{u} + (\nabla \mathbf{u})^\mathrm{T}\right). \label{eq:strain}
\end{equation}
The thermal strain is expressed as
\begin{equation}
\boldsymbol{\varepsilon}_{\mathrm{th}} = \alpha (T - T_0)  \mathbf{I}, \label{eq:thermal_strain}
\end{equation}
where $\mathbf{I}$ is the identity tensor, $\alpha$ is the coefficient of thermal expansion, and $T_0$ is the reference (initial) temperature at which thermal strains are zero.

The boundary conditions on $\partial\Omega = \partial\Omega_D \cup \partial\Omega_N$ are expressed as
\begin{align}
\mathbf{u}(\mathbf{x}) &= \overline {\mathbf{u}}, \quad \forall \mathbf{x} \in \partial\Omega_D, \label{eq:disp_bc} \\
\boldsymbol{\sigma}(\mathbf{x}) \mathbf{n} &= \overline {\mathbf{t}}, \quad \forall \mathbf{x} \in \partial\Omega_N, \label{eq:traction_bc}
\end{align}
where $\overline {\mathbf{u}}$ and $\overline {\mathbf{t}}$ represent the prescribed displacement and traction vectors, respectively.

\subsection{Weak formulation}
\label{sec:weak_formulation}

The weak form of the thermal equilibrium equation is
\begin{equation}
    \label{eq:weak_form_thermal}
  \lambda_k a_T(T,\varphi)+ \int_\Omega Q\varphi\mathrm{~d}\Omega= \ell_T(\varphi),
\end{equation}
where $\varphi$ is the test function (virtual temperature variation). For simplicity, we assume $Q = 0$ (no internal heat generation) in this study.

The corresponding bilinear and linear forms are defined as
\begin{align}
a_T(T, \varphi) &=\lambda\int_{\Omega}  \nabla T \cdot \nabla \varphi \, \mathrm{d}\Omega, \label{eq:bilinear_thermal} \\
\ell_T(\varphi) &= \int_{\partial\Omega_q} \varphi \overline {q} \, \mathrm{d}\Gamma. \label{eq:linear_thermal}
\end{align}
The thermoelastic equilibrium equation becomes: find $\mathbf{u} \in V_u$ such that
\begin{equation}
a_u(\mathbf{u}, \mathbf{v}) = \ell(\mathbf{v}), \quad \forall \mathbf{v} \in V_u^0,
\label{eq:weak_form_mechanical}
\end{equation}
where $\mathbf{u}$ is the displacement field, $\mathbf{v}$ is the test function, $V_u$ is the space of admissible displacement functions satisfying essential boundary conditions, and $V_u^0$ is the space of test functions satisfying homogeneous displacement boundary conditions. The corresponding bilinear and linear forms are defined as~\cite{Wriggers2024virtual}
\begin{align}
a_u(\mathbf{u}, \mathbf{v}) &= \int_{\Omega} \boldsymbol{\varepsilon}(\mathbf{u}) : \mathbb{D} : \boldsymbol{\varepsilon}(\mathbf{v}) \, \mathrm{d}\Omega, \label{eq:bilinear_mechanical} \\
\ell(\mathbf{v}) &= \int_{\Omega} \mathbf{v} \cdot \mathbf{f} \, \mathrm{d}\Omega + \int_{\partial\Omega_N} \mathbf{v} \cdot \overline {\mathbf{t}} \, \mathrm{d}\Gamma + \int_{\Omega} \boldsymbol{\varepsilon}(\mathbf{v}) : \mathbb{D} : \boldsymbol{\varepsilon}_{\mathrm{th}} \, \mathrm{d}\Omega. \label{eq:linear_mechanical}
\end{align}
The linear form includes contributions from body forces, surface tractions, and thermal loads, with the last term in Eq.~\eqref{eq:linear_mechanical} representing the thermomechanical coupling through thermal expansion effects.

\section{Numerical discretization methods}
\label{sec:numerical_discretization}

\subsection{Finite element formulation}
\label{sec:finite_element_formulation}

\subsubsection{Thermal analysis}
\label{sec:thermal_analysis_fem}

For the finite element discretization of the thermal problem, we employ the weak formulation established in Sec.~\ref{sec:weak_formulation}. The domain is discretized into finite elements such that $\Omega\equiv \text{int}\left(\bigcup_i E_i\right)$, where the operator ‘int’ denotes the interior of a set, $\Omega$ represents the union of the interior regions of all elements $E_i$, excluding the boundaries. And element \textit{E\textsubscript{i}} is chosen as a four-node quadrilateral isoparametric element. 
 
According to the Galerkin method, the test functions are chosen as the shape functions, i.e., $\varphi_j = N_j$ for $j=1,2,3,4$. The temperature field within each element is approximated as
\begin{equation}
T(\mathbf{x}) = \sum_{i=1}^4 N_i(\xi, \eta) T_i, 
\label{eq:temp_approximation}
\end{equation}
where $N_i(\xi, \eta)$ are the bilinear shape functions and $(\xi, \eta) \in[-1,1]$ index local coordinates, $T_i$ are the nodal temperatures at the element nodes.

%

The element-wise weak form can be written as
\begin{equation}
    \sum_{i=1}^{4} \int_{E} \lambda (\nabla N_j)^\mathrm{T} \cdot \nabla N_i \, \mathrm{d}E \, T_i = -\int_{\partial E} N_j \overline {q} \, \mathrm{d}\Gamma.
    \label{eq:element_weak_form}
\end{equation}

This leads to the standard finite element matrix equation for each element
\begin{equation}
\mathbf{K}_{\mathrm{FE}}^{th,E} \mathbf{T}^E = \mathbf{f}_{\mathrm{FE}}^{th,E},
\label{eq:element_system}
\end{equation}
where the element thermal stiffness matrix is defined as
\begin{equation}
\mathbf{K}_{\mathrm{FE}}^{th,E} = \int_{E} \lambda \mathbf{B}_T^\mathrm{T} \mathbf{B}_T \, \mathrm{d}E
\label{eq:thermal_stiffness}
\end{equation}
and the element thermal load vector is
\begin{equation}
\mathbf{f}_{\mathrm{FE}}^{th,E} = - \int_{\partial E \cap \partial\Omega_q} \mathbf{N}^\mathrm{T} \bar{q} \, \mathrm{d}\Gamma,
\label{eq:thermal_load}
\end{equation}
where $\mathbf{B}_T$ is the thermal gradient matrix and $\mathbf{N} = [N_1, N_2, N_3, N_4]^\mathrm{T}$ is the shape function vector. The vector $\mathbf{f}_{\mathrm{FE}}^{th,E}$ represents the equivalent nodal loads due to prescribed heat flux boundary conditions.

\label{eq:gaussian_integration}

The global system of equations is assembled from the element contributions using standard finite element assembly
\begin{equation}
\mathbf{K}_{\mathrm{FE}}^{th} \mathbf{T} = \mathbf{f}_{\mathrm{FE}}^{th},
\label{eq:global_thermal_system}
\end{equation}
where the global thermal stiffness matrix and load vector are assembled as
\begin{align}
\mathbf{K}_{\mathrm{FE}}^{th} &= \mathcal{A}_{E} \mathbf{K}_{\mathrm{FE}}^{th,E} \label{eq:stiffness_assembly},\\
\mathbf{f}_{\mathrm{FE}}^{th} &= \mathcal{A}_{E} \mathbf{f}_{\mathrm{FE}}^{th,E}, \label{eq:load_assembly}
\end{align}
where $\mathcal{A}$ denotes the standard assembly operator that maps element degrees of freedom (Dofs) to global Dofs according to element connectivity.

\subsubsection{Thermoelastic analysis}
\label{sec:thermoelastic_analysis_fem}

For the finite element discretization of the mechanical problem, we employ the weak formulation established in Sec.~\ref{sec:weak_formulation}. Within each element, the displacement field is interpolated using the same shape functions as the thermal analysis
\begin{equation}
\mathbf{u}(\mathbf{x}) = \sum_{i=1}^{4} N_i(\xi, \eta) \mathbf{u}_i = \mathbf{N} \mathbf{u}^E,
\label{eq:displacement_approximation}
\end{equation}
where $\mathbf{u}_i$ is the $i$th nodal displacement vector, $\mathbf{u}^E = [\mathbf{u}_1^\mathrm{T}, \mathbf{u}_2^\mathrm{T}, \mathbf{u}_3^\mathrm{T}, \mathbf{u}_4^\mathrm{T}]^\mathrm{T}$ is the element displacement vector, and $\mathbf{N}$ is the shape function matrix
\label{eq:shape_function_matrix}

The strain-displacement relationship is given by
\begin{equation}
\boldsymbol{\varepsilon} = \mathbf{B}_u \mathbf{u}^E,
\label{eq:strain_displacement}
\end{equation}
where $\mathbf{B}_u$ is the strain-displacement matrix that relates the element nodal displacements to the strain components. The matrix $\mathbf{B}_u$ contains the derivatives of the shape functions with respect to the global coordinates, obtained through the coordinate transformation using the Jacobian matrix.

\label{eq:mechanical_stiffness}

The element load vector includes contributions from external loads and thermal effects
\begin{equation}
\mathbf{F}_{\mathrm{FE}}^{u,E} = \mathbf{F}_{\mathrm{ext}}^E + \mathbf{F}_{\mathrm{FE}}^{th,E},
\label{eq:element_total_load}
\end{equation}
where $\mathbf{F}_{\mathrm{ext}}^E$ represents the external load vector from body forces and surface tractions.

After assembly, the global system of equations becomes
\begin{equation}
\mathbf{K}_{\mathrm{FE}}^{u} \mathbf{u} = \mathbf{F}_{\mathrm{FE}}^{u},
\label{eq:global_mechanical_system}
\end{equation}
where $\mathbf{K}_{\mathrm{FE}}^{u}$ and $\mathbf{F}_{\mathrm{FE}}^{u}$ are assembled from the element contributions.

The thermal load vector, arising from the thermomechanical coupling, is computed as
\begin{equation}
\mathbf{F}_{\mathrm{FE}}^{th} = \sum_{E} \int_{E} \mathbf{B}_u^\mathrm{T} \mathbb{D} \boldsymbol{\varepsilon}_{\mathrm{th}} \, \mathrm{d}E = \sum_{E} \int_{E} \mathbf{B}_u^\mathrm{T}  \mathbb{D} \alpha (T - T_0) \mathbf{I} \, \mathrm{d}E,
\label{eq:thermal_load_mechanical}
\end{equation}
where $\mathbf{I} = [1, 1, 0]^\mathrm{T}$ is the thermal expansion vector for plane stress/strain conditions, $T(\xi, \eta) = \sum_{i=1}^4 N_i(\xi, \eta) T_i$ is the temperature field obtained from the thermal analysis.

\subsection{Virtual element formulation}
\label{sec:virtual_element_formulation}

\subsubsection{Thermal analysis}
\label{sec:vem_thermal_analysis}
For the heat conduction problem, the scalar virtual element function space is defined as
\begin{equation}
\mathcal{V}_T = \left\{T \in \mathcal{H}^1(\Omega) : \left.T\right|_E \in \mathcal{V}_T^h(E) \text{ for all } E, \, T = \bar{T} \text{ on } \partial\Omega_T\right\},
\label{eq:global_vem_space}
\end{equation}
where the local virtual element space on element $E$ is defined as
\begin{equation}
\mathcal{V}_T^h(E) = \left\{T \in \mathcal{H}^1(E) : \Delta T = 0 \text{ in } E, \, \left.T\right|_{\partial E} \in \mathcal{P}_k(\partial E)\right\}.
\label{eq:local_vem_space}
\end{equation}
Here, $\mathcal{P}_k(\partial E)$ denotes the space of polynomial functions of degree not exceeding $k$ on the boundary of element $E$. In this work, $k = 1$ is chosen, leading to linear virtual elements.
The scaled polynomial space $\mathcal{P}_k(E)$ on element $E$ is spanned by the monomial basis. For $k = 1$, the polynomial space is~\cite{dhanush2019implementation}
\begin{equation}
\mathcal{P}_1(E) = \{p_1,p_2,p_3\}=\text{span}\{1, \zeta, \wp\},
\label{eq:polynomial_space}
\end{equation}
where the scaled coordinates $\zeta$ and $\wp$ are defined as~\cite{XU2024116826}
\begin{align}
\zeta &= \frac{x - \bar{x}}{h_E}, \quad \wp = \frac{y - \bar{y}}{h_E},
\label{eq:scaled_coordinates}
\end{align}
in which $\bar{\mathbf{x}} = (\bar{x}, \bar{y})$ represents the centroid of element $E$, and $h_E$ denotes the characteristic element size defined as the maximum distance between any two nodes in element $E$.

Since the element shape functions do not have explicit expressions, a projection operator is defined to compute the bilinear forms. The approximate function space $\mathcal{V}_T^h(E)$ on element $E$ is mapped to the polynomial space $\mathcal{P}_k(E)$ through the projection operator
\begin{equation}
\Pi: \mathcal{V}_T^h(E) \rightarrow \mathcal{P}_k(E), 
\label{eq:projection_definition}
\end{equation}
which satisfies the orthogonality condition
\begin{equation}
a_T^E\left(\psi - \Pi \psi, p\right) = 0, \quad \forall p \in \mathcal{P}_k(E)
\label{eq:projection_definition}
\end{equation}
where $\psi$ is the basis function of $\mathcal{V}_T^h(E)$.

Expanding Eq.~\eqref{eq:projection_definition} yields
\begin{equation}
\int_E \nabla\left(\Pi \psi\right) \cdot \nabla p \, \mathrm{d} E = \int_E \nabla \psi \cdot \nabla p \, \mathrm{d} E.
\label{eq:projection_expanded}
\end{equation}

The matrix representation $\tilde{\boldsymbol{\Pi}}$ of the projection operator $\Pi$ can be obtained by~\cite{MENGOLINI2019995,dhanush2019implementation}
\begin{equation}\label{eq:projection_expanded2-thermal-1}
    \tilde{\mathbf{\Pi}}=\mathbf{G}^{-1} \mathbf{B},
\end{equation}
where
\begin{equation}
    \mathbf{G}_{\alpha\beta} = \lambda \int_E \nabla p_\alpha \cdot \nabla p_\beta \, \mathrm{d}E
    \label{eq:G_matrix}
    \end{equation}
    \begin{equation}
    \mathbf{B}_{i\alpha} = \lambda \int_{\partial E} \psi_i \cdot (\nabla p_\alpha) \mathbf{n} \, \mathrm{d}\Gamma
    \label{eq:B_matrix}
    \end{equation}
where $\alpha, \beta = 1, 2, 3$ index the polynomial basis functions $p_\alpha \in \mathcal{P}_1(E)$, and  $i = 1, 2, \ldots, n_v$ indexes the basis function $\psi_i$ of $\mathcal{V}_T^h(E)$ corresponding to the $i$-th vertice of the polygonal element.

For the temperature field $T^h$, the bilinear form on each element can be decomposed as
\begin{align}
a_T^E\left(T^h, \varphi^h\right) &= a_T^E\left(T^h - \Pi T^h + \Pi T^h, \varphi^h - \Pi \varphi^h + \Pi \varphi^h\right) \nonumber \\
&= a_T^E\left(\Pi T^h, \Pi \varphi^h\right) + a_T^E\left(T^h - \Pi T^h, \varphi^h - \Pi \varphi^h\right).
\label{eq:bilinear_decomposition}
\end{align}
In Eq.~\eqref{eq:bilinear_decomposition}, the first term on the right-hand side represents the consistency term computed from the integral of projected functions in element $E$. The second term is the stabilization term that represents the energy contribution of projection residuals  $T^h-\Pi^{\nabla} T^h$ and $\varphi^h-\Pi^{\nabla} \varphi^h$ . This stabilization term is essential in VEM formulations to ensure that the stiffness matrix is positive definite and full rank.

The element thermal stiffness matrix consists of two parts
\begin{equation}
\mathbf{K}_{\mathrm{VE}}^{th,E} = \mathbf{K}_E^{th,\mathrm{c}} + \mathbf{K}_E^{th,\mathrm{s}},
\label{eq:vem_stiffness_decomposition}
\end{equation}
where $\mathbf{K}_E^{th,\mathrm{c}}$ is the consistency stiffness matrix and $\mathbf{K}_E^{th,\mathrm{s}}$ is the stabilization stiffness matrix.



For the heat conduction problem, the matrix $\mathbf{D}_{n_v \times n_k}$ represents the evaluation of polynomial basis functions $p_\alpha$ at the degrees of freedom locations
\begin{equation}
\mathbf{D}_{j\beta} = \text{dof}_j(p_\beta), \quad j = 1, 2, \ldots, n_v, \quad \beta = 1, 2, 3
\label{eq:dof_matrix}
\end{equation}
where $\text{dof}_j(\cdot)$ denotes the $j$-th degree of freedom value (nodal evaluation for the temperature field), $n_v$ is the number of vertices for the polygonal element, and $n_k = 3$ is the dimension of the polynomial space $\mathcal{P}_1(E)$. 

The stabilization stiffness matrix is expressed as~\cite{MENGOLINI2019995}
\begin{equation}
\mathbf{K}_E^{th,\mathrm{s}} = \tau^h \, \mathrm{tr}\left(\mathbf{K}_E^{th,\mathrm{c}}\right) (\mathbf{I} - \boldsymbol{\Pi})^\mathrm{T}  (\mathbf{I} - \boldsymbol{\Pi}),
\label{eq:stabilization_matrix}
\end{equation}
where the projection operator is given by
\begin{equation}
\boldsymbol{\Pi} = \mathbf{D}\tilde{\boldsymbol{\Pi}},
\label{eq:projection_relation}
\end{equation}
representing the mapping from the VEM space to the degrees of freedom space through the polynomial space. Here $\mathrm{tr}(\cdot)$ denotes the trace operator, and $\mathbf{I}$ is the $n_v \times n_v$ identity matrix, $\tau^h = 1/2$ is the stabilization parameter~\cite{MENGOLINI2019995}. According to Reference \cite{artioli2017arbitrary}, VEM results show limited sensitivity to the stabilization parameter $\tau^h$. The numerical solution maintains stable accuracy and convergence across a wide range of $\tau^h$ values, with optimal performance achieved at $\tau^h = 1/2$.

\subsubsection{Thermoelastic analysis}
\label{sec:vem_thermoelastic_analysis}

For elasticity problems, the global virtual element space is defined as
\begin{equation}
\mathcal{V}_{\mathbf{u}} = \left\{\mathbf{v} \in [\mathcal{H}^1(\Omega)]^2 : \left.\mathbf{v}\right|_E \in \mathcal{V}_{\mathbf{u}}^h(E) \text{ for all } E, \, \mathbf{v} = \overline{\mathbf{u}} \text{ on } \partial\Omega_D\right\},
\label{eq:global_vem_vector_space}
\end{equation}
where $\mathbf{v}$ represents any vector function in the space, and the local virtual element space on element $E$ is defined as~\cite{herrera2023numerical}
\begin{equation}
\begin{aligned}
\mathcal{V}_{\mathbf{u}}^h(E) := \bigg\{\mathbf{v}^h \in [\mathcal{H}^1(E)]^2 : &\left.\mathbf{v}^h\right|_{\partial E} \in [\mathcal{C}^0(\partial E)]^2, \\
&\left.\mathbf{v}^h\right|_E \in [\mathbf{M}_k(E)]^2, \\
&\left.\Delta \mathbf{v}^h\right|_E \in [\mathbf{M}_{k-2}(E)]^2\bigg\},
\end{aligned}
\label{eq:local_vem_vector_space}
\end{equation}
where $\mathbf{v}^h$ represents the discrete vector function within the element and $\Delta$ represents the Laplace Operator such that $\Delta\mathbf{v}^h=\nabla^2\mathbf{v}^h=\sum_{i=1}^n\frac{\partial^2\mathbf{v}^h}{\partial x_i^2}$.

For interpolation order $k = 1$, the constraint effectively requires $\Delta \mathbf{v}^h = 0$ within each element, since $\mathbf{M}_{-1}(E)$ contains only the zero polynomial.
The local polynomial space $\mathbf{M}_1(E)$ consists of the following basis functions~\cite{CHEN202388}
\begin{equation}
\mathbf{M}_1(E) = \left\{ \mathbf{m}_1, \mathbf{m}_2, \mathbf{m}_3, \mathbf{m}_4, \mathbf{m}_5, \mathbf{m}_6 \right\},
\label{eq:polynomial_basis}
\end{equation}
where
$$
\mathbf{m}_1=\begin{pmatrix} 1 \\ 0 \end{pmatrix},~ \mathbf{m}_2 = \begin{pmatrix} 0 \\ 1 \end{pmatrix},~ \mathbf{m}_3 = \begin{pmatrix} -\wp \\ \zeta \end{pmatrix},~ \mathbf{m}_4 = \begin{pmatrix} \wp \\ \zeta \end{pmatrix}, ~ \mathbf{m}_5 = \begin{pmatrix} \zeta \\ 0 \end{pmatrix},~ \mathbf{m}_6 = \begin{pmatrix} 0 \\ \wp \end{pmatrix}.
$$

For a polygon with $n_v$ vertices, the Dof for the two-dimensional vector problem are $2n_v$. The local virtual element space $\mathcal{V}_{\mathbf{u}}^h(E)$ is mapped to the polynomial space $\mathbf{M}_k(E)$ through the projection operator 
\begin{equation}
    \Pi^{\nabla}: \mathcal{V}_{\mathbf{u}}^h(E) \rightarrow \mathbf{M}_k(E).
\end{equation}

Assuming $\boldsymbol{\phi}$ is the basis function of $\mathcal{V}_{\mathbf{u}}^h(E)$, the projection operator satisfies the orthogonality condition,
\begin{equation}
a_u^E\left(\boldsymbol{\phi} - \Pi^{\nabla} \boldsymbol{\phi}, \mathbf{m} \right) = 0, \quad \forall \mathbf{m} \in \mathbf{M}_k(E)
\label{eq:orthogonality_condition}
\end{equation}
where the local bilinear form is defined as
\begin{equation}
a_u^E(\boldsymbol{\phi}, \mathbf{m}) = \int_E \boldsymbol{\varepsilon}(\boldsymbol{\phi}) : \mathbb{D} : \boldsymbol{\varepsilon}(\mathbf{m}) \, \mathrm{d}E.
\label{eq:local_bilinear_form_elasticity}
\end{equation}

According to Eq.~\eqref{eq:orthogonality_condition}, we can obtain 
\begin{equation}
a_u^E\left(\boldsymbol{\phi}, \mathbf{m} \right) =a_u^E\left( \Pi^{\nabla} \boldsymbol{\phi}, \mathbf{m} \right) , \quad \forall \mathbf{m} \in \mathbf{M}_k(E)
\label{eq:orthogonality_condition_n}
\end{equation}

Based on the references~\cite{xu2024high,BERBATOV2021351}, $\Pi^{\nabla}\boldsymbol{\phi}$ can be expressed as a linear combination of the polynomial basis functions $\mathbf{m}$ in $\mathbf{M}_k(E)$. Thus,
\begin{equation}
\Pi^{\nabla}\boldsymbol{\phi}=\mathbf{m}\tilde{\boldsymbol{\Pi}}^{*\nabla}  ,
\label{eq:local_bilinear_form_elasticity}
\end{equation}
By using Eq.~\eqref{eq:orthogonality_condition_n}, we can derive the following matrix equation to solve for $\tilde{\boldsymbol{\Pi}}^{*\nabla}$,
\begin{equation}
\mathbf{M} \tilde{\boldsymbol{\Pi}}^{*\nabla} = \overline {\mathbf{B}},
\label{eq:projection_matrix_equation}
\end{equation}
where $\tilde{\boldsymbol{\Pi}}^{*\nabla}$ is the matrix formulation of the Ritz projection operator~\cite{xu2024high}, and
\begin{equation}
\mathbf{M} = \int_E \hat{\boldsymbol{\varepsilon}}^\mathrm{T}(\mathbf{m}) \hat{\mathbb{D}} \hat{\boldsymbol{\varepsilon}}(\mathbf{m}) \, \mathrm{d}E, \quad 
\overline{\mathbf{B}} = \int_E \hat{\boldsymbol{\varepsilon}}^\mathrm{T}(\mathbf{m}) \hat{\mathbb{D}} \hat{\boldsymbol{\varepsilon}}(\boldsymbol{\phi}) \, \mathrm{d}E,
\end{equation}
where $\hat{\square}$ represents the Voigt form. 

For $k=1$, the strain matrix is
\begin{equation}
    \label{r3.varepsilon_m}
    \hat{\boldsymbol{\varepsilon}}(\mathbf{m}) = 
    \begin{bmatrix}
    \frac{\partial \mathbf{m}_1}{\partial x} & \cdots & \frac{\partial \mathbf{m}_6}{\partial x}\\
    \frac{\partial \mathbf{m}_1}{\partial y} & \cdots & \frac{\partial \mathbf{m}_6}{\partial y} \\
    \frac{\partial \mathbf{m}_1}{\partial x}+\frac{\partial \mathbf{m}_1}{\partial y} & \cdots & \frac{\partial \mathbf{m}_6}{\partial x}+\frac{\partial \mathbf{m}_6}{\partial y}
    \end{bmatrix}.
\end{equation}

Besides, the matrix $\overline{\mathbf{B}}$ can be expressed using integration by parts as
\begin{equation}
    \small 
    \overline{\mathbf{B}} = \int_{E}\hat{\boldsymbol{\varepsilon}}^\mathrm{T}(\mathbf{m})\hat{\mathbb{D}}\hat{\boldsymbol{\varepsilon}}(\boldsymbol{\phi}) \, \mathrm{d}E
= -\int_E\left[\nabla\cdot\left(\hat{\mathbb{D}}\hat{\boldsymbol{\varepsilon}}(\mathbf{m})\right)\right]^\mathrm{T}\boldsymbol{\phi} \, \mathrm{d}E +
\int_{\partial E}\left[\hat{\mathbb{D}}\hat{\boldsymbol{\varepsilon}}(\mathbf{m})\right]^\mathrm{T} \mathbf{n} \cdot \boldsymbol{\phi} \, \mathrm{d}\Gamma.
\label{eq:B_matrix_green}
\end{equation}
Since $\mathbf{m} \in \mathbf{M}_k(E)$ and the divergence term involves polynomials of degree $k-2$, for $k = 1$, the first term on the right-hand side vanishes.

The matrix $\overline{\mathbf{D}}_{2n_v \times n_k}$ represents the evaluation of polynomial basis functions $\mathbf{m}_\alpha$ at the degrees of freedom locations
\begin{equation}
\overline{\mathbf{D}}_{i\alpha} = \text{dof}_i(\mathbf{m}_\alpha), \quad i = 1, 2, \ldots, 2n_v, \quad \alpha = 1, 2, \ldots, n_k,
\label{eq:dof_matrix}
\end{equation}
where $\text{dof}_i(\cdot)$ denotes the $i$-th Dof value (nodal evaluation for the displacement field) and the dimension of the polynomial space $\mathbf{M}_k(E)$ for $k=1$ is $n_k = 6$.

The relationship between the matrices is given by the projection system. The projection operator is determined from Eq.~\eqref{eq:projection_matrix_equation}
\begin{equation}
\tilde{\boldsymbol{\Pi}}^{*\nabla} = \mathbf{M}^{-1} \overline{\mathbf{B}} 
\end{equation}
subject to boundary constraint conditions.

For $\mathbf{u}^h, \mathbf{v}^h \in \mathcal{V}_{\mathbf{u}}^h(E)$, the bilinear form within an element can be decomposed as
\begin{align}
a_u^E\left(\mathbf{u}^h, \mathbf{v}^h\right) &= a_u^E\left(\boldsymbol{\Pi}^{\nabla} \mathbf{u}^h, \boldsymbol{\Pi}^{\nabla} \mathbf{v}^h\right) + a_u^E\left(\mathbf{u}^h - \boldsymbol{\Pi}^{\nabla} \mathbf{u}^h, \mathbf{v}^h - \boldsymbol{\Pi}^{\nabla} \mathbf{v}^h\right).
\label{eq:vem_bilinear_decomposition}
\end{align}

The element stiffness matrix is decomposed as
\begin{equation}
\mathbf{K}_{\mathrm{VE}}^{u,E} = \mathbf{K}_E^{u,\mathrm{c}} + \mathbf{K}_E^{u,\mathrm{s}},
\label{eq:vem_stiffness_total}
\end{equation}
where $\mathbf{K}_E^{u,\mathrm{c}}$ is the consistency part and $\mathbf{K}_E^{u,\mathrm{s}}$ is the stabilization part.

The consistency stiffness matrix is expressed as
\begin{equation}
\mathbf{K}_E^{u,\mathrm{c}} = \left(\tilde{\boldsymbol{\Pi}}^{*\nabla}\right)^\mathrm{T} \mathbf{M} \tilde{\boldsymbol{\Pi}}^{*\nabla}.
\label{eq:vem_consistency_stiffness}
\end{equation}

The stabilization stiffness matrix is expressed as~\cite{MENGOLINI2019995}
\begin{equation}
\mathbf{K}_E^{u,\mathrm{s}} = \tau^h \, \mathrm{tr}\left(\mathbf{K}_E^{u,\mathrm{c}}\right) \left(\mathbf{I} - \tilde{\boldsymbol{\Pi}}^{\nabla}\right)^\mathrm{T} \left(\mathbf{I} - \tilde{\boldsymbol{\Pi}}^{\nabla}\right),
\label{eq:vem_stabilization_stiffness}
\end{equation}
where $\tilde{\boldsymbol{\Pi}}^{\nabla} = \overline{\mathbf{D}} \tilde{\boldsymbol{\Pi}}^{*\nabla}$, $\tau^h$ is a stabilization parameter taken as $1/2$.

Considering thermal strain, the total load vector can be expressed as
\begin{equation}
\mathbf{F}_{\mathrm{VE}}^{u} = \mathbf{F}_{\mathrm{ext}} + \mathbf{F}_{\mathrm{VE}}^{th},
\label{eq:vem_total_load}
\end{equation}
where $\mathbf{F}_{\mathrm{ext}}$ represents the external mechanical load and $\mathbf{F}_{\mathrm{VE}}^{th}$ represents the global thermal load, computed as
\begin{equation}
\mathbf{F}_{\mathrm{VE}}^{th} = \sum_{E}\left(\tilde{\boldsymbol{\Pi}}^{*\nabla}\right)^\mathrm{T} \int_E \hat{\boldsymbol{\varepsilon}}^\mathrm{T}(\mathbf{m})\hat{\mathbb{D}} \hat{\boldsymbol{\varepsilon}}_{\mathrm{th}} \, \mathrm{d}E, 
\label{eq:vem_thermal_load}
\end{equation}
where $\tilde{\boldsymbol{\Pi}}^{*\nabla}$ is the Ritz matrix representation of the projection operator  $\tilde{\boldsymbol{\Pi}}^{\nabla}$ can be solved from Eq.~\eqref{eq:projection_matrix_equation}~\cite{xu2024high}. The thermal strain $\hat{\boldsymbol{\varepsilon}}_{\mathrm{th}} = \alpha (T - T_0) [1, 1, 0]^\mathrm{T}$ is expressed in Voigt form. 

\subsection{Comparative analysis of FEM and VEM}
\label{sec:fem_vem_comparison}

This section presents a comparative analysis between the FEM and VEM for thermoelastic problems. The comparison examines the fundamental differences in discretization strategies, computational efficiency, and practical implementation aspects. Table~\ref{tab:fem_vem_comparison} summarizes the key characteristics and trade-offs between these two numerical approaches across various computational and theoretical criteria.

    \begin{table}[htbp]
        \centering
        \caption{Key differences between FEM and VEM for thermomechanical analysis}
        \label{tab:fem_vem_comparison}
        \footnotesize
        \begin{tabular}{lp{4.5cm}p{4.5cm}}
        \toprule
        & \textbf{FEM} & \textbf{VEM} \\ 
        \midrule
        Mesh & Standard elements (tri/quad) & Polygonal elements of arbitrary shape \\
        Implementation & Direct stiffness assembly & Requires projection $\Pi^\nabla$ and stabilization \\
        Best suited for & Regular structured domains & Complex geometries with irregular interfaces \\
        Dofs per element & Fixed (e.g., 4 for Q4) & Variable ($n_v$ vertices) \\
        Computational cost & More efficient for simple geometries & Additional overhead from projections \\
        \bottomrule
        \end{tabular}
    \end{table}


The comparative analysis reveals that FEM possesses a well-established theoretical foundation and demonstrates superior computational efficiency for regular meshes and standard geometric configurations. The method is particularly effective for standard engineering problems where the geometry can be adequately discretized using conventional element shapes. In contrast, VEM excels in handling geometrically complex problems through its ability to accommodate arbitrary polygonal elements, providing exceptional mesh flexibility and geometric adaptability. This capability is particularly valuable when dealing with irregular domains or when mesh generation becomes challenging with traditional finite elements.

To leverage the complementary strengths highlighted in Table~\ref{tab:fem_vem_comparison}, this study proposes a FE-VE coupling approach that strategically combines both methods. The coupling strategy assigns FEM to regular regions where computational efficiency is paramount, while utilizing VEM in geometrically complex areas where mesh flexibility is essential. This hybrid approach aims to achieve optimal computational efficiency while maintaining solution accuracy for multi-scale geometric structures encountered in electronic packaging applications.

\section{FE-VE coupling strategy}
\label{sec:fem_vem_coupling_strategy}

This section presents a comprehensive FE-VE coupling strategy specifically designed to address geometric multi-scale challenges in electronic packaging applications. The methodology strategically combines the computational efficiency of finite elements with the geometric flexibility of virtual elements through systematic domain decomposition and robust interface treatment.
The coupling approach is based on three key components: (i) domain partitioning criteria that identify regions suitable for each method, (ii) interface compatibility conditions that ensure continuity across FE-VE boundaries, and (iii) solution transfer mechanisms that maintain accuracy during the coupling process. The strategy aims to optimize computational resources while preserving solution quality for complex geometries encountered in modern electronic systems.

\subsection{Domain decomposition}
\label{sec:domain_decomposition}

\begin{figure}[htbp]
    \centering
    \begin{minipage}{0.41\textwidth}
        \centering
        \includegraphics[width=\textwidth]{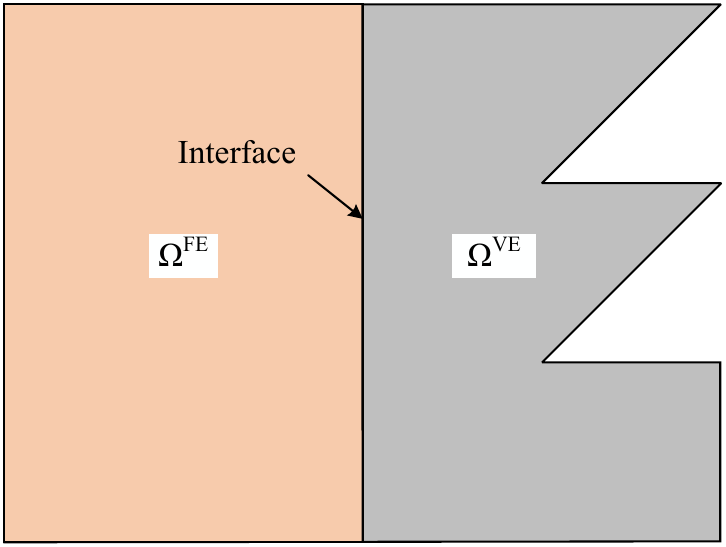}
        \subcaption{Domain partitioning}
        \label{fig:domain_partition}
    \end{minipage}
    \hspace{0.02\textwidth}
    \begin{minipage}{0.44\textwidth}
        \centering
        \includegraphics[width=\textwidth]{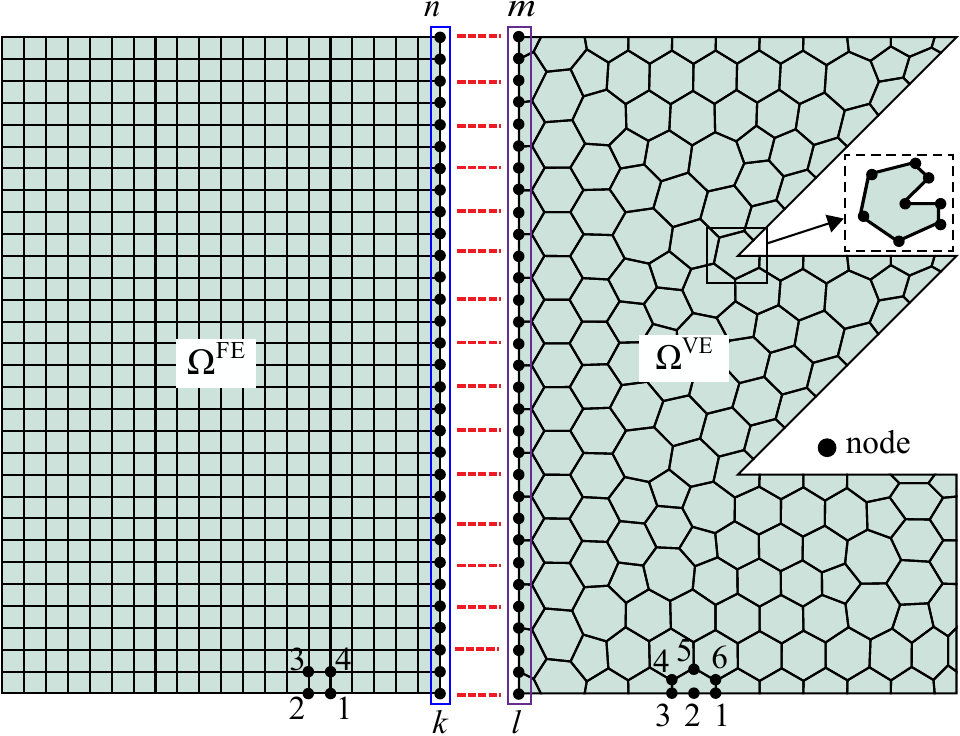}
        \subcaption{Interface mesh discretization}
        \label{fig:mesh_coupling}
    \end{minipage}
    \caption{FE-VE coupling approach. (a) The computational domain is partitioned into finite element and virtual element regions based on geometric complexity. (b) The interface mesh discretization ensures compatibility between FE and VE domains through coincident nodes along the coupling boundary.}
    \label{fig:fem_vem_coupling}
\end{figure}

As illustrated in Fig.~\ref{fig:domain_partition}, the computational domain $\Omega$ is partitioned into two complementary regions based on geometric complexity: the finite element domain $\Omega^{\mathrm{FE}}$ and the virtual element domain $\Omega^{\mathrm{VE}}$, such that $\Omega = \Omega^{\mathrm{FE}} \cup \Omega^{\mathrm{VE}}$ with $\Omega^{\mathrm{FE}} \cap \Omega^{\mathrm{VE}} = \emptyset$. The FEM is deployed in $\Omega^{\mathrm{FE}}$ to efficiently handle regular geometric regions with structured quadrilateral meshes, while the VEM is employed in $\Omega^{\mathrm{VE}}$ to accommodate complex geometric configurations using arbitrary polygonal elements.

The coupling between the two domains is achieved through interface treatment at the shared boundary $\Gamma^{\mathrm{inter}} = \partial\Omega^{\mathrm{FE}} \cap \partial\Omega^{\mathrm{VE}}$. Fig.~\ref{fig:mesh_coupling} demonstrates the discretization strategy, where the finite element region utilizes a regular quadrilateral mesh while the virtual element region employs polygonal elements that conform to the complex geometric boundaries.
The interface compatibility is established through coincident nodes along $\Gamma^{\mathrm{inter}}$, where finite element boundary nodes and their corresponding virtual element boundary nodes share identical spatial coordinates. This nodal correspondence ensures displacement continuity across the coupling boundary, expressed as
\begin{equation}
\mathbf{u}^{\mathrm{FE}}|_{\Gamma^{\mathrm{inter}}} = \mathbf{u}^{\mathrm{VE}}|_{\Gamma^{\mathrm{inter}}}.
\label{eq:interface_continuity}
\end{equation}

As shown in Fig.~\ref{fig:mesh_coupling}, the interface nodes are strategically positioned to ensure compatibility along the coupling boundary. The figure illustrates typical element configurations: a four-node finite element (numbered 1, 2, 3, 4) in the structured quadrilateral mesh region and a six-node virtual element (numbered 1, 2, 3, 4, 5, 6) in the polygonal mesh region. The actual interface nodes are indexed from $k$ to $n$ on the finite element side and from $l$ to $m$ on the virtual element side, where coincident nodal positioning ensures displacement continuity across the FE-VE boundary.

\subsection{Interface treatment}
\label{sec:interface_treatment}

The coupling between finite element and virtual element domains requires systematic treatment of interface variables to ensure displacement continuity and force equilibrium. This section presents the partitioning strategy and assembly procedure for establishing the coupled system.

\subsubsection{Matrix partitioning for interface treatment}

Both finite element and virtual element domains are partitioned to separate interior and interface Dofs based on the nodal positions illustrated in Fig.~\ref{fig:mesh_coupling}.

\paragraph{Finite element domain}
The finite element system is partitioned into block matrix form
\begin{equation}
\begin{bmatrix}
\mathbf{K}_{FF} & \mathbf{K}_{FI} \\
\mathbf{K}_{IF} & \mathbf{K}_{II}^{\mathrm{FE}}
\end{bmatrix}
\begin{Bmatrix}
\mathbf{u}_{F} \\
\mathbf{u}_{I}^{\mathrm{FE}}
\end{Bmatrix}
=
\begin{Bmatrix}
\mathbf{F}_{F} \\
\mathbf{F}_{I}^{\mathrm{FE}}
\end{Bmatrix},
\label{eq:fem_block_system}
\end{equation}
where $\mathbf{u}_{I}^{\mathrm{FE}} = [u_k, u_{k+1}, \ldots, u_n]^\mathrm{T}$ represents the displacement variables at interface nodes (nodes $k$ to $n$ in Fig.~\ref{fig:mesh_coupling}), $\mathbf{u}_{F}$ represents interior node displacements, and $\mathbf{K}_{II}^{\mathrm{FE}}$ is the interface stiffness contribution from the finite element domain.

\paragraph{Virtual element domain}
Similarly, the virtual element system is partitioned as
\begin{equation}
\begin{bmatrix}
\mathbf{K}_{VV} & \mathbf{K}_{VI} \\
\mathbf{K}_{IV} & \mathbf{K}_{II}^{\mathrm{VE}}
\end{bmatrix}
\begin{Bmatrix}
\mathbf{u}_{V} \\
\mathbf{u}_{I}^{\mathrm{VE}}
\end{Bmatrix}
=
\begin{Bmatrix}
\mathbf{F}_{V} \\
\mathbf{F}_{I}^{\mathrm{VE}}
\end{Bmatrix},
\label{eq:vem_block_system}
\end{equation}
where $\mathbf{u}_{I}^{\mathrm{VE}} = [u_l, u_{l+1}, \ldots, u_m]^\mathrm{T}$ represents interface node displacements (nodes $l$ to $m$ in Fig.~\ref{fig:mesh_coupling}), $\mathbf{u}_{V}$ represents interior virtual element nodes, and $\mathbf{K}_{II}^{\mathrm{VE}}$ is the interface stiffness contribution from the virtual element domain.

\subsubsection{Strong coupling implementation}

At the FE-VE interface, a strong coupling scheme is employed where corresponding nodes are matched to ensure displacement compatibility
\begin{equation}
\mathbf{u}_{I}^{\mathrm{FE}} = \mathbf{u}_{I}^{\mathrm{VE}} = \mathbf{u}_{I}.
\label{eq:displacement_compatibility}
\end{equation}

The strong coupling implementation follows a systematic assembly procedure. As established in Eqs.~\eqref{eq:fem_block_system} and \eqref{eq:vem_block_system}, the finite element and virtual element systems are constructed independently in their respective domains with appropriate partitioning for interior and interface degrees of freedom. 

Due to the displacement compatibility condition at the interface (Eq.~\eqref{eq:displacement_compatibility}), the interface degrees of freedom are unified. The coupling is achieved by recognizing that interface nodes belong simultaneously to both domains, and their stiffness and load contributions must be properly combined. By assembling the partitioned systems from Eqs.~\eqref{eq:fem_block_system} and \eqref{eq:vem_block_system}, a coupled global system is obtained
\begin{equation}
\begin{bmatrix}
\mathbf{K}_{FF} & \mathbf{K}_{FI} & \mathbf{0} \\
\mathbf{K}_{IF} & \mathbf{K}_{II}^{\mathrm{FE}} + \mathbf{K}_{II}^{\mathrm{VE}} & \mathbf{K}_{IV} \\
\mathbf{0} & \mathbf{K}_{VI} & \mathbf{K}_{VV}
\end{bmatrix}
\begin{Bmatrix}
\mathbf{u}_F \\
\mathbf{u}_I \\
\mathbf{u}_V
\end{Bmatrix}
=
\begin{Bmatrix}
\mathbf{F}_F \\
\mathbf{F}_I^{\mathrm{FE}} + \mathbf{F}_I^{\mathrm{VE}} \\
\mathbf{F}_V
\end{Bmatrix},
\label{eq:coupled_system}
\end{equation}
where subscripts $F$, $I$, and $V$ denote finite element interior nodes, interface nodes, and virtual element interior nodes, respectively. The key features of this coupling formulation are
(i) The off-diagonal blocks $\mathbf{K}_{FI}$ and $\mathbf{K}_{IF}$ represent the coupling between FE interior nodes and interface nodes, arising from finite elements adjacent to the interface.
(ii) Similarly, $\mathbf{K}_{IV}$ and $\mathbf{K}_{VI}$ represent the coupling between interface nodes and VE interior nodes from virtual elements adjacent to the interface.
(iii) The diagonal interface block $\mathbf{K}_{II}^{\mathrm{FE}} + \mathbf{K}_{II}^{\mathrm{VE}}$ represents the additive stiffness contributions from both domains. This linear superposition is mathematically valid because both methods discretize the same underlying continuum problem, and the interface nodes serve as connection points where element-level contributions are naturally summed during standard assembly procedures.
(iv) The interface load vector $\mathbf{F}_I^{\mathrm{FE}} + \mathbf{F}_I^{\mathrm{VE}}$ combines load contributions from both domains, including external mechanical loads and thermal expansion effects, ensuring proper load transfer across the coupling boundary.
(v) The zero blocks in the upper-right and lower-left corners reflect the absence of direct coupling between FE and VE interior nodes, as these nodes interact only through the shared interface.

This strong coupling approach ensures that both displacement continuity and traction equilibrium are satisfied at the interface, maintaining the physical consistency of the coupled solution. The resulting system can be solved using standard sparse linear solvers, with the coupling naturally handled through the unified interface degrees of freedom.

\subsection{Implementation workflow}
\label{sec:implementation_workflow}

\begin{figure}[htbp]
    \centering
    \includegraphics[width=0.95\textwidth]{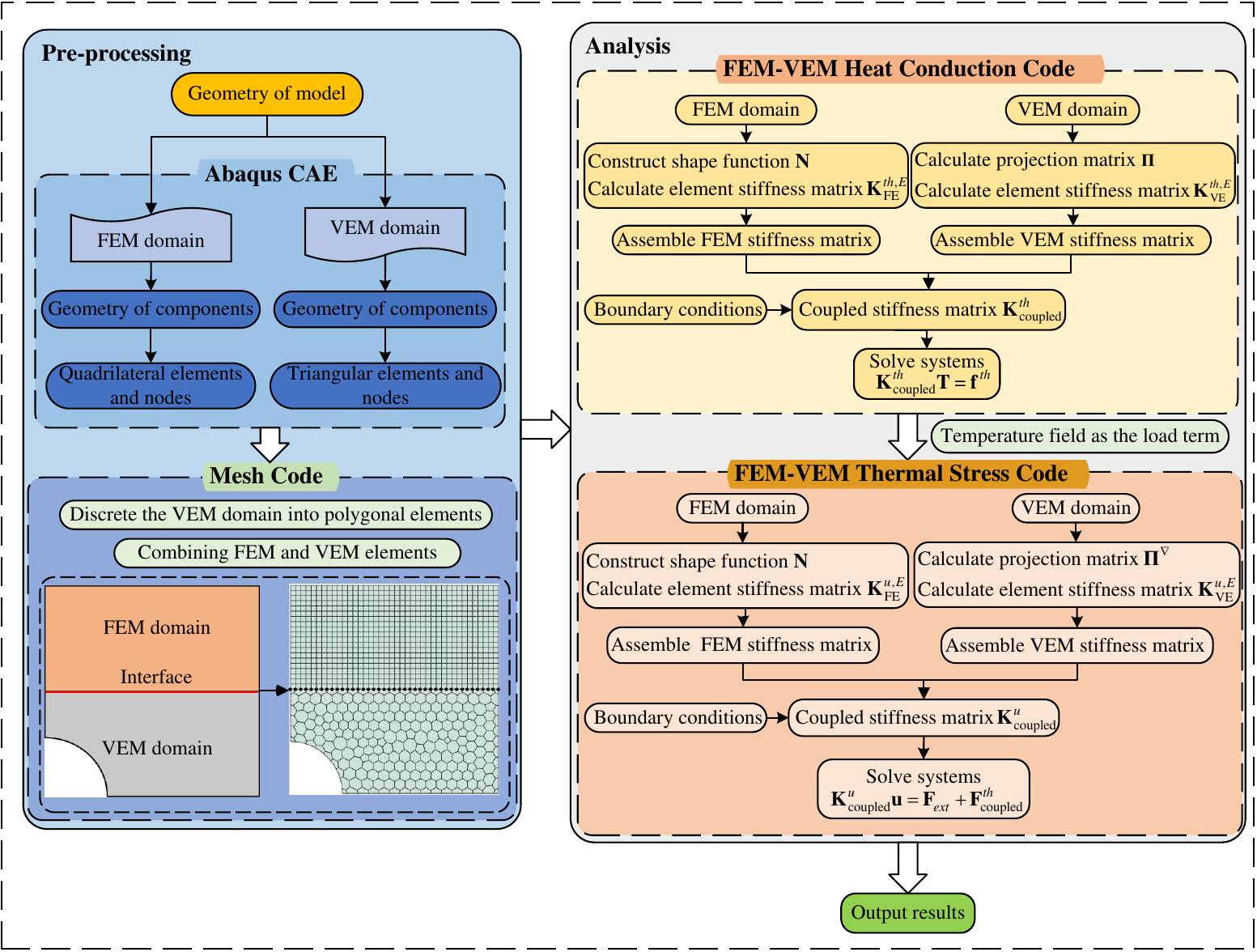}
    \caption{Thermomechanical coupling analysis based on FE-VE coupled method.}
    \label{fig:fem_vem_workflow}
\end{figure}

The comprehensive implementation workflow for the FE-VE coupled thermomechanical analysis is illustrated in Fig.~\ref{fig:fem_vem_workflow}, demonstrating the seamless integration of preprocessing, analysis, and post-processing phases. The computational procedure follows a systematic three-stage approach.

\subsubsection{Pre-processing stage}
\label{sec:preprocessing_stage}

The preprocessing phase establishes the computational framework through strategic domain partitioning and mesh generation. The geometry is initially defined and systematically partitioned into FE and VE domains based on geometric complexity considerations. The FE domain $\Omega^{\mathrm{FE}}$ employs structured quadrilateral elements to efficiently handle regular geometric regions, while the VE domain $\Omega^{\mathrm{VE}}$ utilizes arbitrary polygonal elements to accommodate complex geometric configurations.

The mesh generation process discretizes each domain using appropriate element types and establishes proper nodal connectivity at the FE-VE interface. This stage ensures geometric compatibility between the two domains and prepares the foundation for the subsequent coupling implementation detailed in Sec.~\ref{sec:interface_treatment}.

\subsubsection{Thermal analysis phase}
\label{sec:thermal_analysis_phase}

The thermal analysis phase implements the heat conduction formulation established in Sec.~\ref{sec:weak_formulation}. For the FE domain, standard shape functions are constructed and element thermal stiffness matrices $\mathbf{K}_{\mathrm{FE}}^{th,E}$ are computed using the procedures described in Sec.~\ref{sec:thermal_analysis_fem}. Simultaneously, the VE domain employs projection operators to calculate element stiffness matrices $\mathbf{K}_{\mathrm{VE}}^{th,E}$ following the VEM formulation.

The individual domain contributions are assembled into the coupled global thermal stiffness matrix $\mathbf{K}_{\mathrm{coupled}}^{th}$ according to the interface treatment methodology. After applying appropriate thermal boundary conditions, the coupled thermal system
\begin{equation}
\mathbf{K}_{\mathrm{coupled}}^{th}\mathbf{T} = \mathbf{f}_{\mathrm{coupled}}^{th}
\label{eq:coupled_thermal_system}
\end{equation}
is solved to obtain the temperature distribution $\mathbf{T}$ throughout the computational domain.

\subsubsection{Thermoelastic analysis phase}
\label{sec:thermoelastic_analysis_phase}

The mechanical analysis phase utilizes the computed temperature field as thermal loading input, implementing the thermomechanical coupling described in Sec.~\ref{sec:weak_formulation}. Following procedures analogous to the thermal analysis, shape functions are constructed and element mechanical stiffness matrices are computed for both domains: $\mathbf{K}_{\mathrm{FE}}^{u,E}$ for the finite element region and $\mathbf{K}_{\mathrm{VE}}^{u,E}$ for the virtual element region.

The thermal load vector $\mathbf{F}_{\mathrm{coupled}}^{th}$ is computed based on the temperature field obtained from the thermal analysis, incorporating thermal expansion effects. After assembling the coupled global mechanical stiffness matrix and applying appropriate displacement boundary conditions, the coupled mechanical system
\begin{equation}
\mathbf{K}_{\mathrm{coupled}}^u\mathbf{u} = \mathbf{F}_{\mathrm{ext}} + \mathbf{F}_{\mathrm{coupled}}^{th}
\label{eq:coupled_mechanical_system}
\end{equation}
is solved to obtain the displacement field $\mathbf{u}$ and subsequently the stress distribution throughout the domain.

For the electronic packaging structures considered in this work, material properties do not exhibit significant temperature dependence within the operating range, allowing us to employ this decoupled approach without loss of accuracy. By solving the thermal and mechanical problems sequentially, we can utilize specialized solvers optimized for each individual field, which enhances numerical stability and computational efficiency compared to solving the fully coupled system. 

\subsubsection{Post-processing and output}
\label{sec:postprocessing}

The final stage involves post-processing to extract engineering quantities of interest. Stress fields are computed using the established constitutive relationships, and von Mises equivalent stresses are evaluated for reliability assessment. The workflow concludes with the generation of visualization outputs including temperature distributions, displacement fields, and stress contours, providing insight into the thermomechanical behavior of the analyzed system.

\section{Numerical Examples}
\label{sec:numerical_examples}

\subsection{Square plate with circular hole}
\label{sec:square_plate_hole}

\begin{figure}[htbp]
    \centering
    \includegraphics[width=0.4\textwidth]{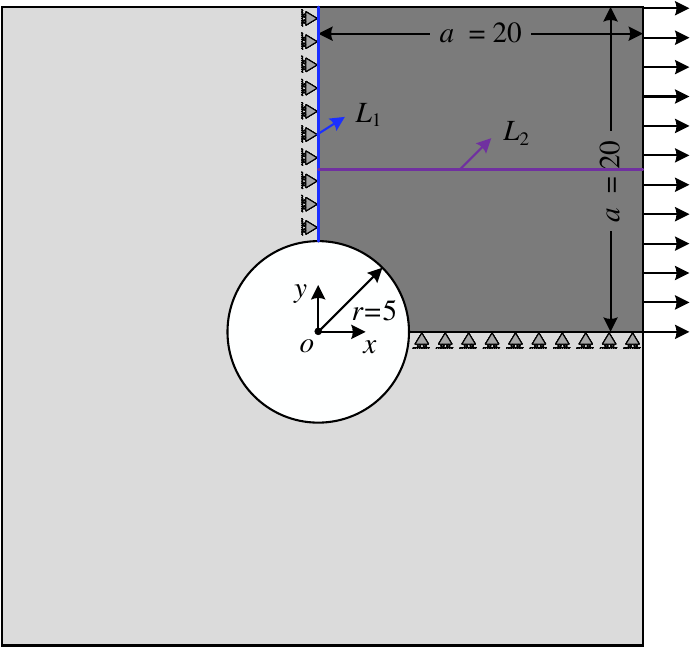}
    \caption{Geometric configuration and boundary conditions.}
    \label{fig:plate_geometry}
\end{figure}

To demonstrate the computational efficiency and convergence characteristics of the proposed algorithm, a square plate with a circular hole subjected to uniform loading is presented as the first numerical example. Due to the symmetric nature of the model, only one quarter of the domain is analyzed by exploiting symmetry conditions, as illustrated in Fig.~\ref{fig:plate_geometry}. 
The model parameters are specified as follows: circular hole radius $r = 5$~mm, plate dimension $a = 20$~mm, applied load intensity $F = 5$~MPa, elastic modulus $E = 10$~MPa, and Poisson's ratio $\nu = 0.3$. The problem is solved under plane stress conditions.

\begin{figure}[htbp]
    \centering
    \includegraphics[width=0.8\textwidth]{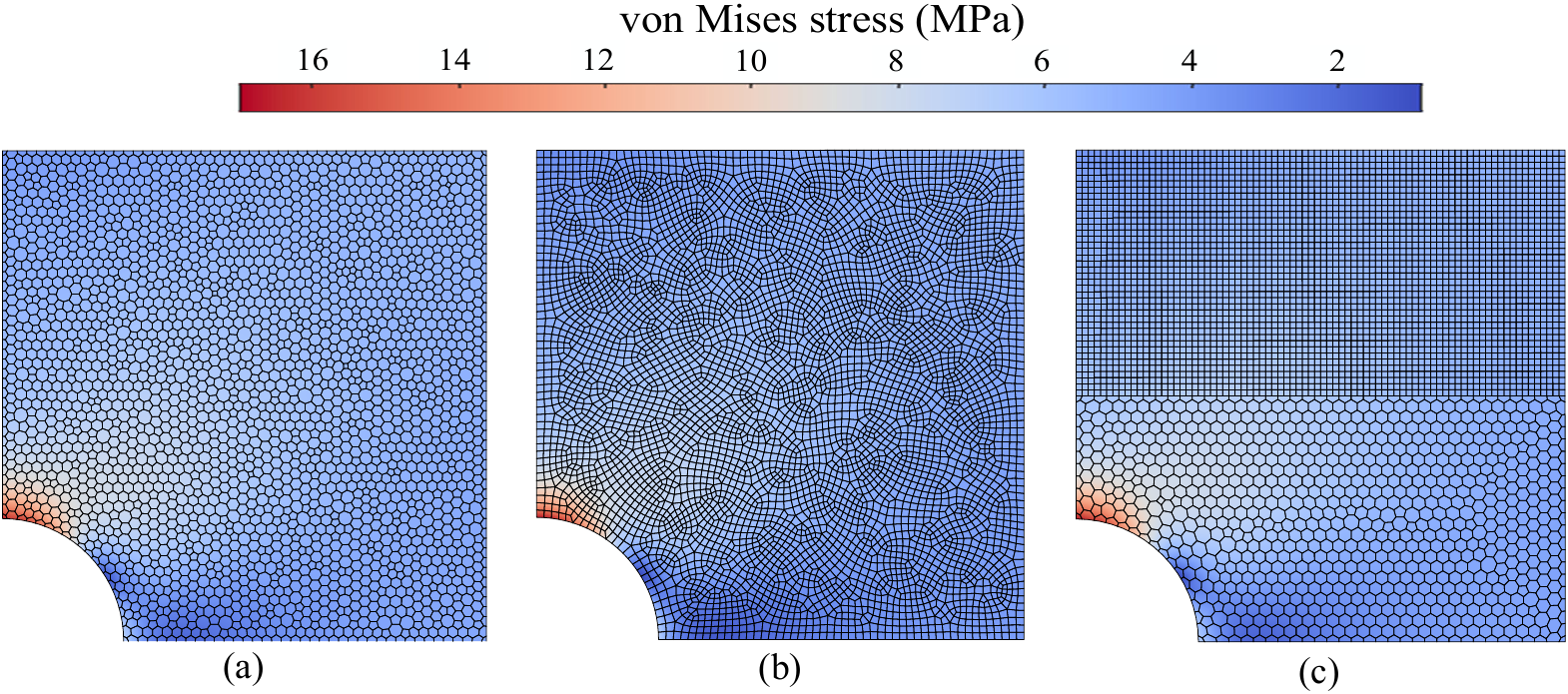}
    \caption{von Mises stress distribution comparison for square plate with circular hole: (a) VEM solution, (b) FEM solution, (c) FE-VE coupling solution.}
    \label{fig:stress_contours}
\end{figure}

Three different numerical approaches were employed to solve this problem: VEM, FEM, and the proposed FE-VE coupling method. The domain is partitioned such that the regular regions employ finite elements while the area near the circular hole utilizes virtual elements to handle the geometric complexity. The resulting von Mises stress distribution contours are presented in Fig.~\ref{fig:stress_contours}. The stress contours demonstrate excellent agreement among all three methods, particularly around the critical stress concentration region near the circular hole. Fig.~\ref{fig:stress_contours}a shows the VEM solution with its characteristic polygonal mesh that naturally adapts to the curved boundary geometry. Fig.~\ref{fig:stress_contours}b presents the FEM solution employing a structured quadrilateral mesh, while Fig.~\ref{fig:stress_contours}c illustrates the FE-VE coupling approach that effectively combines the advantages of both methods. The stress magnitude distributions and overall stress patterns are highly consistent across all three approaches, confirming the accuracy and reliability of the proposed coupling methodology.
\begin{figure}[htbp]
    \centering
    \includegraphics[width=0.5\textwidth]{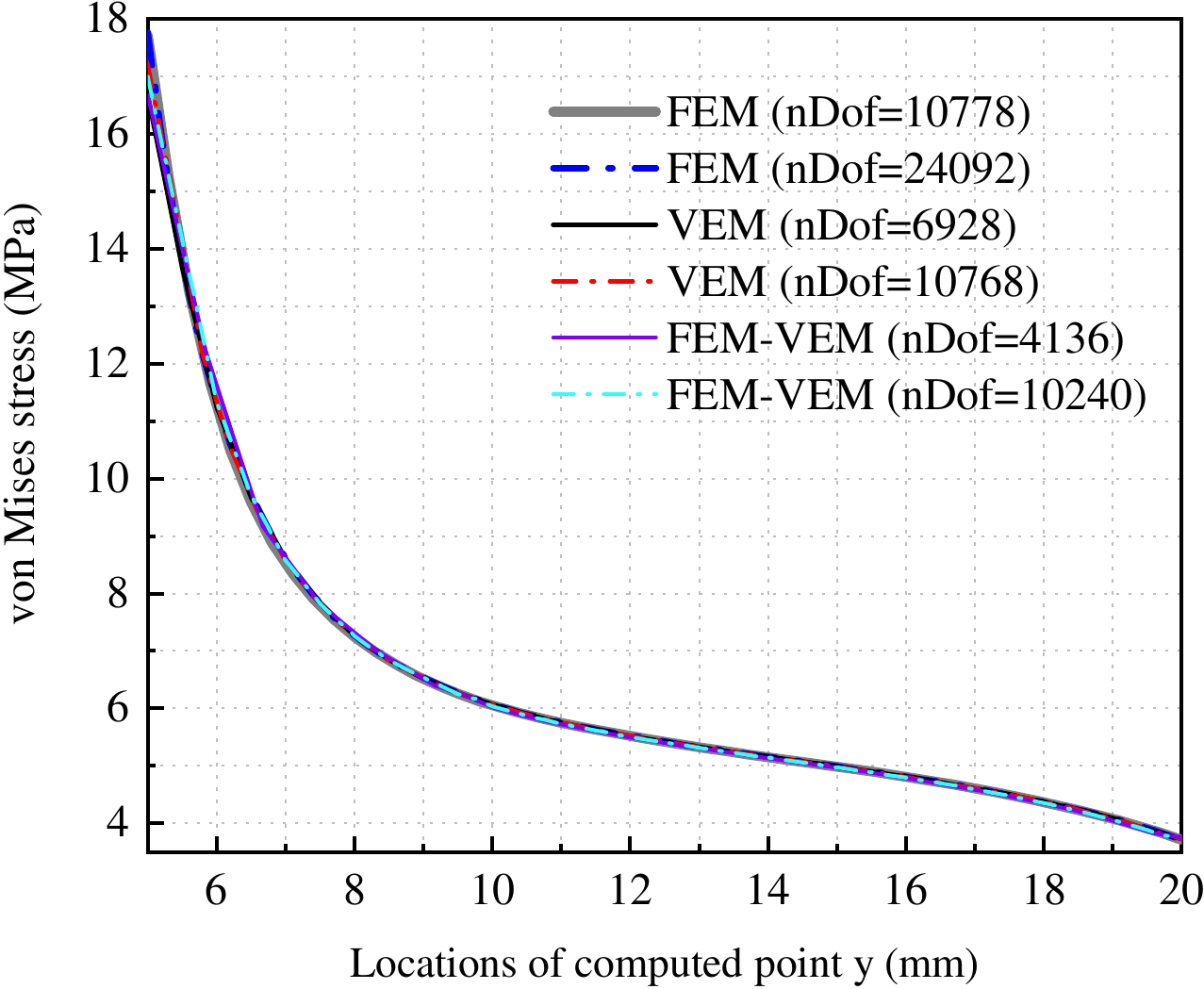}
    \caption{Stress distribution along line $L_1$ for different methods with varying nDofs.}
    \label{fig:stress_distribution_x0}
\end{figure}

The stress distributions along the line $L_1$ (given in Fig.~\ref{fig:plate_geometry}) were extracted for FEM, VEM, and FE-VE method under different Dofs, as illustrated in Fig.~\ref{fig:stress_distribution_x0}. The results demonstrate excellent agreement among all three approaches, confirming the accuracy and consistency of the proposed coupling methodology. Additionally, the computational times for the three methods were compared, as presented in Table~\ref{tab:computational_time_comparison}. From Table~\ref{tab:computational_time_comparison}, it is evident that computational time increases with the number of degrees of freedom (nDofs) for all three methods, which is expected behavior. When comparing methods with similar nDof (approximately 10,700), FEM outperforms VEM with computational times of 0.61 seconds and 0.78 seconds, respectively. This performance gap arises primarily from the complex internal projection operations and stabilization terms required by VEM, which lead to significantly higher computational cost per element. Notably, the coupled FE-VE method achieves the highest efficiency, requiring only 0.47 seconds. By leveraging VEM’s ability to handle complex geometries with fewer elements and FEM’s computational efficiency in regular regions, this hybrid strategy effectively circumvents the high cost of pure VEM in complex domains and the need for excessive mesh refinement with pure FEM in simple regions, resulting in superior overall performance. 

\begin{table}[htbp]
    \scriptsize
    \centering
    \caption{Computational time comparison for VEM, FEM, and FE-VE coupled method.}
    \label{tab:computational_time_comparison}
    \begin{tabular}{cccc}
        \toprule
        \textbf{Method} & \textbf{Number of elements} & \textbf{nDof} & \textbf{ Time (s)} \\
        \midrule
        \multirow{2}{*}{VEM} & 1655 & 6928 & 0.33\\
                             & 2595 & 10768 & 0.778\\
        \midrule
        \multirow{2}{*}{FEM} & 5263 & 10778 & 0.61\\
                             & 11850 & 24092 & 2.19\\
        \midrule
        \multirow{2}{*}{FEM-VEM} & 1611 & 4136 & 0.16\\
                                 & 4083 & 10240 & 0.47\\
        \bottomrule
    \end{tabular}
\end{table}

To investigate the convergence characteristics of the FE-VE algorithm, calculations were performed using meshes of varying densities. The coupling interface along the line $L_2$ (shown in Fig.~\ref{fig:plate_geometry}) was selected as the computational output and comparison path for detailed analysis. The mean relative error (MRE) is defined as
\begin{equation}
\text{MRE} = \frac{1}{N} \sum_{n=1}^{N}\left|\frac{\boldsymbol{\sigma}_{\text{num}} - \boldsymbol{\sigma}_{\text{ref}}}{\boldsymbol{\sigma}_{\text{ref}}}\right|,
\label{eq:mre_definition}
\end{equation}
where $N$ represents the total number of nodes in the calculation, $\boldsymbol{\sigma}_{\text{num}}$ denotes the numerical solution of nodal stress, and $\boldsymbol{\sigma}_{\text{ref}}$ represents the reference solution of nodal stress obtained from mesh refinement (nDof =334,254) in Abaqus.
\begin{figure}[htbp]
    \centering
    \includegraphics[width=0.55\textwidth]{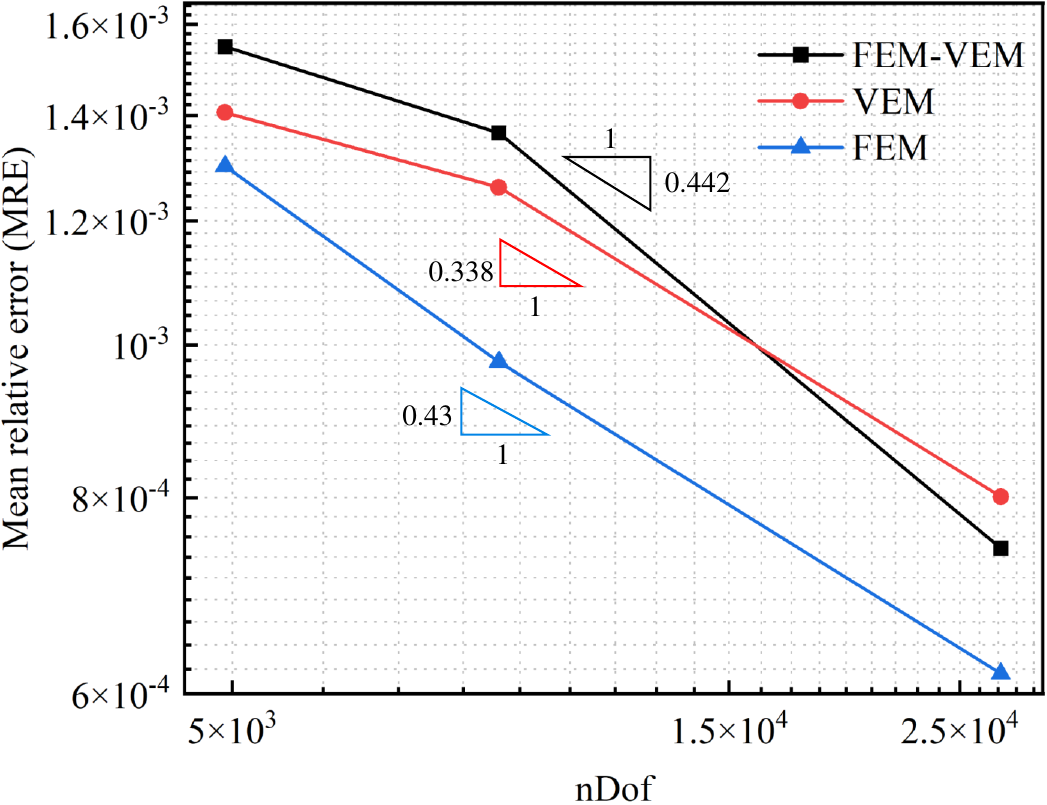}
    \caption{Convergence analysis: MRE comparison of different numerical algorithm.}
    \label{fig:convergence_analysis}
\end{figure}

Fig.~\ref{fig:convergence_analysis} presents the MRE convergence comparison for FEM, VEM, and FE-VE coupling method as the nDof increases from approximately 5,000 to 25,000. The convergence analysis reveals distinct performance characteristics for each method, with convergence rates indicated by the slope annotations in the figure.
The convergence rates indicate that the FE-VE coupling method actually achieves slightly better convergence behavior than pure FEM, while VEM shows the most gradual convergence. All three methods demonstrate stable convergence trends, confirming the reliability of the proposed coupling strategy. The FE-VE approach maintains good accuracy while providing the additional advantage of geometric flexibility for complex domain geometries.


\subsection{Thermoelastic analysis of thick-walled cylindrical structure}

\begin{figure}[htbp]
    \centering
    \subfloat[]{\centering
        \includegraphics[width=0.35\textwidth]{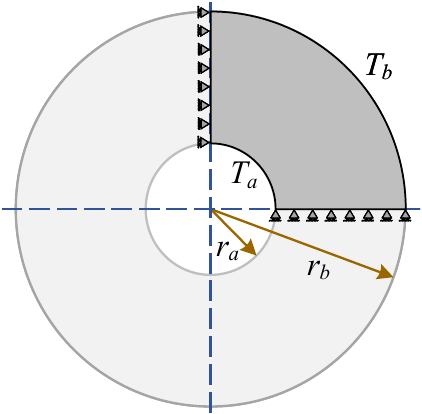}
        \label{fig: thick-walled cylindrical model}
    }
    \hfill
    \subfloat[]{\centering
     \includegraphics[width=0.35\textwidth]{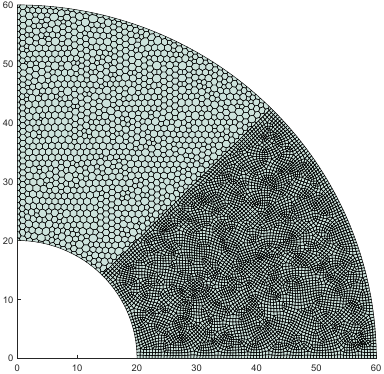}
        \label{fig:thick-walled cylindrical mesh}
    }
    \caption{Thick-walled cylinder configuration: (a) Full geometry with inner radius $r_a$, outer radius $r_b$, and applied temperatures $T_a$ and $T_b$; 
    (b) Mesh employed in the FE-VE coupled method}
    \label{fig:cylinder_model}
\end{figure}

We analyze a homogeneous thick-walled cylinder under thermal loading to validate the FE-VE coupling method. The analytical solution for this problem serves as a benchmark for accuracy assessment. Fig.~\ref{fig:cylinder_model} illustrates the geometry and boundary conditions of this verification case. 

The cylinder has inner radius $r_a = 20$ mm and outer radius $r_b = 60$ mm. Material properties include Young's modulus $E = 460000$ MPa, Poisson's ratio $\nu = 0.3$, thermal expansion coefficient $\alpha = 7.4 \times 10^{-6}$ K$^{-1}$, and thermal conductivity $\lambda = 20$ W/(m$\cdot$K). The thermal boundary conditions consist of $T_a = 0$ K at the inner surface and $T_b = 500$ K at the outer surface under steady-state conditions. Due to symmetry, only one quarter of the cylinder is modeled to optimize computational efficiency, as shown in Fig.~\ref{fig:cylinder_model}.

Given the material parameters and boundary conditions, the analytical temperature field is
\begin{equation}\label{eq:temperature_analytical}
T = T_a + \frac{T_b - T_a}{\ln(r_b/r_a)}\ln\left(\frac{r}{r_a}\right),
\end{equation}
where $r$ represents the radial coordinate. 

The method's convergence characteristics are assessed using the root mean square (RMS) $L^2$ error as the accuracy metric. This normalized error measure is defined as
\begin{equation}\label{eq:rms_error}
\varepsilon_{\text{rms}} = \frac{1}{\max(|X_e|)}\sqrt{\frac{1}{N}\sum_{i=1}^N |X - X_e|^2}
\end{equation}
where $N$ is the number of evaluation points, and $X_e$ and $X$ denote the analytical and numerical solutions, respectively.

\begin{figure}[htbp]
    \centering
    \includegraphics[width=0.5\textwidth]{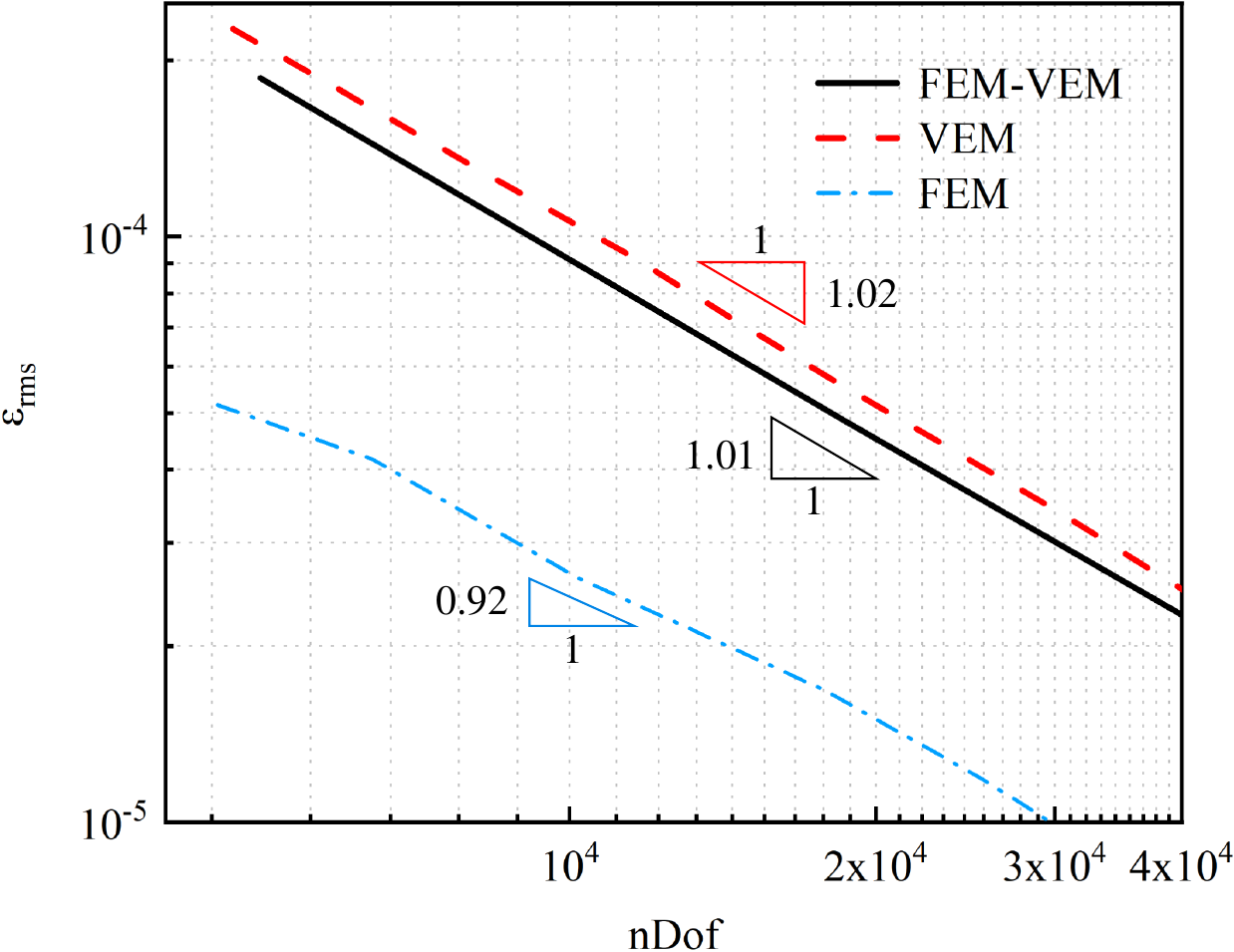}
    \caption{Temperature accuracy versus mesh refinement for thick-walled cylinder}
    \label{fig:EX2_temperature_convergence}
\end{figure}

Fig.~\ref{fig:EX2_temperature_convergence} presents the convergence analysis comparing FE-VE coupling, pure VEM, and pure FEM methods, showing RMS errors for temperature as functions of degrees of freedom. The convergence rates, indicated by the slope annotations in the figure, demonstrate that the VEM achieves a convergence rate of 1.02, while the FE-VE coupling method exhibits a rate of 1.01, both outperforming pure FEM with a rate of 0.92. The proposed FE-VE coupling method maintains nearly optimal convergence characteristics comparable to pure VEM, while providing the additional advantage of geometric flexibility through strategic domain partitioning. All three methods show consistent error reduction with mesh refinement, confirming the stability and reliability of the coupling approach.

\subsection{Simplified sintered silver model}
\label{sec:sintered_silver_model}

\begin{figure}[htbp]
    \centering
    \includegraphics[width=0.48\textwidth]{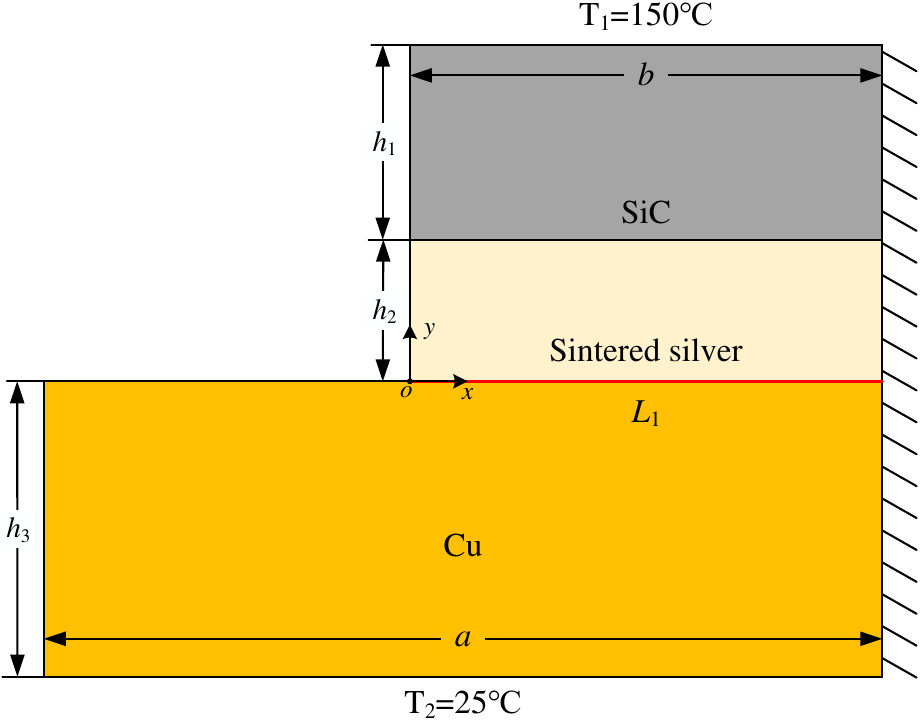}
    \caption{Geometric configuration and thermal boundary conditions.}
    \label{fig:sintered_silver_model}
\end{figure}

The third numerical example investigates a sintered silver interconnect structure commonly employed in high-power electronic packaging applications. As illustrated in Fig.~\ref{fig:sintered_silver_model}, the model adopts a typical `sandwich' configuration consisting of three distinct layers: a SiC chip (top), sintered silver interconnect layer (middle), and copper substrate (bottom). This multilayer structure is representative of modern power electronic modules where effective thermal management and mechanical reliability are critical design considerations~\cite{yu2021thermal}. The geometric dimensions are specified as follows: the SiC chip has width $b = 1.8$~mm and thickness $h_1 = 0.5$~mm; the sintered silver layer maintains the same width as the SiC chip with thickness $h_2 = 0.3$~mm; and the copper substrate has width $a = 3$~mm and thickness $h_3 = 0.8$~mm. The material properties for each component are detailed in Table~\ref{tab:material_properties}, which reveal significant differences in thermal and mechanical characteristics that contribute to complex thermomechanical stress development during thermal loading.

\begin{table}[htbp]
    \scriptsize
    \centering
    \caption{Material properties of sintered silver interconnect model components.}
    \label{tab:material_properties}
    \begin{tabular}{lcccc}
    \toprule
    Material & \begin{tabular}{c} Young's modulus \\ (GPa) \end{tabular} & Poisson's ratio & \begin{tabular}{c} Thermal conductivity \\ (W/m$\cdot$K) \end{tabular} & \begin{tabular}{c} CTE $\alpha$ \\ ( $^\circ$C$^{-1}$) \end{tabular} \\
    \midrule
    SiC chip & 410 & 0.14 & 370 & 4.5 $\times 10^{-6}$ \\
    Sintered silver & 12.9 & 0.38 & 278 & 19.0 $\times 10^{-6}$ \\
    Copper substrate & 110 & 0.38 & 400 & 16.5 $\times 10^{-6}$ \\
    \bottomrule
    \end{tabular}
\end{table}

The top surface is subjected to a temperature of $T_1 = 150$ $^\circ$C representing heat generation from the semiconductor device, while the bottom surface is maintained at ambient temperature $T_2 = 25$ $^\circ$C. All lateral surfaces are treated as thermally insulated. Mechanically, the right edge is fully constrained ($u_x = u_y = 0$) and other boundaries of the model are set as free boundaries.To demonstrate the capabilities of the FE-VE coupling approach, the upper layers (SiC chip and sintered silver) are discretized using virtual elements, while the copper substrate employs conventional finite elements. 

\begin{figure}[htbp]
    \centering
    \includegraphics[width=0.6\textwidth]{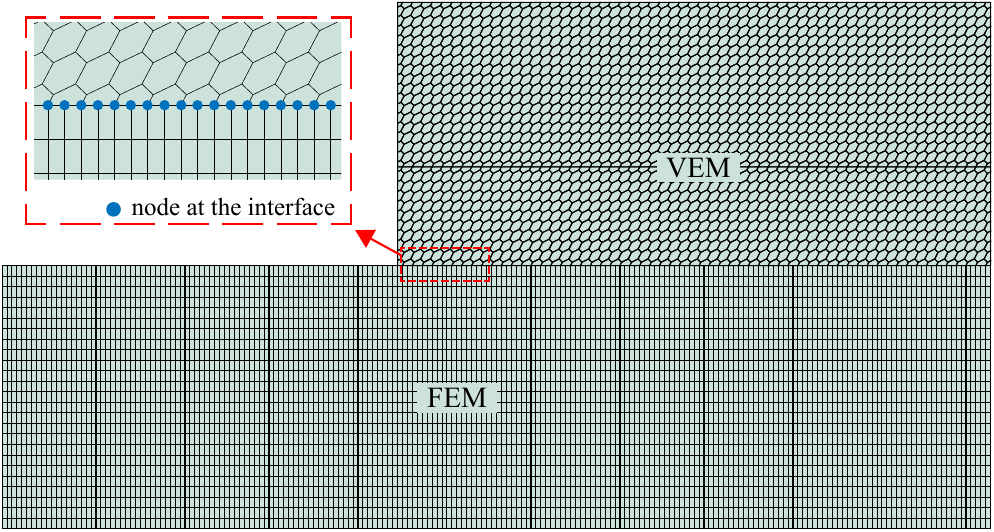}
    \caption{Hybrid mesh discretization strategy for FE-VE coupling analysis}.
    \label{fig:mesh_discretization}
\end{figure}

For the numerical computation, the model is discretized using Abaqus preprocessing capabilities combined with custom MATLAB code, where the virtual element domain (upper layers) employs polygonal elements to leverage VEM's geometric flexibility, while the finite element domain (copper substrate) utilizes conventional four-node quadrilateral elements for computational efficiency. 
A critical aspect of the implementation is ensuring nodal compatibility at the coupling interface, where nodes from both domains must be properly matched to maintain solution continuity, as illustrated in Fig.~\ref{fig:mesh_discretization}. The figure demonstrates the distinct meshing approaches for each domain and the interface treatment, where polygonal elements transition to structured quadrilateral elements while maintaining nodal correspondence. 
This interface treatment ensures that displacement and stress continuity conditions are satisfied across the domain boundary, which is essential for the accuracy and stability of the FE-VE coupling scheme.

\begin{figure}[htbp]
    \centering
    \subfloat[FE-VE coupling solution]{\centering
        \includegraphics[width=0.46\textwidth]{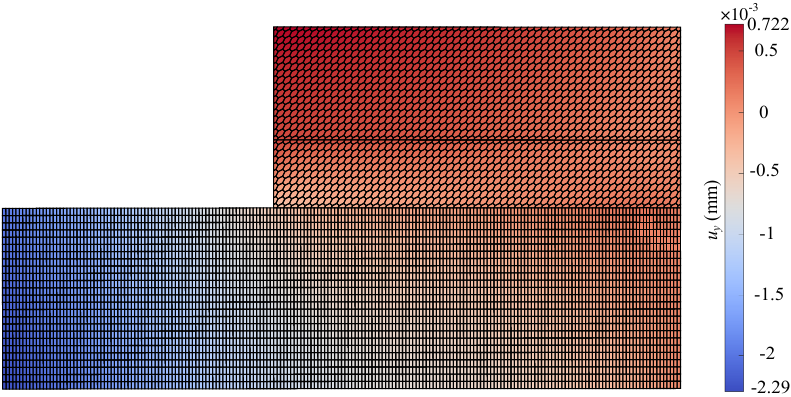}
        \label{fig:displacement_coupling}
    }
    \hfill
    \subfloat[FEM reference solution]{\centering
        \includegraphics[width=0.46\textwidth]{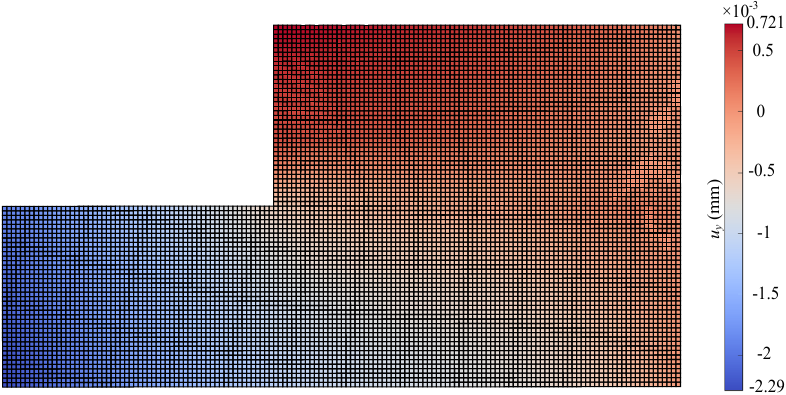}
        \label{fig:displacement_fem}
    }
    \caption{Comparison of displacement ($u_y$) distributions.}
    \label{fig:displacement_comparison}
\end{figure}

\begin{figure}[htbp]
    \centering
    \subfloat[Stress on sintered silver side of interface]{\centering
        \includegraphics[width=0.46\textwidth]{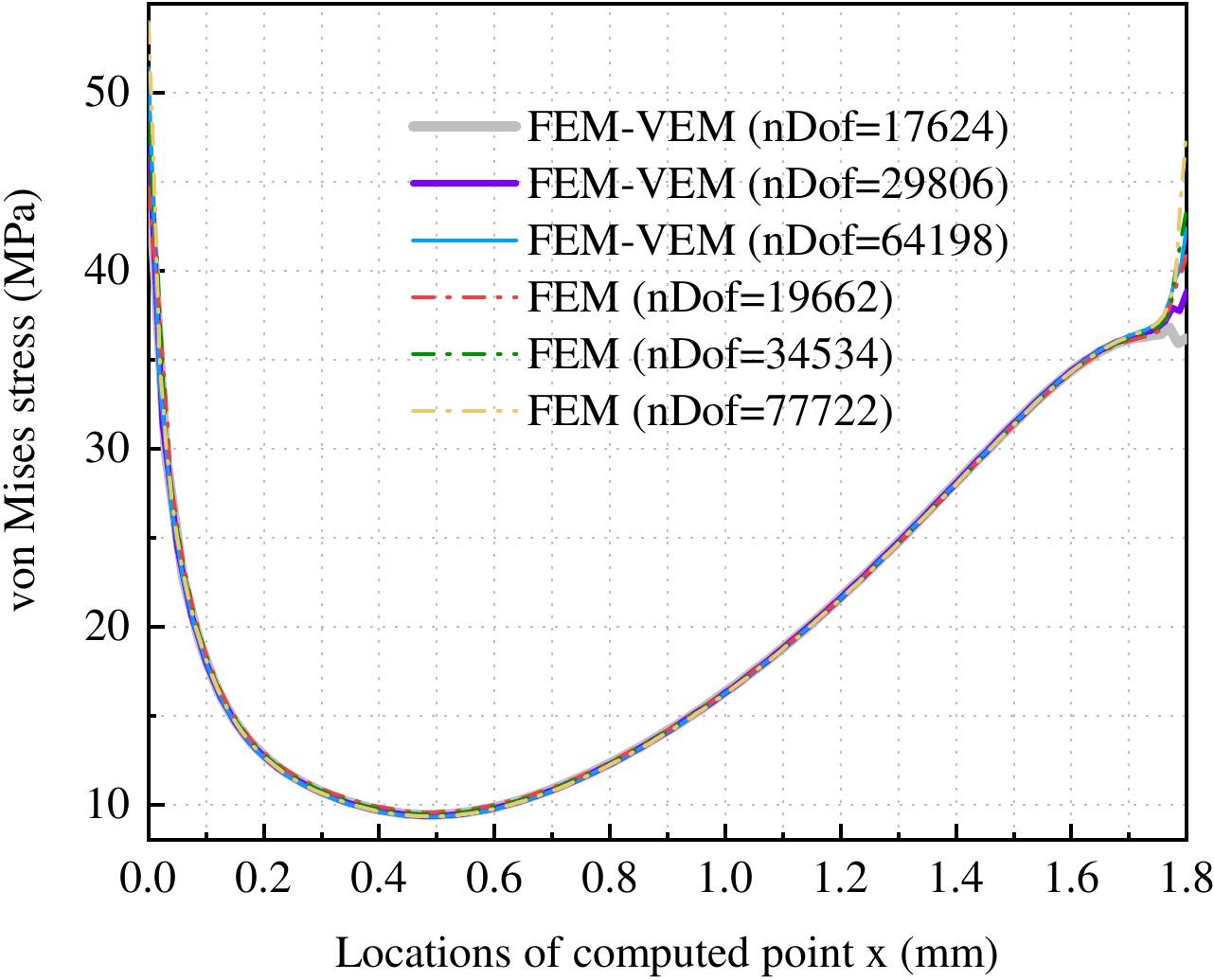}
        \label{fig:stress_silver}
    }
    \hfill
    \subfloat[Stress on copper substrate side of interface]{\centering
        \includegraphics[width=0.46\textwidth]{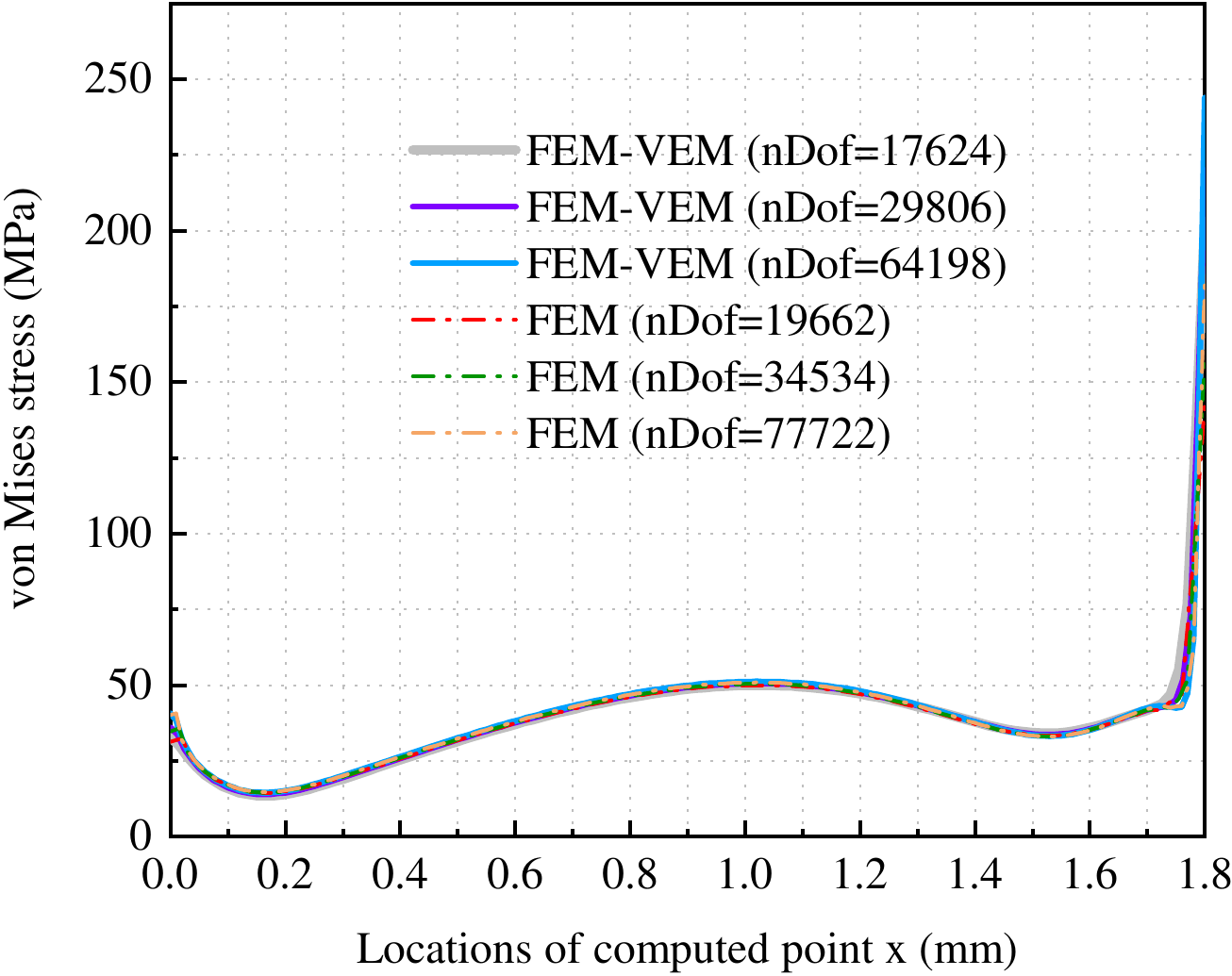}
        \label{fig:stress_copper}
    }
    \caption{Stress distributions along both sides of coupling interface $L_1$ (shown in Fig.~\ref{fig:sintered_silver_model}).}
    \label{fig:stress_interface_comparison}
\end{figure}

The computational accuracy of the FE-VE coupling algorithm is validated through comparison with a refined finite element solution. Fig.~\ref{fig:displacement_comparison} presents the overall displacement ($u_y$) contours, demonstrating that both methods yield very similar results across the entire domain. The displacement patterns exhibit expected thermal expansion behavior, with upper layers (SiC and sintered silver) experiencing upward displacement due to thermal loading while the constrained right boundary maintains zero displacement. 

To provide detailed accuracy assessment, von Mises stress values at each node along coupling interface line $L_1$ (given in Fig.~\ref{fig:sintered_silver_model}) were extracted for both sintered silver and copper substrate layers, as shown in Fig.~\ref{fig:stress_interface_comparison}. The comparison reveals good agreement between the FE-VE coupling algorithm and the reference solution across different mesh refinement levels, with stress distributions showing consistent convergence trends. 
The stress profiles demonstrate expected material-dependent behavior, with peak stresses reaching approximately 50 MPa at the sintered silver interface and up to 260 MPa at the copper substrate interface. The significantly higher stresses in the copper substrate can be attributed to its higher elastic modulus (110 GPa vs 12.9 GPa for sintered silver), which results in greater stress concentrations under thermal loading conditions. These results confirm the capability of the proposed coupling methodology to accurately capture critical stress distributions in multi-material electronic packaging applications.

\begin{figure}[htbp]
    \centering
    \includegraphics[width=0.6\textwidth]{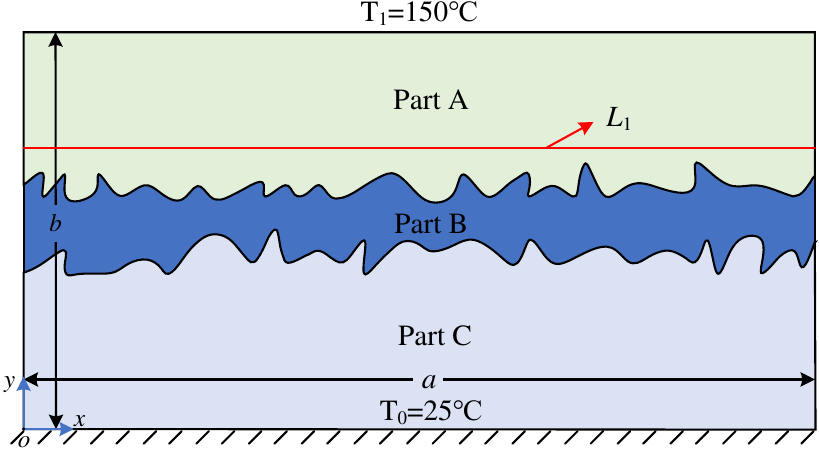}
    \caption{Simplified model with complex-shaped interface featuring three material regions (Parts A, B, C) and irregular interface geometry.}
    \label{fig:Complex_interface_model}
\end{figure}

Previous analyses examined the thermomechanical behavior of sintered silver structures with flat interfaces. However, in practical applications, interface morphology often exhibits complex geometric characteristics such as random roughness rather than ideal smoothness. This geometric complexity can induce significant local stress concentrations and alter load transfer paths within the structure, thereby critically influencing overall mechanical performance. To more accurately predict the mechanical response of such structures, we establish and analyze a simplified model with a randomly complex-shaped interface, whose geometric configuration is shown in Fig.~\ref{fig:Complex_interface_model}.

The model has geometric parameters $a = 40$ mm and $b = 20$ mm and consists of three distinct parts with material properties listed in Table~\ref{tab:Complex_interface_materials}. Temperature boundary conditions are applied at the top surface ($T_1 = 150$ $^\circ$C) and bottom surface ($T_0 = 25$ $^\circ$C), with all other surfaces thermally insulated. The bottom surface is subjected to a fully fixed mechanical constraint and other boundaries of the model are set as free boundaries. 

\begin{table}[htbp]
    \scriptsize
    \centering
    \caption{Material properties of Complex interface model components.}
    \label{tab:Complex_interface_materials}
    \begin{tabular}{lcccc}
    \toprule
    Material & \begin{tabular}{c} Young's modulus \\ (GPa) \end{tabular} & Poisson's ratio & \begin{tabular}{c} Thermal conductivity \\ (W/m$\cdot$K) \end{tabular} & \begin{tabular}{c} CTE $\alpha$ \\ ( $^\circ$C$^{-1}$) \end{tabular} \\
    \midrule
    Part A& 69.5& 0.33& 234& 23.2 $\times 10^{-6}$\\
    Part B& 1.2& 0.47& 3.2& 92 $\times 10^{-6}$\\
    Part C& 187& 0.28& 147& 2.6 $\times 10^{-6}$\\
    \bottomrule
    \end{tabular}
\end{table}

\begin{figure}[htbp]
    \centering
    \includegraphics[width=0.8\textwidth]{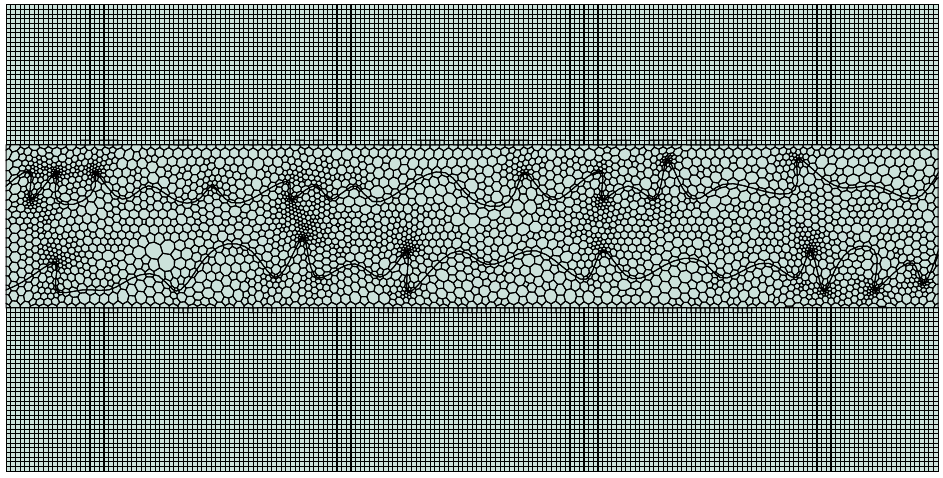}
    \caption{Hybrid mesh discretization strategy for the model with complex interface: quadrilateral elements in regular regions and polygonal elements near irregular interface}
    \label{fig:Complex_interface_mesh}
\end{figure}

To demonstrate the advantages of the FE-VE method in handling complex geometries, a coupled analysis is performed where the virtual element method discretizes regions around the complex interface while the FEM is applied to regular regions. The finite element domain employs four-node quadrilateral elements, whereas the virtual element domain utilizes arbitrary polygonal elements that flexibly conform to the irregular interface geometry. The corresponding mesh discretization is illustrated in Fig.~\ref{fig:Complex_interface_mesh}.

\begin{figure}[htbp]
    \centering
    \subfloat[Displacement: FE-VE coupling solution]{\centering
    \includegraphics[width=0.46\textwidth]{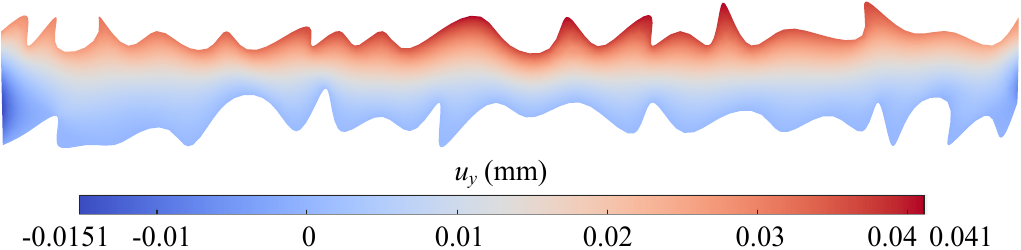}
    \label{fig:disp_coupling}
    }
    \subfloat[Displacement: FEM solution]{\centering
    \includegraphics[width=0.46\textwidth]{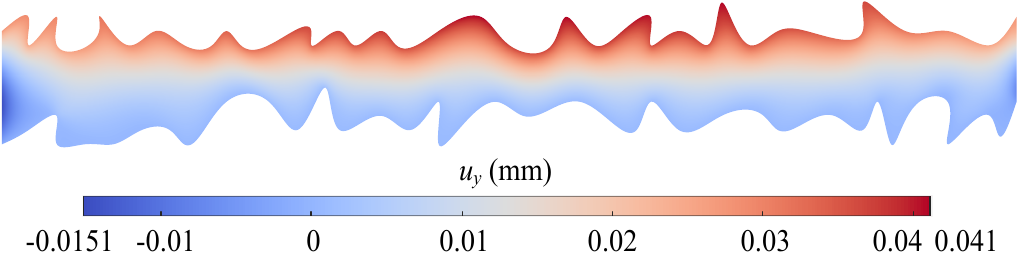}
    \label{fig:disp_FEM}
    }\\
    \subfloat[Stress: FE-VE coupling solution]{\centering
        \includegraphics[width=0.46\textwidth]{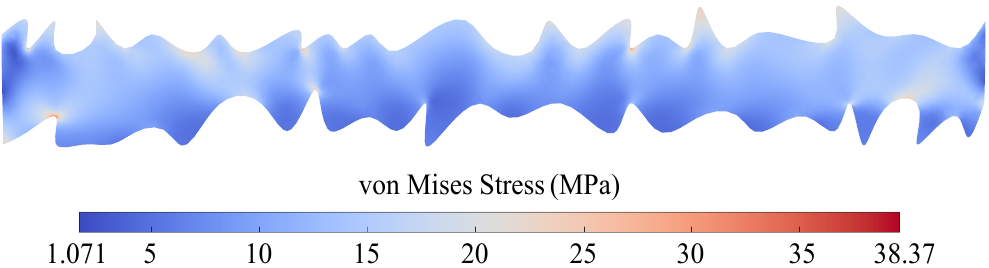}
        \label{fig:stress_coupling}
    }
    \subfloat[Stress: FEM reference solution]{\centering
        \includegraphics[width=0.46\textwidth]{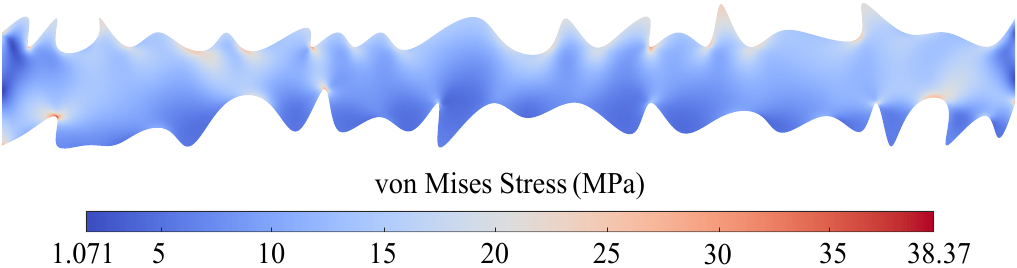}
        \label{fig:stress_fem}
    }
    \caption{Comparison of displacement and von Mises stress distributions in Part B.}
    \label{fig:disp_stress_comparison}
\end{figure}

To validate the computational accuracy of the FE-VE method, results are compared with a refined pure finite element analysis. Fig.~\ref{fig:disp_stress_comparison} presents the displacement and von Mises stress distributions in Part B near the complex interface, revealing excellent agreement between both methods. Both approaches accurately capture the stress concentration phenomena at the irregular interface, thereby verifying the effectiveness and reliability of the FE-VE method for modeling structures with complex geometries.

\begin{figure}[htbp]
    \centering
    \subfloat[Horizontal displacement ($u_x$)]{\centering
        \includegraphics[width=0.32\textwidth]{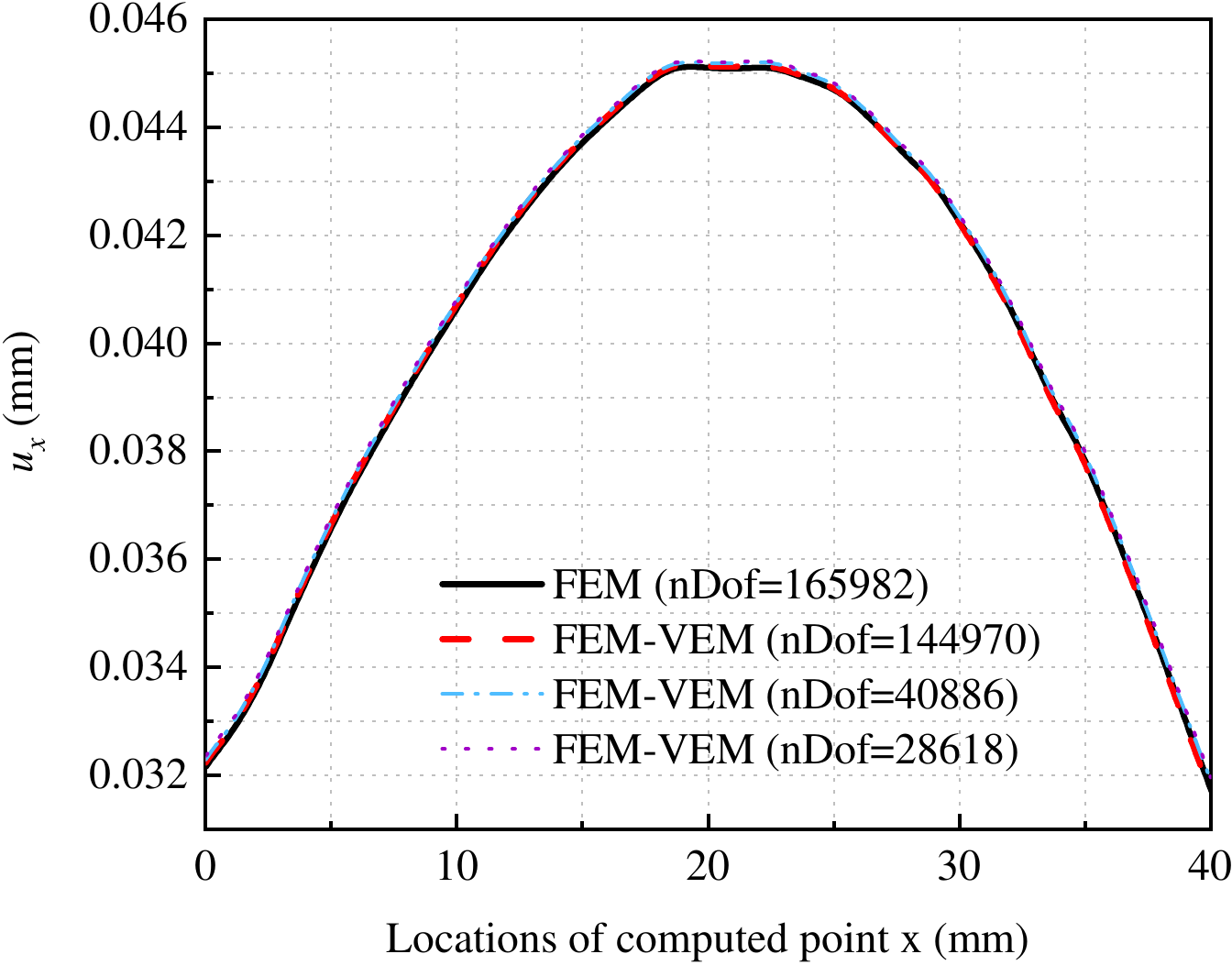}
        \label{fig:Complex_interface_ux_l1}
    }
    \subfloat[Vertical displacement ($u_y$)]{\centering
        \includegraphics[width=0.32\textwidth]{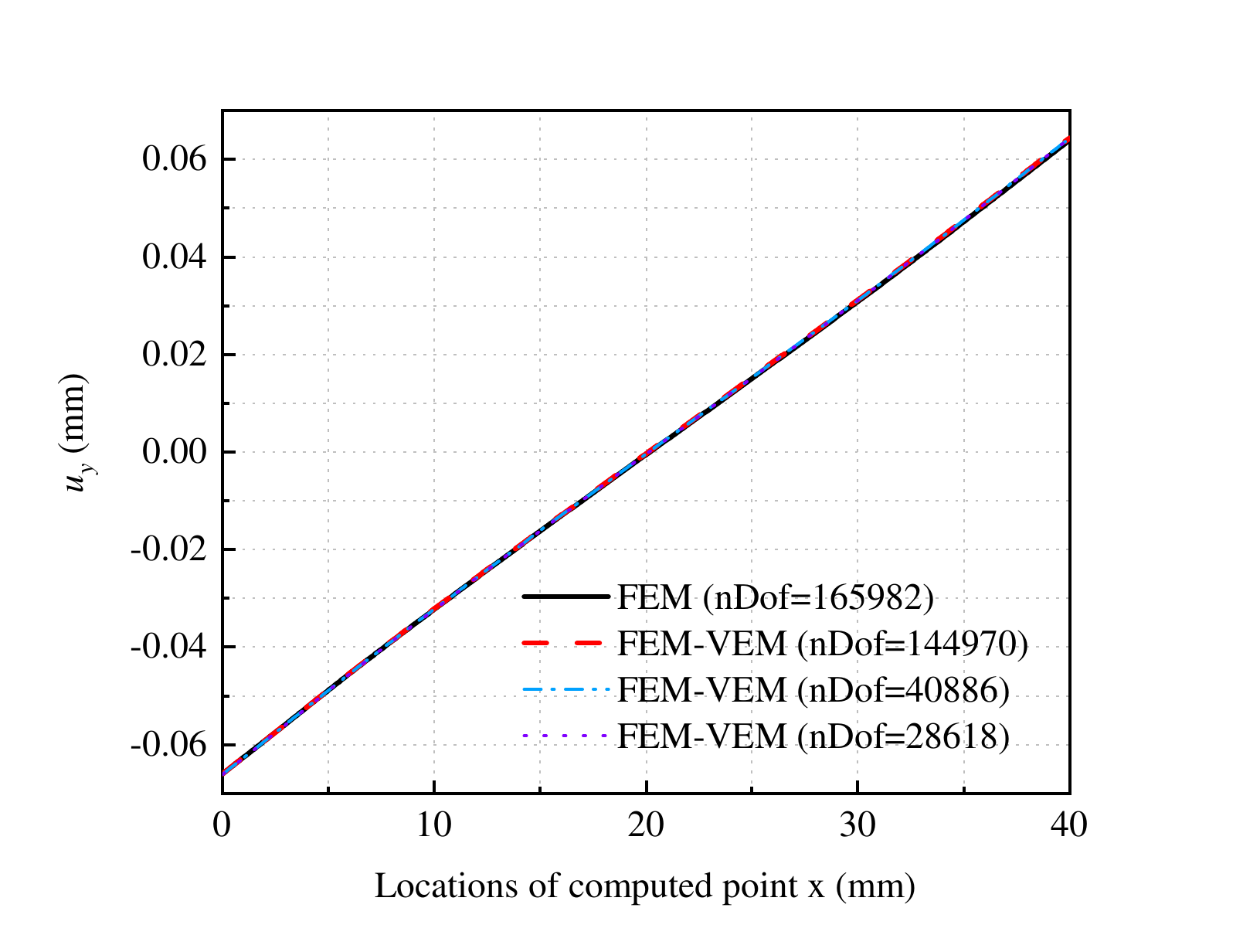}
        \label{fig:Complex_interface_uy_l1}
    }
    \subfloat[von Mises stress distribution]{\centering
        \includegraphics[width=0.32\textwidth]{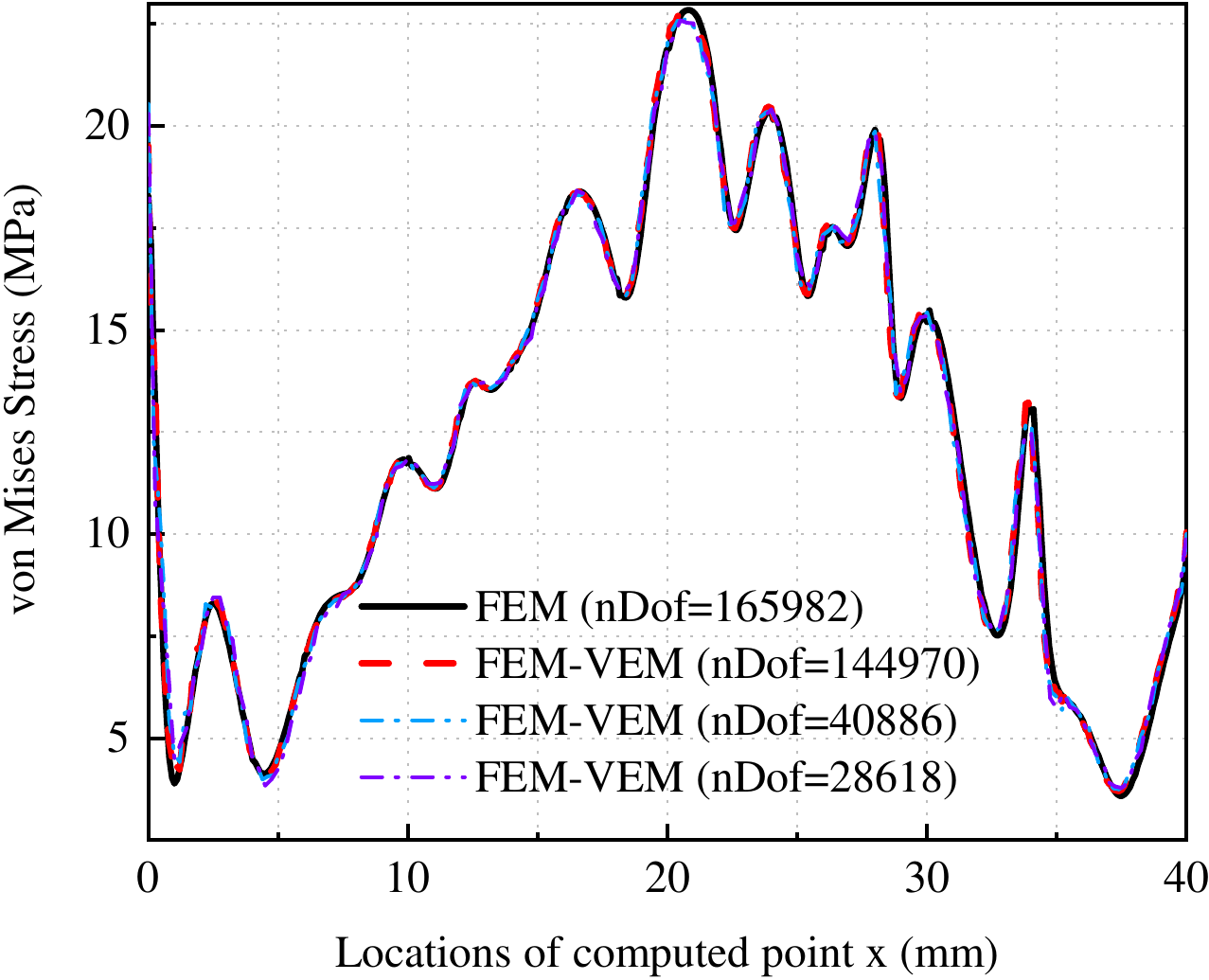}
        \label{fig:Complex_interface_stress_l1}
    }
    \caption{Displacement and stress distributions along coupling interface line $L_1$.}
    \label{fig:Complex_interface_results_l1}
\end{figure}

To quantitatively assess the coupling accuracy, displacement and stress distributions along the interface line $L_1$ (shown in Fig.~\ref{fig:Complex_interface_model}) are extracted and compared with the refined FEM reference solution (nDof = 165,982). Fig.~\ref{fig:Complex_interface_ux_l1} shows that the horizontal displacement peaks near the midpoint at approximately 0.045 mm. Fig.~\ref{fig:Complex_interface_uy_l1} reveals a nearly linear vertical displacement distribution ranging from –0.065 mm to 0.065 mm. Fig.~\ref{fig:Complex_interface_stress_l1} demonstrates that the von Mises stress exhibits significant fluctuations between 4 MPa and 23 MPa due to the complex interface geometry. 

The FE-VE solutions show excellent agreement with the FEM reference across all three mesh densities. As the mesh is refined, the coupled method results converge systematically toward the reference solution, demonstrating good accuracy and convergence characteristics for problems involving complex geometries.

\subsection{FC-BGA Model}
\label{sec:fc_bga_model}

\begin{figure}[htbp]
    \centering
    \includegraphics[width=0.8\textwidth]{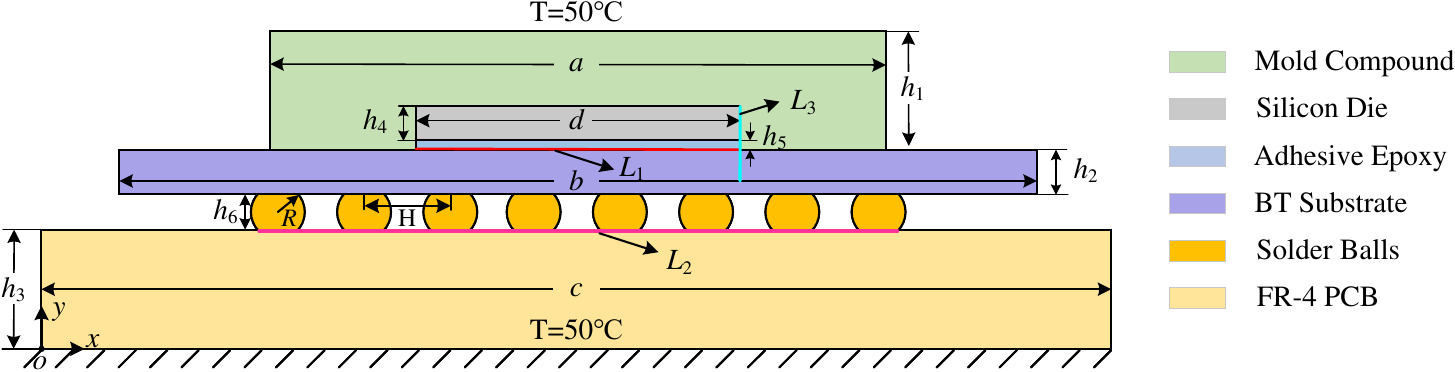}
    \caption{Geometric configuration and component layout of FC-BGA model.}
    \label{fig:fc_bga_model}
\end{figure}
\begin{table}[htbp]
\centering
\scriptsize
    \caption{Geometric dimensions of FC-BGA module components (mm).}
    \label{tab:fc_bga_dimensions}
    \begin{tabular}{lcc}
    \toprule
    Component & Length& Thickness \\
    \midrule
    Mold Compound& $a = 9.00$  &  $h_1 = 1.20$  \\
    BT Substrate& $b = 11.00$  & $h_2 = 0.40$  \\
    FR-4 PCB& $c = 13.50$  & $h_3 = 0.80$\\
    Silicon Die& $d = 5.00$  & $h_4 = 0.30$  \\
    Adhesive Epoxy& $d = 5.00$  & $h_5 = 0.10$  \\
    Solder Balls& $R = 0.76$  & $h_6 = 0.56$  \\
    \bottomrule
    \end{tabular}
\end{table}

The fourth numerical example examines the thermomechanical performance of a Flip-Chip Ball Grid Array (FC-BGA) packaging structure under thermal loading conditions. FC-BGA packaging technology offers advantages of high integration density and superior performance with excellent product reliability and stability, making it widely adopted in high-performance processors and advanced semiconductor devices. As illustrated in Fig.~\ref{fig:fc_bga_model}, the computational model consists of six distinct components: mold compound, silicon die, adhesive epoxy, BT substrate, solder balls, and FR-4 PCB. The geometric dimensions of each component are detailed in Table~\ref{tab:fc_bga_dimensions}. The solder balls are uniformly distributed with a pitch of H=1.24 mm.

\begin{table}[htbp]
    \scriptsize
    \centering
    \caption{Material properties of FC-BGA model components.}
    \label{tab:fc_bga_materials}
    \begin{tabular}{lcccc}
    \toprule
    Material & \begin{tabular}{c} Young's modulus \\ (GPa) \end{tabular} & Poisson's ratio & \begin{tabular}{c} Thermal conductivity \\ (W/m$\cdot$K) \end{tabular} & \begin{tabular}{c} CTE $\alpha$ \\ ( $^\circ$C$^{-1}$) \end{tabular} \\
    \midrule
    Mold Compound & 24 & 0.25 & 2.1 & 10 $\times 10^{-6}$\\
    Silicon Die & 165.5 & 0.25 & 119 & 2.8 $\times 10^{-6}$\\
    Solder Balls & 11 & 0.11 & 73 & 35 $\times 10^{-6}$\\
    Adhesive Epoxy & 2.6 & 0.3 & 0.188 & 90 $\times 10^{-6}$\\
    BT Substrate & 26 & 0.19 & 14.5 & 14 $\times 10^{-6}$\\
    FR-4 PCB & 22 & 0.28 & 6.5 & 18 $\times 10^{-6}$\\
    \bottomrule
    \end{tabular}
\end{table}

The material properties detailed in Table~\ref{tab:fc_bga_materials} show significant variations in thermal and mechanical characteristics that contribute to complex thermomechanical interactions. The boundary conditions are defined as follows: the bottom surface is fully constrained and other boundaries of the model are set as free boundaries, the silicon die temperature is set to $T = 500~^\circ$C, the upper surface of the mold compound and lower surface of the FR-4 PCB are maintained at the given temperature $T = 50~^\circ$C, and all remaining surfaces are treated as thermally insulated. This configuration creates a severe thermal gradient across the package with substantial thermal expansion mismatches ranging from 2.8 $\times$ 10$^{-6}$ $^\circ$C$^{-1}$ for silicon to 90 $\times$ 10$^{-6}$ $^\circ$C$^{-1}$ for adhesive epoxy, generating significant thermomechanical stresses that make this an excellent test case for evaluating the FE-VE coupling methodology in complex multi-material electronic packaging applications where reliability assessment is critical.

\begin{figure}[htbp]
    \centering
    \includegraphics[width=0.9\textwidth]{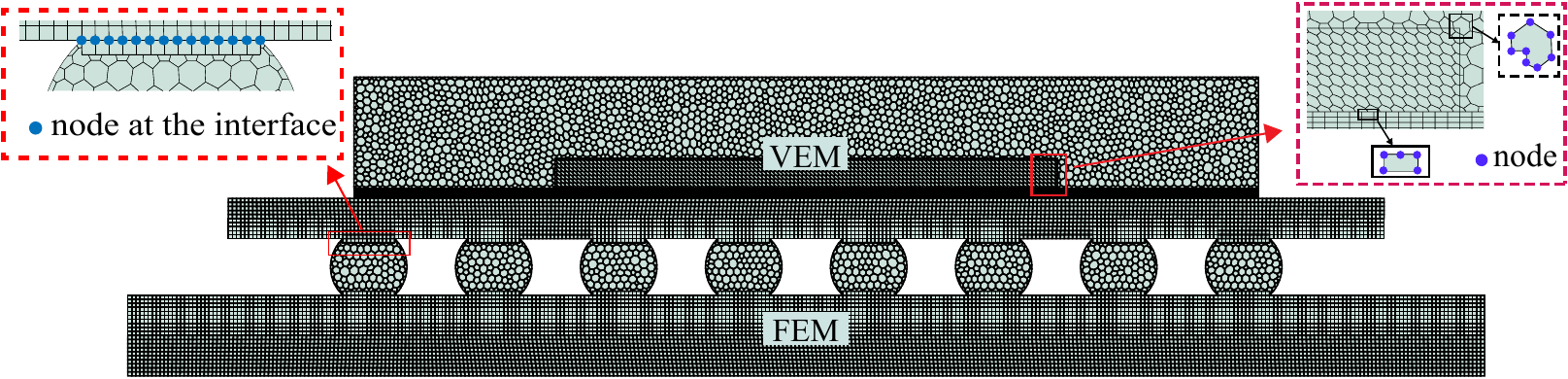}
    \caption{Hybrid mesh discretization strategy for FC-BGA model.}
    \label{fig:fc_bga_mesh}
\end{figure}

In this numerical example, the regular geometrical regions consisting of the BT substrate and FR-4 PCB are discretized using the finite element domain, while the remaining components are assigned to the virtual element domain. The finite element domain employs four-node quadrilateral elements for discretization, while the virtual element domain utilizes regular polygonal elements at interface boundaries and arbitrary polygonal elements in non-interface regions, with non-matching mesh techniques used to integrate the VEM domain meshes into a unified structure. 
To ensure nodal compatibility at the coupling interface, careful attention is given to maintaining nodal correspondence between domains. Fig.~\ref{fig:fc_bga_mesh} illustrates the final mesh discretization scheme, demonstrating the hybrid meshing strategy where complex geometries of the upper package components (mold compound, silicon die, adhesive epoxy) and solder balls are effectively handled by the flexible VEM approach, while structured lower layers benefit from the computational efficiency of conventional FEM. 
The interface region clearly shows how nodal correspondence between VE and FE domains is maintained to ensure solution continuity across the coupling boundary, demonstrating the seamless integration of both numerical methods in this complex multi-material electronic packaging configuration.

\begin{figure}[htbp]
    \centering
    \subfloat[FE-VE coupling solution]{\centering
        \includegraphics[width=0.46\textwidth]{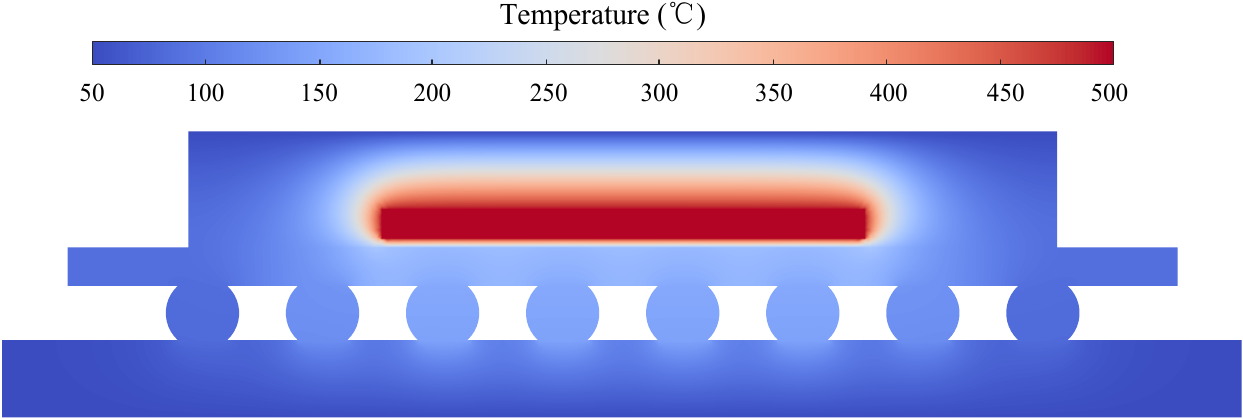}
        \label{fig:fcbga_temp_coupling}
    }
    \hfill
    \subfloat[FEM reference solution]{\centering
        \includegraphics[width=0.46\textwidth]{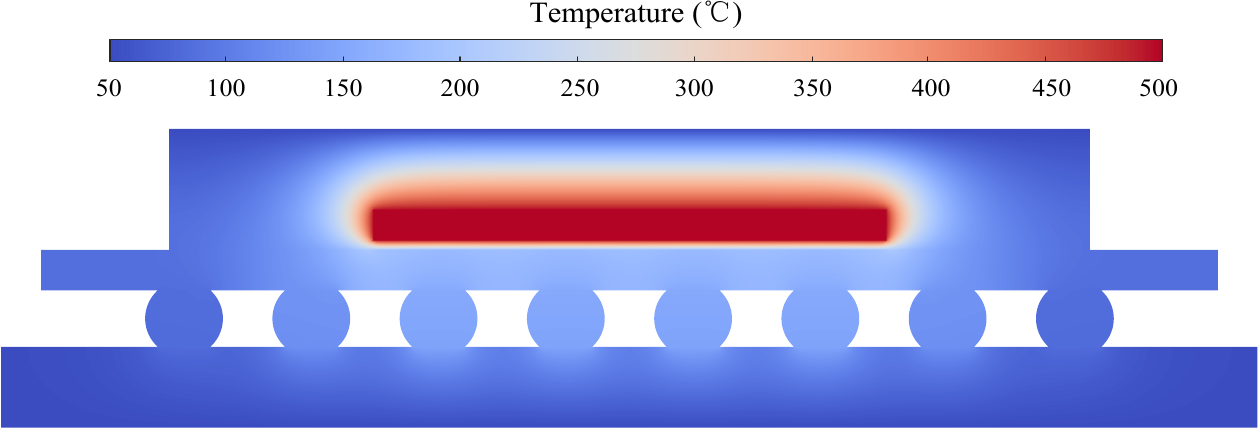}
        \label{fig:fcbga_temp_fem}
    }
    \caption{Temperature distribution comparison for FC-BGA.}
    \label{fig:fcbga_temp_comparison}
\end{figure}
 
Fig.~\ref{fig:fcbga_temp_comparison} presents the temperature field distributions computed by the two methods. The FEM solution with mesh refinement (nDof = 162,444) is shown for comparison. Both methods demonstrate consistent temperature distributions ranging from 50 $^\circ$C to 500 $^\circ$C.
 
\begin{figure}[htbp]
    \centering
    \subfloat[Horizontal displacement ($u_x$)]{\centering
        \includegraphics[width=0.33\textwidth]{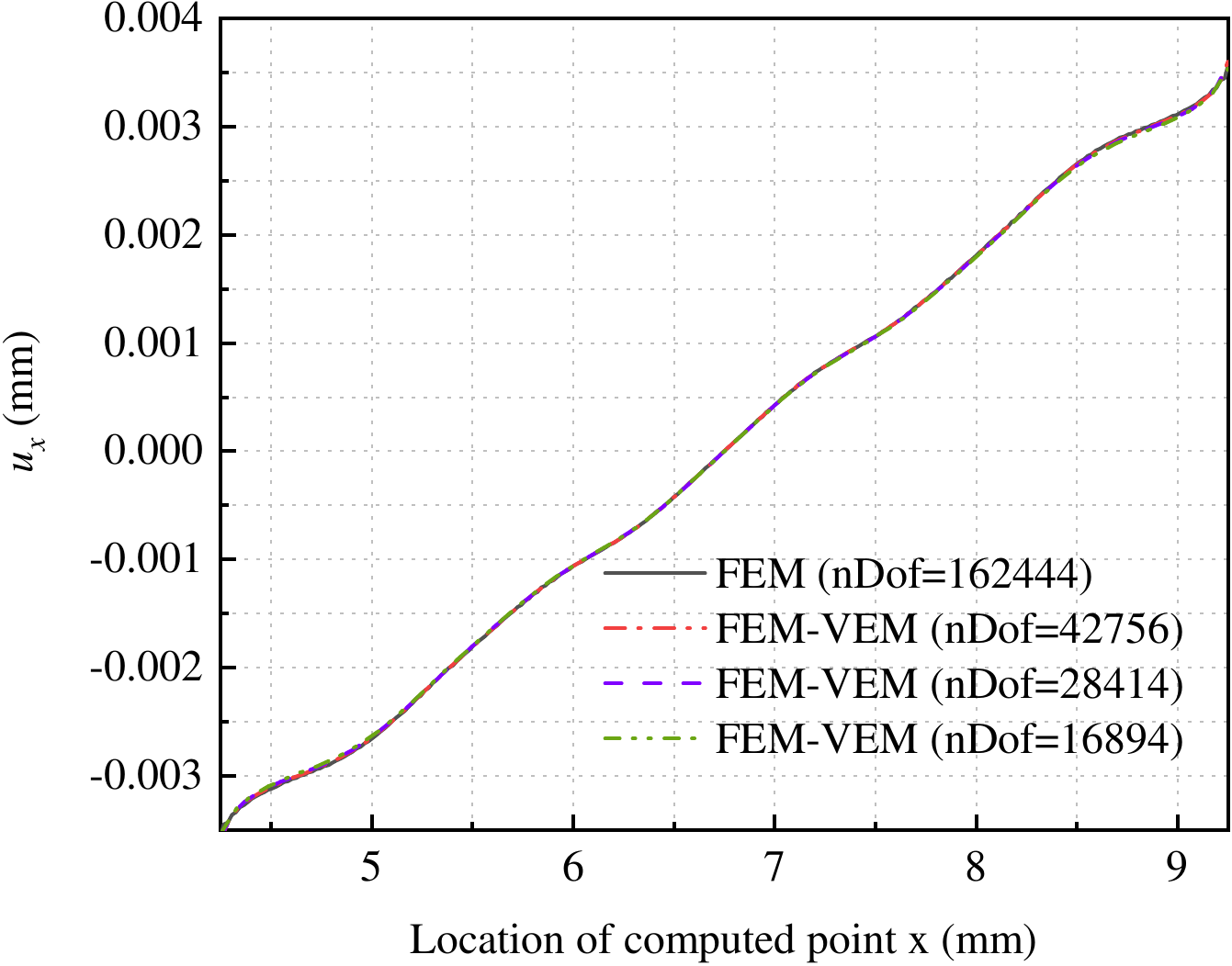}
        \label{fig:fc_bga_ux_l1}
    }
    \subfloat[Vertical displacement ($u_y$)]{\centering
        \includegraphics[width=0.33\textwidth]{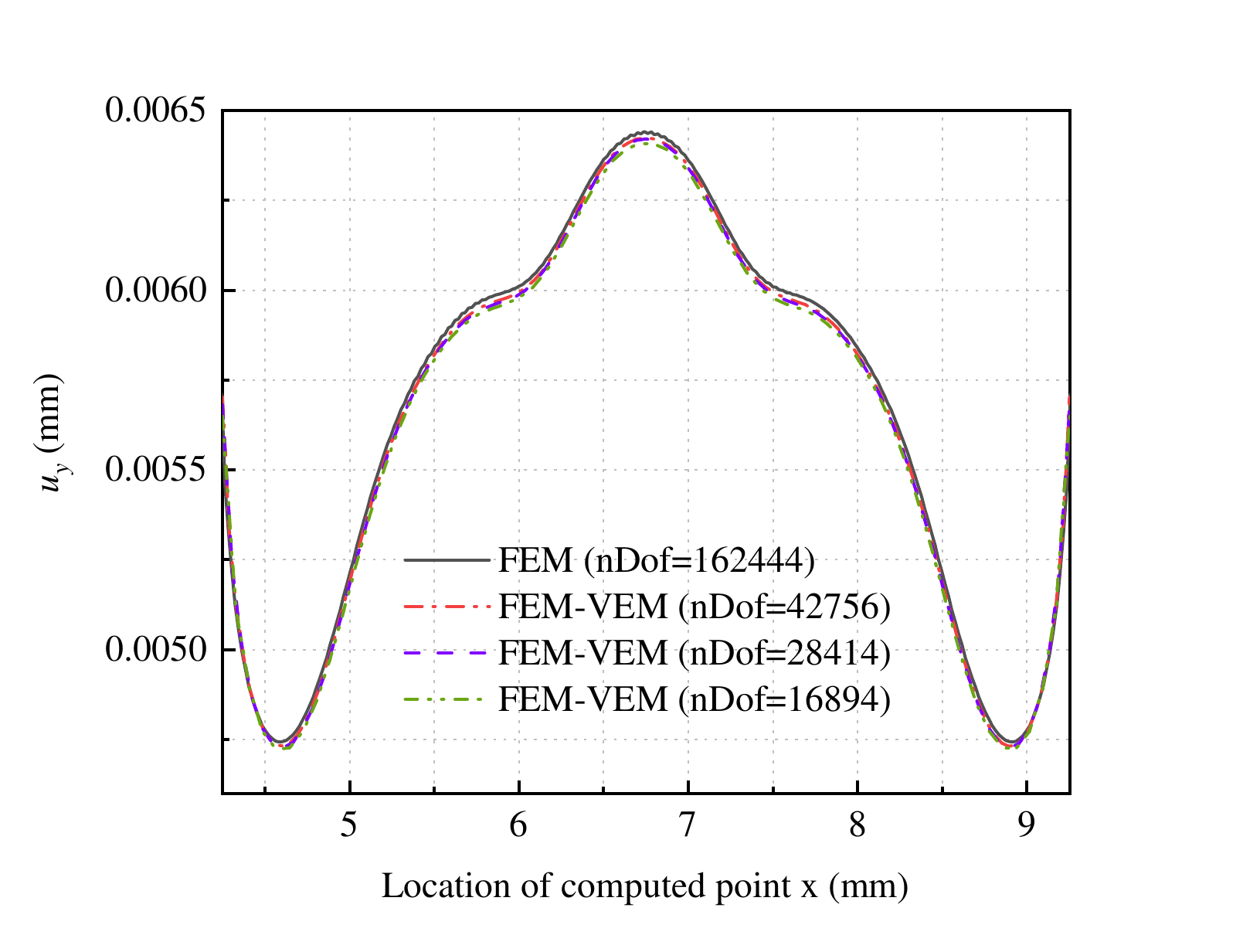}
        \label{fig:fc_bga_uy_l1}
    }
    \subfloat[von Mises stress distribution]{\centering
        \includegraphics[width=0.31\textwidth]{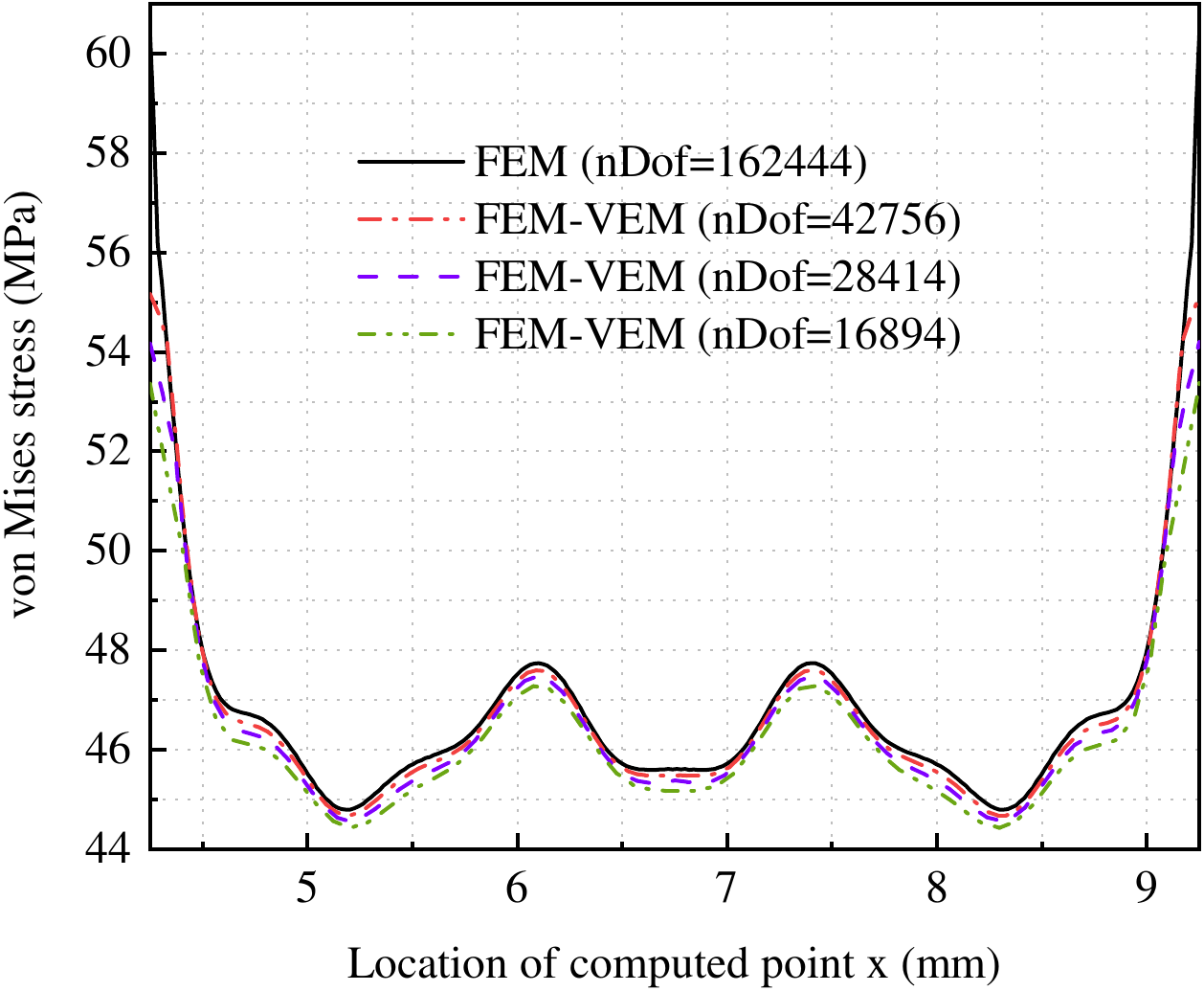}
        \label{fig:fc_bga_stress_l1}
    }
    \caption{Displacement and stress distributions along interface line $L_1$.}
    \label{fig:fc_bga_results_l1}
\end{figure}

\begin{figure}[htbp]
    \centering
    \includegraphics[width=0.6\textwidth]{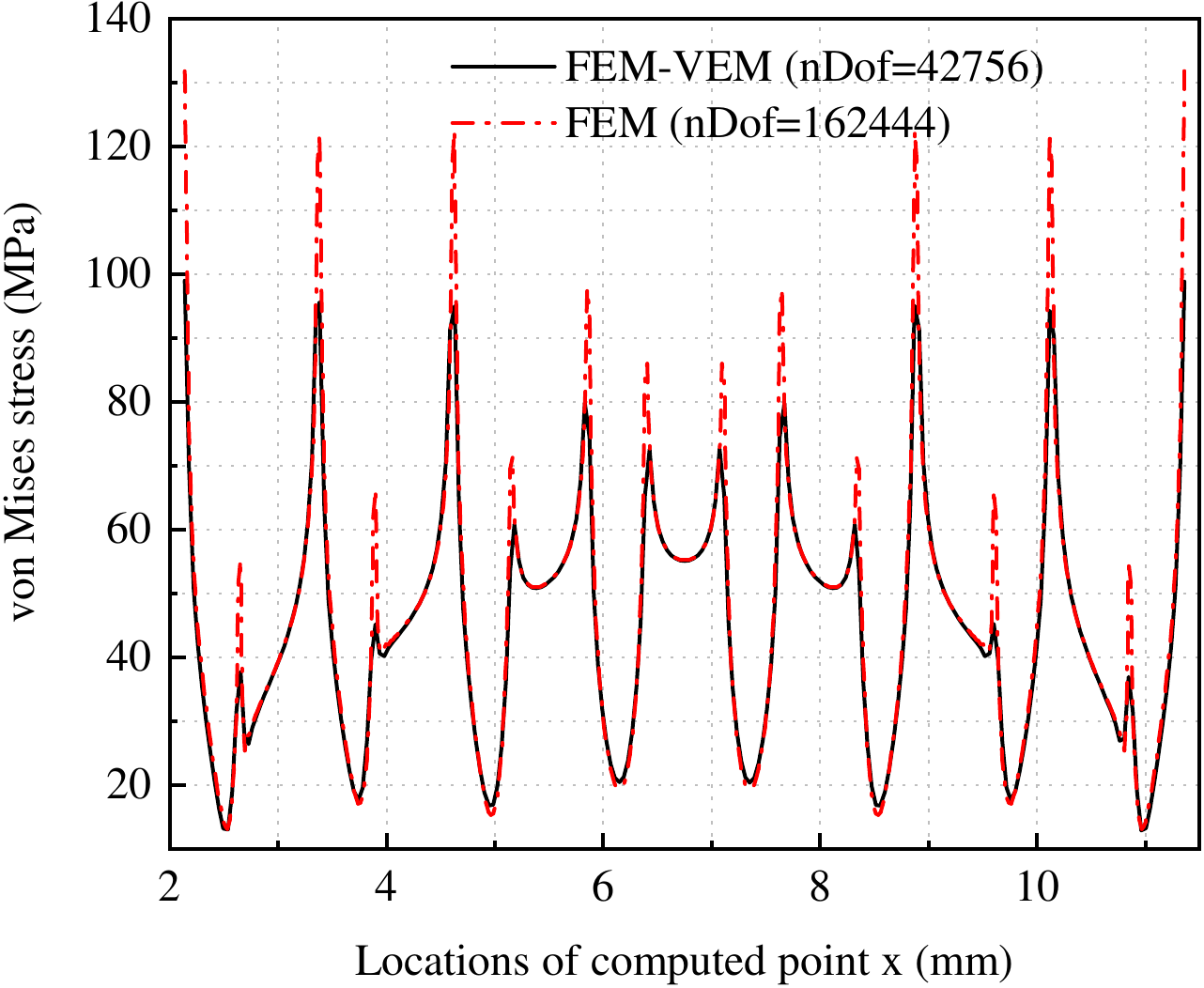}
    \caption{von Mises stress distribution along line $L_2$.}
    \label{fig:fc_bga_results_l2}
\end{figure}

To further validate the computational accuracy of the FE-VE coupling algorithm, a refined finite element analysis with mesh densification (nDof = 162,444) was performed as the reference solution. The displacement and stress profiles were extracted along two critical lines $L_1$ and $L_2$ (given in Fig.~\ref{fig:fc_bga_model}), as presented in Figs.~\ref{fig:fc_bga_results_l1} and~\ref{fig:fc_bga_results_l2}. Fig.~\ref{fig:fc_bga_results_l1} demonstrates excellent agreement between the coupling algorithm and the reference FEM solution across all three response quantities along interface $L_1$. Fig.~\ref{fig:fc_bga_ux_l1} shows the displacement $u_x$ exhibiting a nearly linear variation, indicating consistent thermal expansion behavior with virtually identical results across all mesh refinement levels. Fig.~\ref{fig:fc_bga_uy_l1} presents a bell-shaped distribution of $u_y$, with the most significant thermal expansion effect at the center due to the temperature gradient. Fig.~\ref{fig:fc_bga_stress_l1} illustrates the stress distribution with periodic variations, showing peak stresses occurring at the interface boundaries and lower stresses in the central region.

\begin{figure}[htbp]
    \centering
   \begin{subfigure}[b]{0.9\textwidth}
     \centering
        \includegraphics[width=0.9\textwidth]{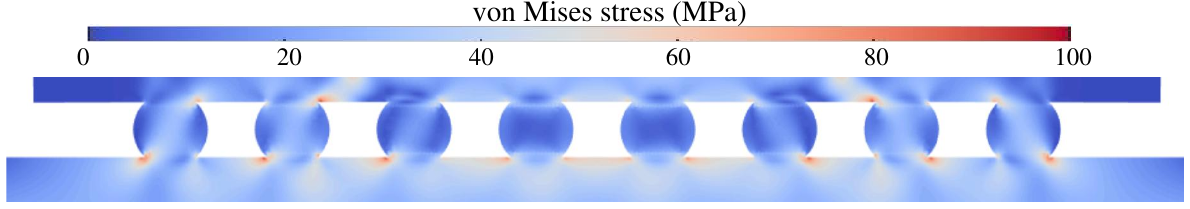}
            \caption{FE-VE coupling solution in nDof=42,756}
        \label{fig:fcbga_mises_coupling}
     \end{subfigure}
    \hfill
      \begin{subfigure}[b]{0.9\textwidth}
     \centering
        \includegraphics[width=0.9\textwidth]{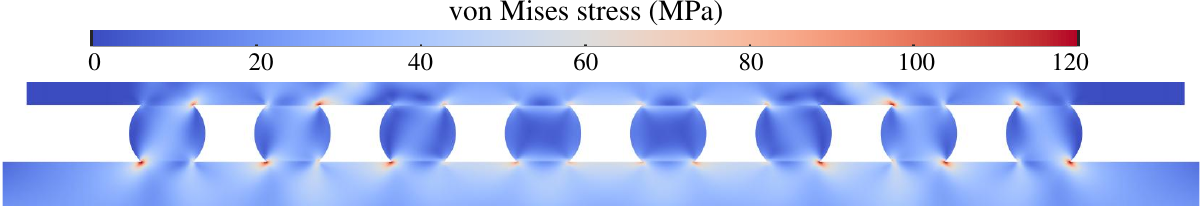}
         \caption{FEM reference solution in nDof=162,444}
        \label{fig:fcbga_mises_fem}
        \end{subfigure}
    
    \caption{Temperature distribution comparison for FC-BGA.}
    \label{fig:fcbga_mises_comparison}
\end{figure}

Fig.~\ref{fig:fc_bga_results_l2} presents the von Mises stress distribution along line $L_2$ (shown in Fig.~\ref{fig:fc_bga_model}), which passes through the solder ball array. The stress profile exhibits pronounced periodic variations corresponding to the solder ball locations.
Fig.~\ref{fig:fcbga_mises_comparison} displays the von Mises stress distributions around the solder joints obtained using both methods. The comparison between FE-VE coupling (nDof = 42,756) and reference FEM solution (nDof = 162,444) shows reasonable agreement in capturing the periodic stress pattern, though some differences are observed at stress concentration locations due to the coarser discretization in the coupling method. Despite the lower mesh density, the FE-VE approach successfully captures the overall stress distribution characteristics and solder ball interface behavior.

These comprehensive results demonstrate that the FE-VE coupling algorithm maintains computational accuracy in both regular and complex geometrical regions, effectively combining the computational efficiency of finite elements with the geometric flexibility of virtual elements for multi-scale electronic packaging applications.

\begin{figure}[htbp]
    \centering
    \includegraphics[width=0.5\textwidth]{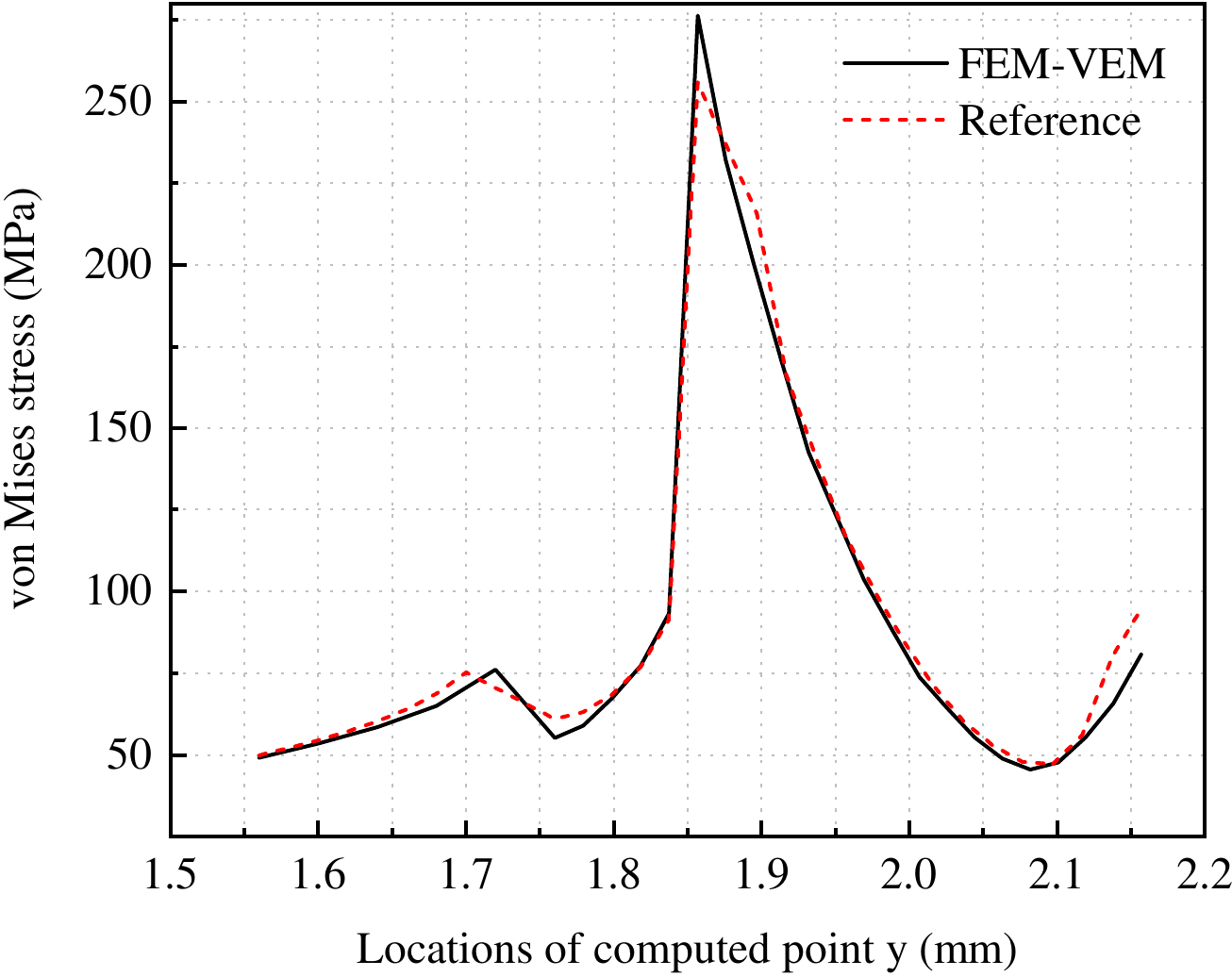}
    \caption{von Mises stress distribution along cross-material line $L_3$.}
    \label{fig:fc_bga_stress_l3}
\end{figure}

Fig.~\ref{fig:fc_bga_stress_l3} presents the nodal stress distribution along cross-material line $L_3$ (the blue line in Fig.~\ref{fig:fc_bga_model}), which traverses multiple material interfaces. The FE-VE coupling algorithm shows close agreement with the reference finite element solution across different material regions. 
The stress distribution 
demonstrating that the coupling algorithm maintains solution continuity when transitioning between materials with different thermal and mechanical properties. These results confirm the capability of the FE-VE methodology for multi-material thermomechanical problems in electronic packaging applications.

\subsection{IGBT Module}
\label{sec:igbt_module}

\begin{figure}[htbp]
    \centering
    \includegraphics[width=0.8\textwidth]{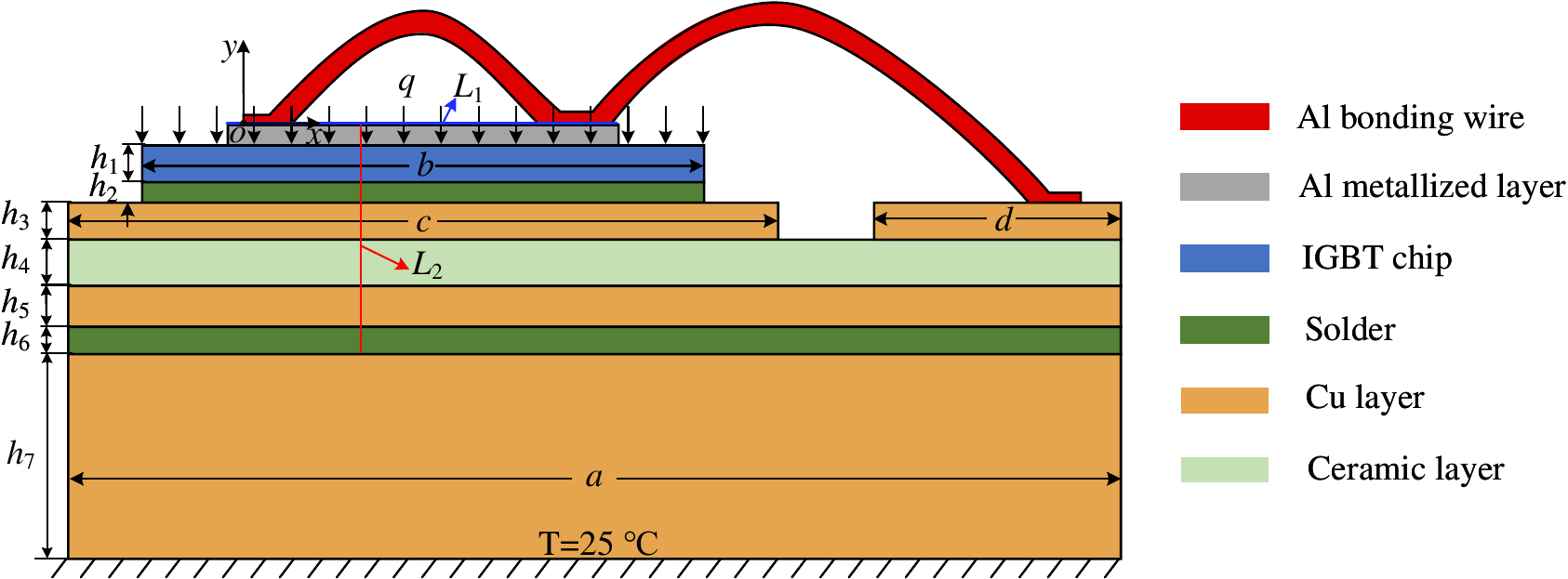}
    \caption{Geometric configuration and layer structure of IGBT module 2D model.}
    \label{fig:igbt_model}
\end{figure}

\begin{figure}[htbp]
    \centering
    \includegraphics[width=0.7\textwidth]{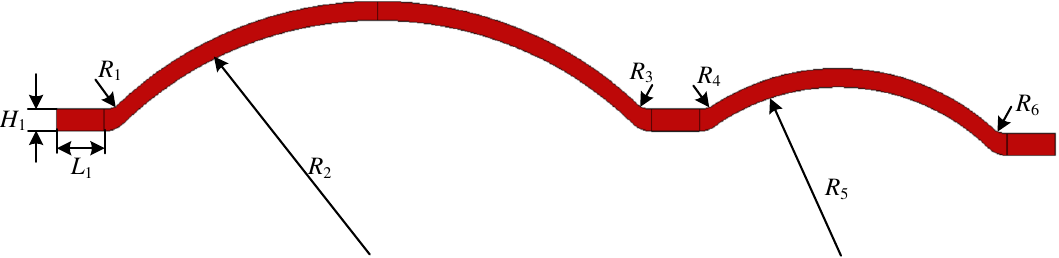}
    \caption{Geometric configuration and dimensional parameters of Al bonding wire.}
    \label{fig:al_bonding_wire}
\end{figure}

The final numerical example investigates the thermomechanical coupling behavior of an Insulated Gate Bipolar Transistor (IGBT) module, which represents a critical component in power electronics applications. The multi-layer IGBT structure experiences significant thermal stresses due to material property mismatches, making thermal stress analysis essential for reliability assessment.
Fig.~\ref{fig:igbt_model} illustrates the structural configuration of the IGBT model, comprising nine distinct layers from top to bottom: Al bonding wires, Al metallized layer, IGBT chip, chip solder layer, upper Cu layers, ceramic layer, lower Cu layer, substrate solder layer, and Cu baseplate. The geometric dimensions of each component are detailed in Table~\ref{tab:igbt_dimensions}, while the Al bonding wire geometry is shown in Fig.~\ref{fig:al_bonding_wire} with specific dimensions listed in Table~\ref{tab:al_wire_dimensions}. This complex multi-layer structure with curved bonding wire geometries presents an ideal test case for the FE-VE coupling methodology, where VEM's geometric flexibility effectively handles complex wire profiles while FEM efficiently manages regular substrate layers.

\begin{table}[htbp]
    \scriptsize
    \centering
\caption{Geometric dimensions of IGBT module components (mm).}
    \label{tab:igbt_dimensions}
    \begin{tabular}{lcc}
    \toprule
    Component & length& Thickness\\
    \midrule
    IGBT surface Al metallized layer &11.0 & 0.004  \\
    IGBT chip & $b = 13.0$  & $h_1 = 0.20$  \\
    IGBT chip solder layer & $b = 13.0$  & $h_2 = 0.15$  \\
    Upper Cu layer 1 & $c = 15.0$  & $h_3 = 0.30$  \\
    Upper Cu layer 2 & $d = 18.0$  & $h_3 = 0.30$  \\
    Ceramic layer & $a = 18.0$  & $h_4 = 0.38$  \\
    Lower Cu layer & $a = 18.0$  & $h_5 = 0.30$  \\
    Substrate solder layer & $a = 18.0$  & $h_6 = 0.15$  \\
    Cu baseplate & $a = 18.0$  & $h_7 = 3.00$  \\
    \bottomrule
    \end{tabular}
\end{table}

\begin{table}[htbp]
    \scriptsize
    \centering
    \caption{Dimensional parameters of Al bonding wire (mm).}
    \label{tab:al_wire_dimensions}
    \begin{tabular}{ccccccccc}
    \toprule
    Parameter & $H_1$ & $L_1$ & $R_1$ & $R_2$ & $R_3$ & $R_4$ & $R_5$ & $R_6$ \\
    \midrule
    Dimension & 0.32 & 0.70& 0.40& 5.36 & 0.40& 0.40& 3.16 & 0.40\\
    \bottomrule
    \end{tabular}
\end{table}

\begin{table}[htbp]
    \scriptsize
    \centering
    \caption{Material properties of IGBT module components.}
    \label{tab:igbt_materials}
    \begin{tabular}{lcccc}
    \toprule
    Material & \begin{tabular}{c} Thermal conductivity \\ (W/m$\cdot$K) \end{tabular} & \begin{tabular}{c} Young's modulus \\ (GPa) \end{tabular} & Poisson's ratio & \begin{tabular}{c} CTE $\alpha$\\ ($^\circ$C$^{-1}$) \end{tabular} \\
    \midrule
    Al & 237 & 70.6 & 0.33 & 21.0 $\times 10^{-6}$ \\
    Cu & 400 & 100 & 0.34 & 16.4 $\times 10^{-6}$ \\
    Al$_2$O$_3$ & 20 & 300 & 0.22 & 6.4 $\times 10^{-6}$ \\
    IGBT Chip& 148 & 112 & 0.22 & 2.5 $\times 10^{-6}$ \\
    Solder & 57 & 10.6 & 0.35 & 22.4 $\times 10^{-6}$\\
    \bottomrule
    \end{tabular}
\end{table}

Material parameters are assigned to each component as listed in Table~\ref{tab:igbt_materials}, which reveals significant variations in thermal and mechanical characteristics that contribute to complex thermomechanical interactions. The boundary conditions simulate realistic operating conditions: a heat flux of $q = 1000$ mW/mm$^2$ is applied to the upper surface of the IGBT chip to represent power dissipation during operation, the lower surface of the Cu baseplate is maintained at $T = 25^\circ$C representing heat sink cooling, and all remaining surfaces are treated as thermally insulated with $\nabla T \cdot \mathbf{n} = 0$. The bottom surface is subjected to a fully fixed mechanical constraint (ux = uy = 0) and other boundaries of the model are set as free boundaries.
This configuration creates a realistic thermal gradient from the heat-generating chip through the multi-layer package structure to the cooled baseplate. The substantial differences in thermal expansion coefficients
generate significant thermal stresses that make this an excellent validation case for the FE-VE coupling methodology in handling complex multi-material electronic packaging with curved geometries and high thermal gradients.


\begin{figure}[htbp]
    \centering
    \includegraphics[width=0.8\textwidth]{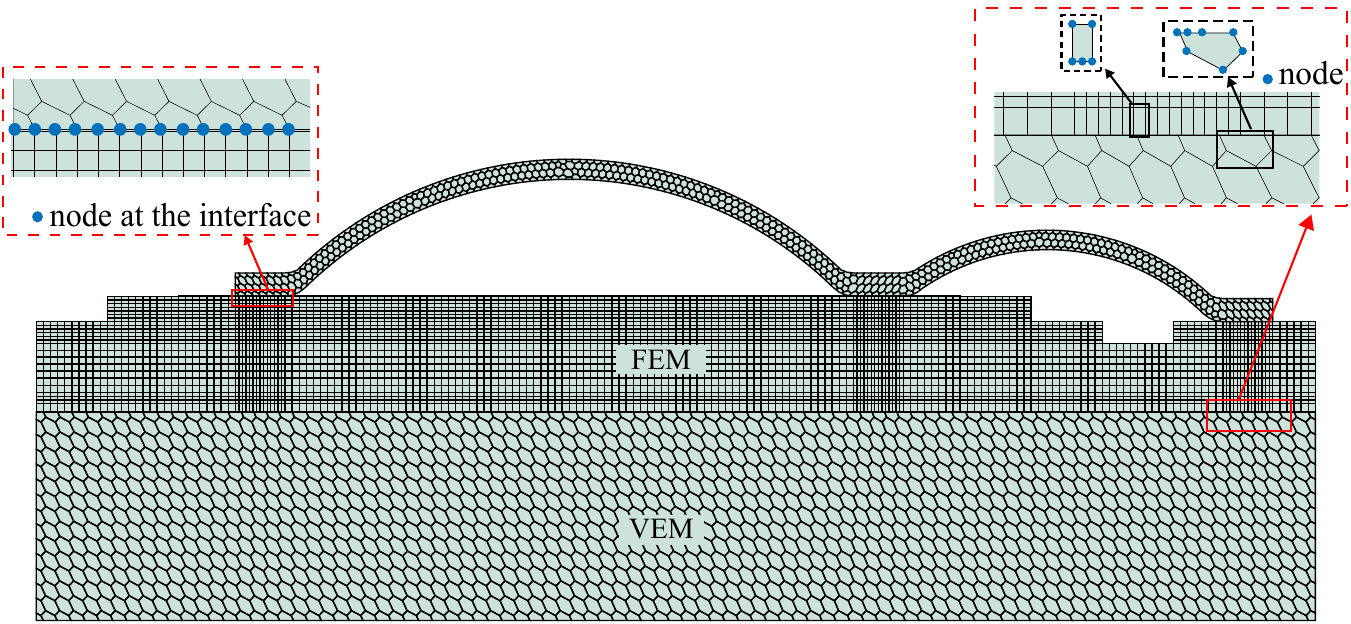}
    \caption{Hybrid mesh discretization strategy for IGBT module.}
    \label{fig:igbt_mesh}
\end{figure}

Here, the complex Al bonding wire component is assigned to the virtual element domain, where polygonal elements facilitate mesh generation for intricate curved geometries, while the Cu baseplate is also discretized using the virtual element domain to reduce element count and computational cost. The remaining components are assigned to the finite element domain, employing four-node quadrilateral elements with careful attention to ensuring shared nodes at coupling interfaces to maintain solution continuity. Fig.~\ref{fig:igbt_mesh} illustrates the hybrid mesh discretization strategy, demonstrating how the approach effectively combines VEM's geometric flexibility for complex wire profiles with FEM's computational efficiency for regular substrate layers. The interface regions show where nodal compatibility between VEM and FEM domains is maintained, enabling seamless integration of both numerical methods in this challenging multi-material electronic packaging configuration.

\begin{figure}[htbp]
    \centering
    \subfloat[FE-VE coupling solution]{\centering
        \includegraphics[width=0.46\textwidth]{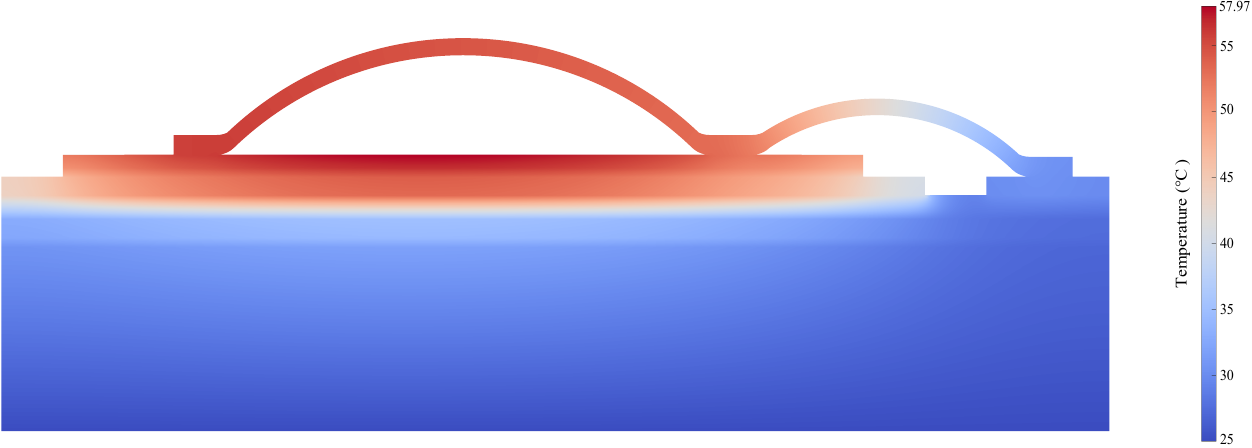}
        \label{fig:igbt_temp_coupling}
    }
    \hfill
    \subfloat[FEM reference solution]{\centering
        \includegraphics[width=0.46\textwidth]{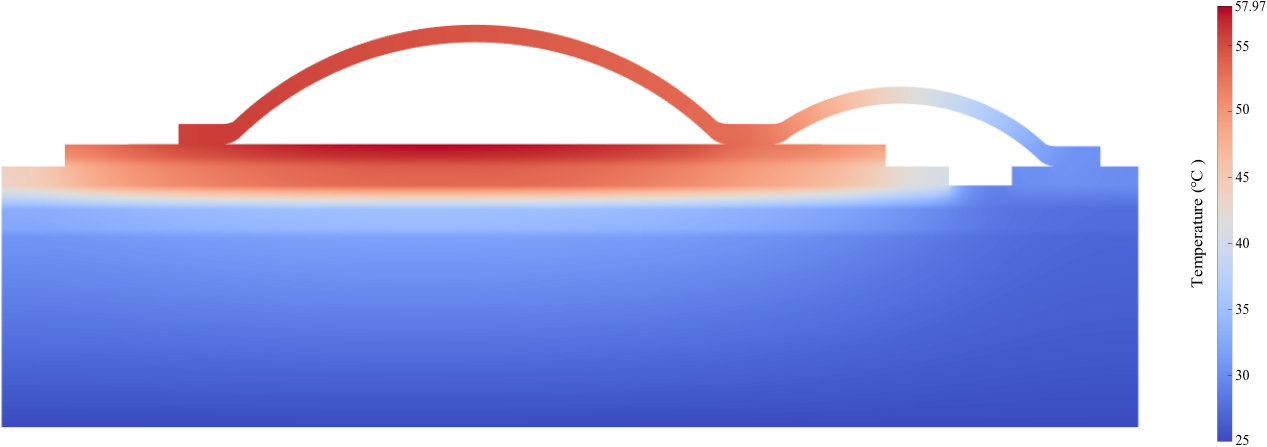}
        \label{fig:igbt_temp_fem}
    }
    \caption{Temperature distribution comparison for IGBT module.}
    \label{fig:igbt_temp_comparison}
\end{figure}

Fig.~\ref{fig:igbt_temp_comparison} displays the temperature field distributions under heat flux. Both methods show excellent agreement, with temperatures ranging from 25 $^\circ$C to 57.97 $^\circ$C.

\begin{figure}[htbp]
    \centering
    \subfloat[FE-VE coupling solution]{\centering
        \includegraphics[width=0.46\textwidth]{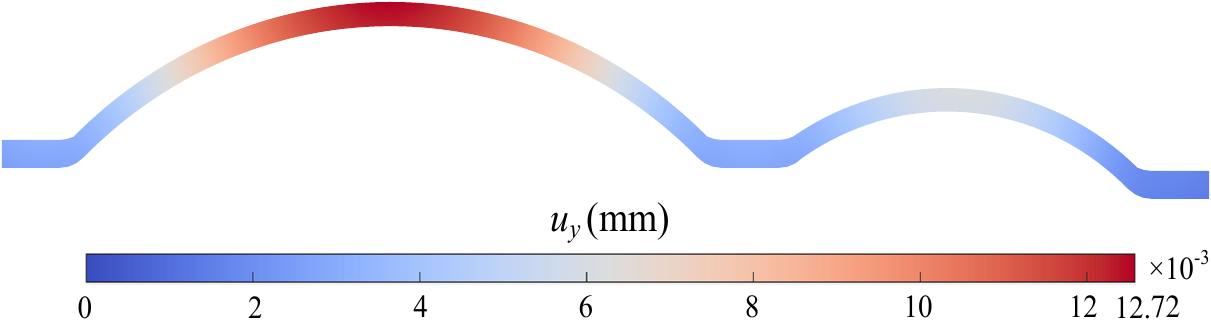}
        \label{fig:igbt_uy_coupling}
    }
    \hfill
    \subfloat[FEM reference solution]{\centering
        \includegraphics[width=0.46\textwidth]{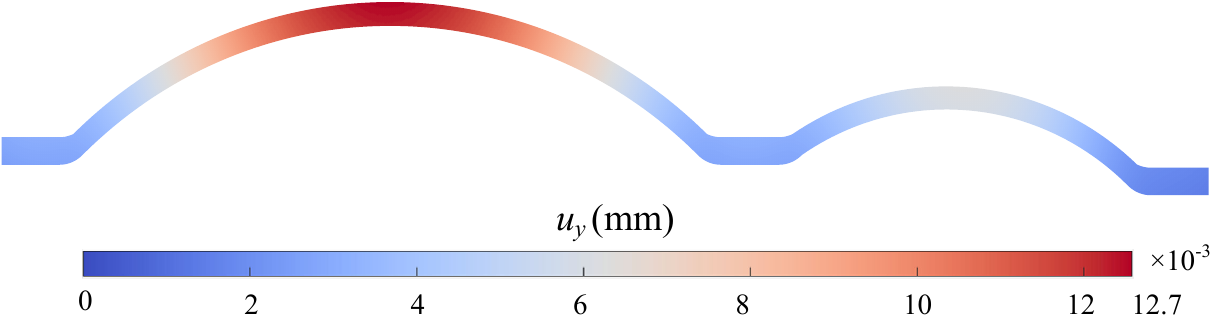}
        \label{fig:igbt_uy_fem}
    }
    \caption{Vertical displacement ($u_y$) distribution comparison for Al bonding wire.}
    \label{fig:igbt_uy_comparison}
\end{figure}

Fig.~\ref{fig:igbt_uy_comparison} presents the vertical displacement distribution contours for the Al bonding wire, showing close agreement between both methods. The displacement values range from 0 to approximately 0.01272 mm for the coupling method and 0 to 0.0127 mm for the FEM method, demonstrating consistent displacement patterns throughout the complex wire geometry due to thermal expansion effects.

\begin{figure}[htbp]
    \centering
    \subfloat[Horizontal displacement ($u_x$)]{\centering
        \includegraphics[width=0.31\textwidth]{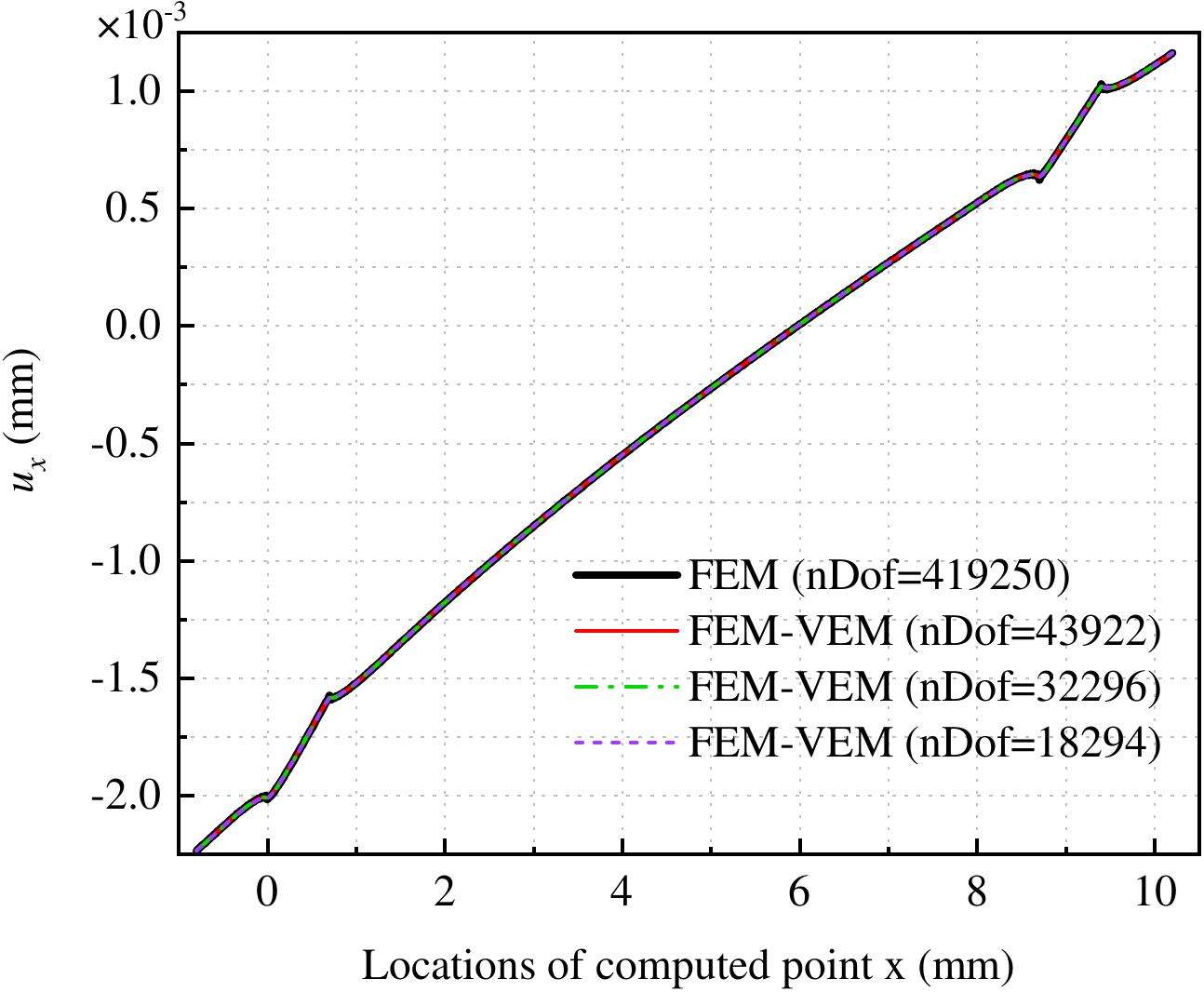}
        \label{fig:igbt_ux_l1}
    }
    \subfloat[Vertical displacement ($u_y$)]{\centering
        \includegraphics[width=0.31\textwidth]{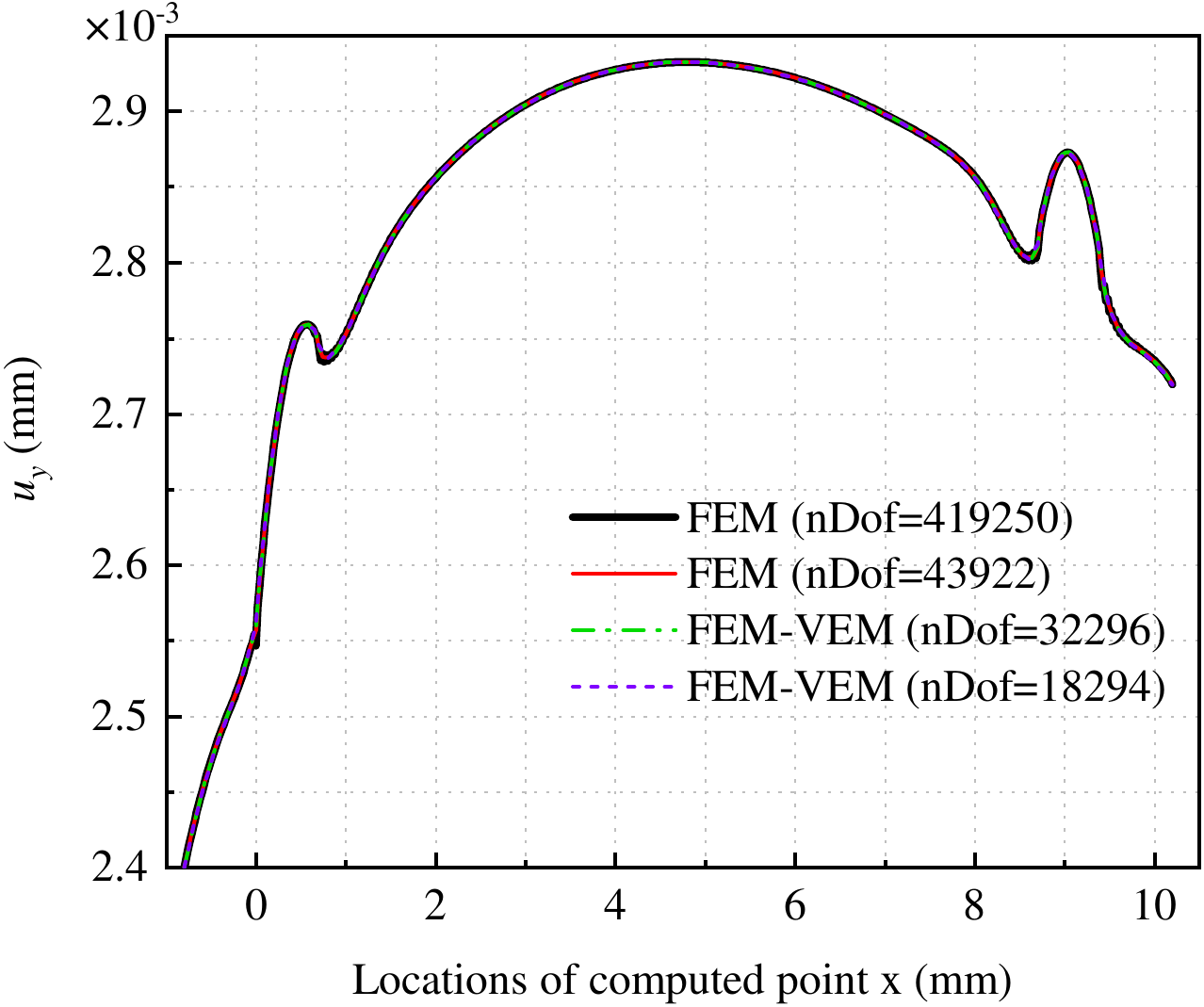}
        \label{fig:igbt_uy_l1}
    }
    \subfloat[von Mises stress distribution]{\centering
        \includegraphics[width=0.31\textwidth]{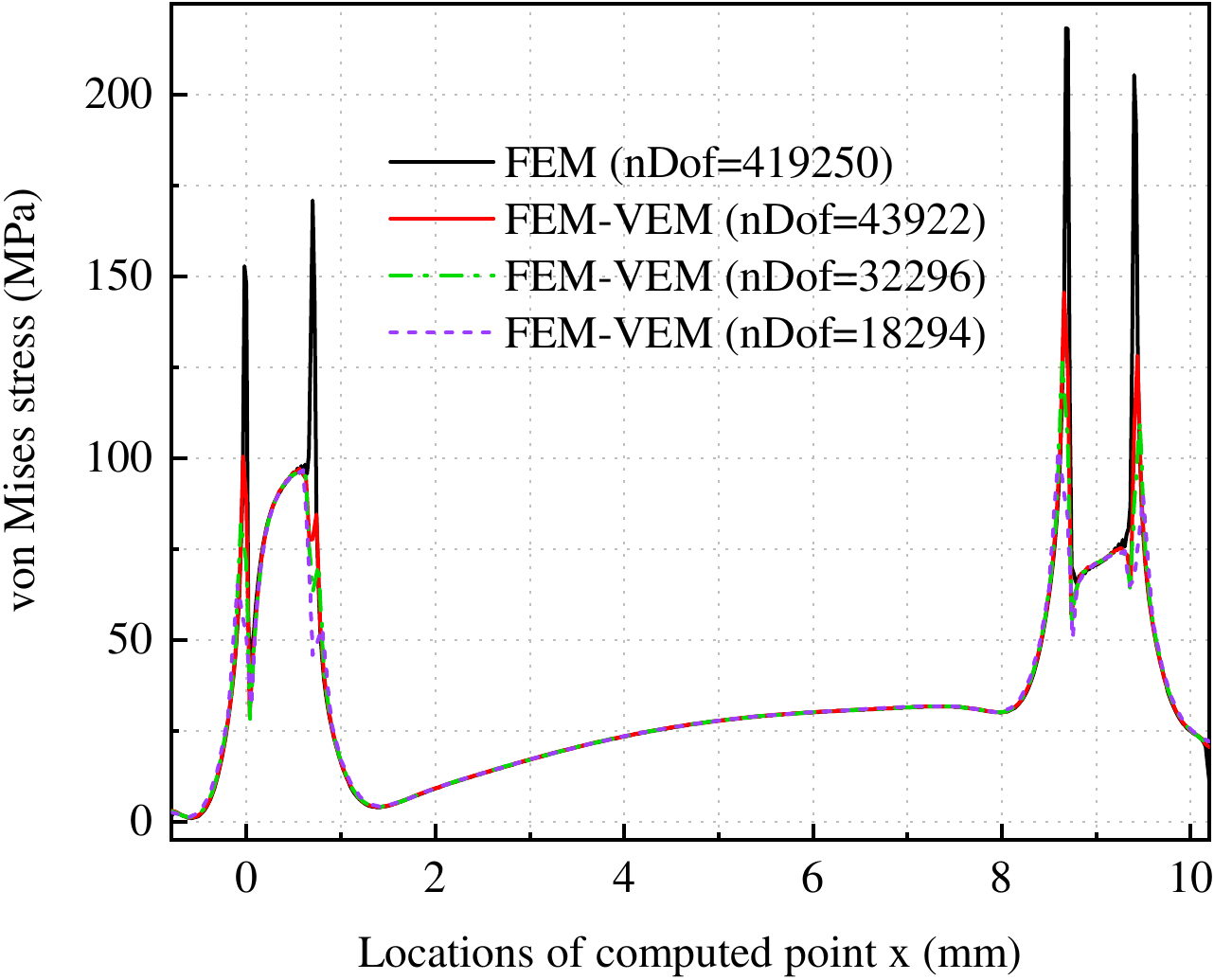}
        \label{fig:igbt_stress_l1}
    }
    \caption{Displacement and stress distributions along interface line $L_1$ (shown in Fig.~\ref{fig:igbt_model}).}
    \label{fig:igbt_results_l1}
\end{figure}

\begin{figure}[htbp]
    \centering
    \subfloat[von Mises stress distribution]{\centering
        \includegraphics[width=0.46\textwidth]{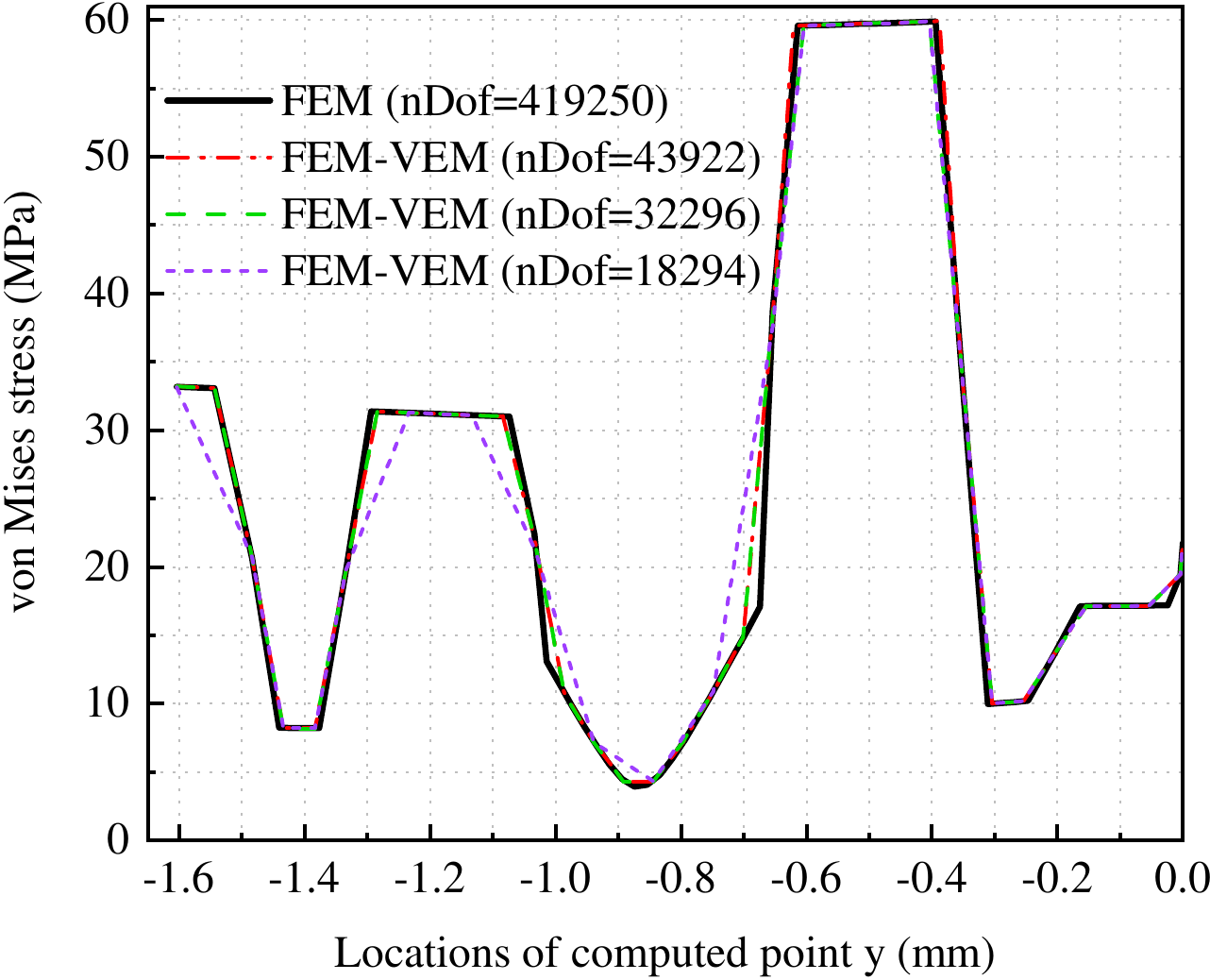}
        \label{fig:igbt_stress_l2}
    }
    \hfill
    \subfloat[Vertical displacement ($u_y$)]{\centering
        \includegraphics[width=0.46\textwidth]{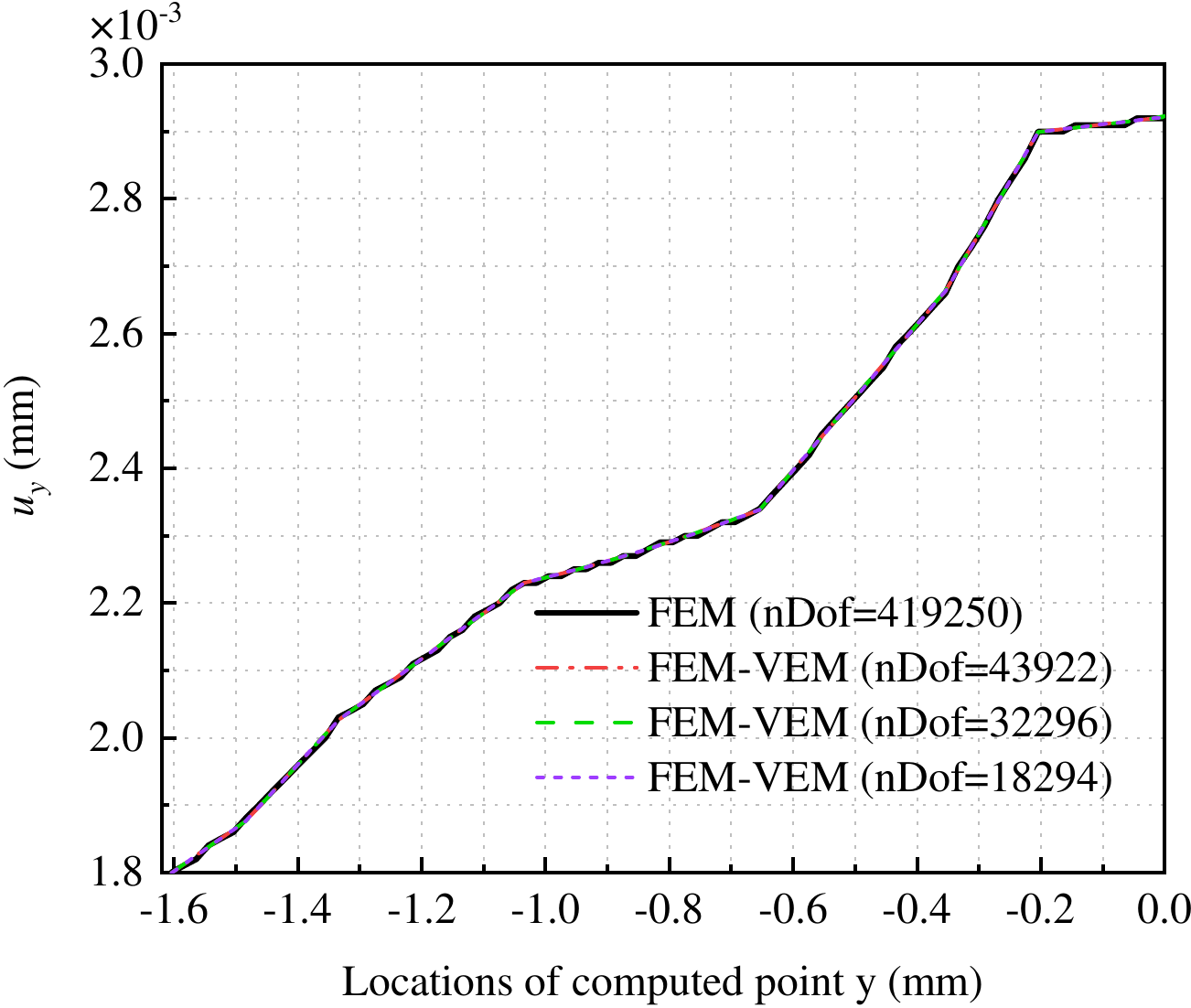}
        \label{fig:igbt_uy_l2}
    }
    \caption{Stress and displacement distributions along cross-material line $L_2$.}
    \label{fig:igbt_results_l2}
\end{figure}

To provide detailed validation, nodal values were extracted along interface line $L_1$ (the blue line in Fig.~\ref{fig:igbt_model}) using refined finite element analysis (nDof = 419,250) as the reference solution, with results presented in Fig.~\ref{fig:igbt_results_l1}. The coupling algorithm shows good consistency across all response quantities.

Fig.~\ref{fig:igbt_results_l2} shows results along cross-material line $L_2$ (the red line in Fig.~\ref{fig:igbt_model}), where stress concentrations reach 60 MPa and vertical displacement varies from 1.8 $\times$ 10$^{-3}$ mm to 2.9 $\times$ 10$^{-3}$ mm. The coupling algorithm maintains consistency with the FEM reference solution across different mesh refinement levels, confirming the effectiveness of the proposed methodology for complex IGBT thermal stress analysis.

\section{Conclusions}
\label{sec:conclusions}
The developed FE-VE coupled method effectively captures thermal distributions and stress concentrations in multi-material systems while maintaining solution continuity across domain interfaces through strategic nodal correspondence. The approach successfully accommodates arbitrary polygonal elements in geometrically complex regions while leveraging computational efficiency in regular domains, demonstrating consistent performance across different mesh refinement levels. Validation through comprehensive numerical examples spanning sintered silver interconnects, FC-BGA packages, and IGBT modules demonstrates reasonable agreement with reference solutions and acceptable convergence characteristics. The coupled method maintains adequate computational precision while providing practical advantages through strategic domain partitioning, establishing a viable computational framework for thermomechanical analysis of electronic packaging structures with complex multi-scale geometries.

This FE-VE coupled framework advances computational mechanics by providing enhanced flexibility for analyzing electronic packaging technologies with complex geometric features. Future research directions include: (i) extending the approach to 3D thermomechanical problems, (ii) incorporating more sophisticated material models including plasticity and damage mechanics, (iii) investigating adaptive mesh refinement strategies for optimal domain partitioning, and (iv) developing fully coupled solutions for strongly bidirectional thermomechanical interaction problems. To facilitate reproducibility and further development, source codes and implementation examples are publicly available at https://github.com/yanpeng-gong/FeVeCoupled-ElectronicPackaging.

\section*{Acknowledgments}
This research was supported by the National Natural Science Foundation of China (No. 12002009). 








\bibliographystyle{elsarticle-num}  
\bibliography{FEMVEMcoupling_ref_lib.bib}  
\end{document}